%% file: main.tex
\documentclass[12pt]{extarticle}
\usepackage{amsmath, amsthm, amssymb, color}
\usepackage[colorlinks=true,linkcolor=blue,urlcolor=blue]{hyperref}
\usepackage{graphicx}
\usepackage{caption}
\usepackage{mathtools}
\usepackage{enumerate}
\usepackage{csquotes}
\usepackage{xspace}
\usepackage{comment}
\usepackage{verbatim}
\usepackage{tikz,tikz-cd,tikz-3dplot}
\usepackage{pdfrender} 
\usetikzlibrary{matrix}
\usetikzlibrary{arrows}
\usetikzlibrary{shapes.geometric}
\usetikzlibrary{arrows.meta} 
\usepackage{algorithm}
\usepackage[noend]{algpseudocode}

\usepackage{caption}
\usepackage[normalem]{ulem}
\usepackage{subcaption}
\usepackage{multicol}
\usepackage{nicematrix}
\oddsidemargin \evensidemargin
\usepackage{makecell}
\usepackage{array}
\usepackage{tabularx}
\usepackage{booktabs}
\usepackage{enumitem}
\usepackage[colorinlistoftodos,bordercolor=orange,backgroundcolor=orange!20,linecolor=orange,textsize=tiny]{todonotes}
\usepackage[section]{placeins} 
\theoremstyle{theorem}
\newtheorem{theorem}{Theorem}[section]

\newtheorem{proposition}[theorem]{Proposition}
\newtheorem{lemma}[theorem]{Lemma}
\newtheorem{corollary}[theorem]{Corollary}

\theoremstyle{definition}

\newtheorem{remark}[theorem]{Remark}

\newtheorem{example}[theorem]{Example}

\definecolor{sevencolor}{RGB}{201,98,21}

\newcommand\figurepdf[2][]{\raisebox{-.5\height}{\includegraphics[trim=72 72 72 72,#1]{figures_standalone/#2}}}

  \DeclareRobustCommand{\celltype}[1]{\texorpdfstring{\text{%
  \textpdfrender{TextRenderingMode=Fill,LineWidth=.05pt}%
{\normalfont\bfseries\oldstylenums{#1}}}}{#1}}
  \DeclareRobustCommand{\raytype}[1]{\texorpdfstring{\text{%
  \textpdfrender{TextRenderingMode=Fill,LineWidth=.1pt}{\boldmath$#1$}}}{#1}}

\newcommand{\Ggroup}{\operatorname{G}_2}
\newcommand{\Gliealgebra}{\mathfrak{g}_2}
\newcommand{\cartan}{\mathfrak{h}} 

\newcommand{\PP}{\mathbb{P}}
\newcommand{\RR}{\mathbb{R}}
\newcommand{\QQ}{\mathbb{Q}}
\newcommand{\CC}{\mathbb{C}}
\newcommand{\KK}{\mathbb{K}} 
\newcommand{\LL}{\mathbb{L}} 
\newcommand{\ZZ}{\mathbb{Z}}
\newcommand{\NN}{\mathbb{N}}
\newcommand{\FF}{\mathbb{F}}

\newcommand\pointspace{\mathcal{Q}}
\newcommand\linespace{\mathcal{G}}
\newcommand\flagspace{\mathcal{F}}
\newcommand\hexagonspace{\mathcal{H}}
\newcommand\modulispace{\mathcal{M}}
\newcommand\divisor{\mathcal{K}}
\newcommand\Dtwelve{{\operatorname{D}_{12}}}
\newcommand\hex{\operatorname{H}} 

\renewcommand{\epsilon}{\varepsilon}
\newcommand{\dominant}{\mathcal{D}}

\DeclareMathOperator\GL{GL} 
\DeclareMathOperator\SL{SL} 
\DeclareMathOperator\PGL{PGL} 
\DeclareMathOperator\Sym{S} 
\DeclareMathOperator\grassmannian{Gr} 
\DeclareMathOperator\orthgrassmannian{OGr} 
\DeclareMathOperator\SO{SO} 
\DeclareMathOperator\trop{trop} 
\DeclareMathOperator\val{val} 
\DeclareMathOperator\initial{in} 
\DeclareMathOperator\characteristic{char} 
\DeclareMathOperator\diag{diag} 
\DeclareMathOperator\relint{relint} 
\DeclareMathOperator\degree{degree} 
\DeclareMathOperator\Hom{Hom} 
\newcommand\one{\mathbf{1}} 
\newcommand\shortroots{\mathrm{short}} 
\newcommand\longroots{\mathrm{long}} 
\newcommand\lineset{\mathcal{A}} 

\newcommand\polymake{\texttt{polymake}\xspace}
\newcommand\oscar{\texttt{Oscar}\xspace}
\newcommand\gfan{\texttt{Gfan}\xspace}
  
\title{\bf Tropical Geometry of $\mathbf{\Ggroup}$}

\author{Michael Joswig \and Kevin Kuehn \and Bernd Sturmfels \and Mathieu Vall\'ee}
  
\date{}

\begin{document}
\maketitle

\begin{abstract} \noindent
  We tropicalize the Grassmannian and the flag variety of the algebraic group of type~$\Ggroup$.
  The two tropical varieties are strictly contained in the associated Dressians.
  Their points are trees in the $\Ggroup$-forest and decorated trees in the $\Ggroup$-orchard.
  These form hexagons in incidence geometry.
  We classify all toric degenerations and  their Newton--Okounkov polytopes.
  Mukai theory leads us to tropical canonical curves of genus~10.
  This article serves as a blueprint for studies of reductive groups in tropical geometry.
\end{abstract}

\section{Introduction}
\label{sec1}

The exceptional group $\Ggroup$ is an algebraic group of dimension $14$,
defined classically as the automorphism group of the octonions~\cite[\S2.2--2.3]{SpringerVeldkamp:2000}.
We work with its split form over an arbitrary field~$\KK$. Following  Aschbacher~\cite{Aschbacher:1987},
our starting point is the alternating trilinear~form
\begin{equation}
\label{eq:trilinear}
\Phi \,\, = \,\, e_{123} \,+\, e_{456} \,+\,
e_{147} \, + \, e_{257} \, + \, e_{367}\,\, \in \,\,
\wedge_3 \KK^7.
\end{equation}
Here, we identify $\KK^7$ with its dual via the standard basis $e_1,\ldots,e_7$, and we write $e_{ijk} = e_i \wedge e_j \wedge e_k$ for the induced basis vectors of $\wedge_3 \KK^7$.
We define $\Ggroup \subset \SL_7$ as the stabilizer of $\Phi$, so that its group of
$\KK$-points $\Ggroup(\KK)$ consists of the $7 \times 7$-matrices in $\SL_7(\KK)$ that fix $\Phi$.
Over every field $\KK$, this is the split simple algebraic group of type~$\Ggroup$,
as shown in~\cite[\S3]{Aschbacher:1987}.

We find that $\Ggroup$ is a subgroup of the
special orthogonal group $\SO(7,q)$ for
\begin{equation}
\label{eq:quadric}
q(x) \,\, = \,\, x_1 x_4 \,+\,x_2 x_5 \, + \,x_3 x_6 \,-\, x_7^2.
\end{equation}
This split quadratic form has signature  $(3,4)$ over $\RR$.
To derive $q(x)$, we write $\Phi$ as a skew-symmetric $7 \times 7$ matrix whose
entries are linear forms in $7$ variables $x_1,\ldots,x_7$, namely
 \begin{equation}
 \label{eq:Phi2}  \begin{matrix} 
\Phi(x, \,\cdot\,,\,\cdot \,) 
 &\! =\! &  \!
     x_1 \cdot (e_{23} + e_{47}) 
   \,-\, x_2 \cdot (e_{13} - e_{57}) 
  \,+\, x_3 \cdot (e_{12} + e_{67})  
 \, -\, x_4 \cdot (e_{17} - e_{56})  \smallskip \\ & &
 \, -\, x_5 \cdot (e_{27} + e_{46}) 
 \, -\, x_6 \cdot (e_{37} - e_{45}) 
 \, +\, x_7 \cdot (e_{14} + e_{25} + e_{36}) 
    \,\,\,\, \in \,\, \wedge_2 \KK^7.
   \end{matrix}
   \end{equation}
   Up to sign,    the $6 \times 6$ Pfaffians of (\ref{eq:Phi2}) 
   are $x_1 q(x), x_2 q(x), \ldots, x_7 q(x)$, so
   $\Phi$ determines  $q(x)$.

The group $\Ggroup$ has an interesting history, described by Agricola~\cite{Agricola:2008} and Manivel~\cite{Manivel:2025}.
It was discovered by Wilhelm Killing, who communicated it to Friedrich Engel in a letter dated May 23, 1887.
Engel worked with Felix Klein and Sophus Lie, both professors at Leipzig.
Engel presented $\Ggroup$ in a lecture at the Saxonian Academy of Sciences on June 11, 1900.

According to Klein, geometries emerge from groups,
and the lifeblood of a group is the spaces it acts on.
Our group $\Ggroup$ acts naturally on three spaces:
the quadric $\pointspace = V(q) \subset \PP^6$,
a Grassmannian $\linespace$ of certain  lines  on $\pointspace$, 
 and the flag variety $\flagspace \subset \pointspace \times \linespace$ which links~them.
Both $\pointspace$ and $\linespace$ have dimension $5$ and
 $\flagspace$ has dimension~$6$.
 The $\Ggroup$-Grassmannian  $\linespace$ is contained in
 the $7$-dimensional orthogonal Grassmannian $\orthgrassmannian(2,7)$
 of all lines on $\pointspace$, and hence in the $10$-dimensional
   Grassmannian $\grassmannian(2,7) \subset \PP^{20}$.
   In this paper we tropicalize these objects.
 
 We now describe the polynomials that define our varieties.
The ideal of the Grassmannian $\grassmannian(2,7)$ in Pl\"ucker coordinates
 is generated by the $4 \times 4$ subpfaffians
$p_{ij} p_{kl} - p_{ik} p_{jl} + p_{il} p_{jk}$
of a skew-symmetric $7 \times 7$ matrix $(p_{ij})$.
The tropical Grassmannian $\trop(\grassmannian(2,7))$ is the space
of phylogenetic trees with $7$ leaves~\cite[Theorem 4.3.5]{MaclaganSturmfels:2015}.
We here also need the linear space
\begin{equation}
\label{eq:7linear} \begin{matrix}
\!\!\! \mathcal{L} \,\, = \,\, \{ \,p \in \PP^{20} \,:\quad
 p_{23} + p_{47} \, = \, 
   p_{13} - p_{57} \, = \, 
  p_{12} + p_{67} \, = \,   
  p_{17} - p_{56}  \qquad \\   \qquad \qquad \qquad \quad \,\,=\,
  p_{27} + p_{46} \, = \, 
  p_{37} - p_{45} \, = \, \,
  p_{14} + p_{25} + p_{36} \,=\,0 \, \}.
  \end{matrix}
 \end{equation} 
 These $7$ linear equations are the coefficients of (\ref{eq:Phi2}) when viewed as a linear form in $x_1,\ldots,x_7$.
 The intersection in $\PP^{20}$ of the Grassmannian with this linear space is the   {\em $\Ggroup$-Grassmannian}:
 \begin{equation}
 \label{eq:grass1}
  \linespace \,\, = \,\, \grassmannian(2,7) \,\cap\, \mathcal{L} .
\end{equation}  
We have $\dim(\linespace) = 5$ and
$\degree(\linespace) = 18$.
Each point in $\linespace$ is a line in $\PP^6$, and
the union of these lines is the quadric $\pointspace$.
The {\em flag variety} $\flagspace$ for $\Ggroup$
is their incidence correspondence.
This is a $6$-dimensional variety,
defined by $34$ bilinear equations inside $\,\pointspace \times \linespace \subset \PP^6 \times \PP^{20}$.

\begin{figure}[ht]
  \centering
  \figurepdf{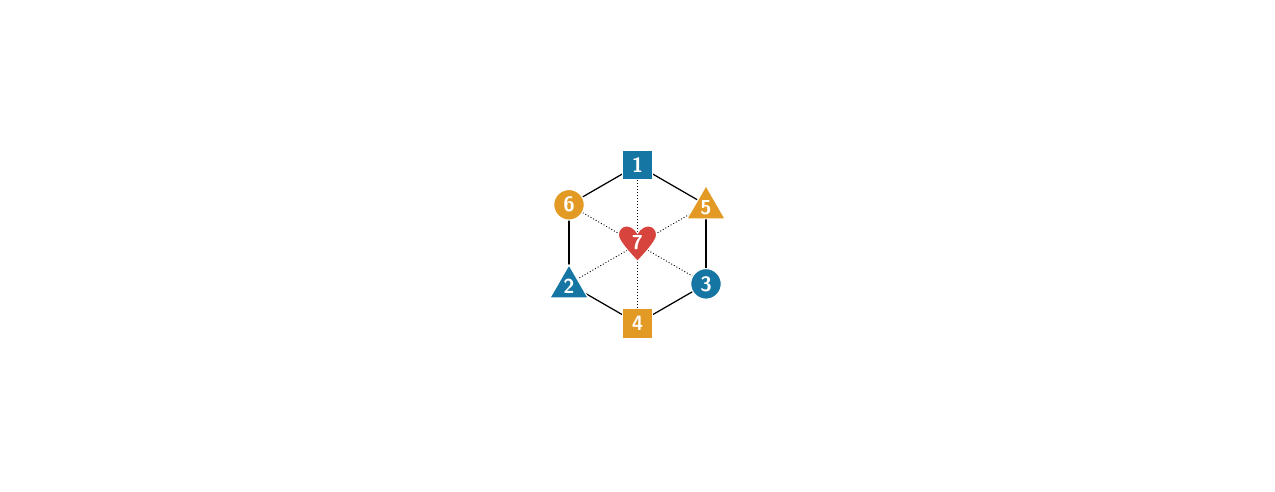}
  \caption{Colorful redrawing of the hexagon (\ref{eq:hexagon}). Its symmetry group, the dihedral group $\Dtwelve$
    of order $12$,   is the {\em Weyl group} of $\Ggroup$.
    This group acts on our classical and tropical varieties.}
  \label{fig:shape_hexagon}
\end{figure}

In Tits' theory of buildings~\cite{Tits:1974}, the triple $(\pointspace,\linespace,\flagspace)$ 
serves as a {\em generalized hexagon}~\cite{VanMaldeghem:1998}.
In our coordinates, a distinguished hexagon is given by a cycle of six points and lines:
\begin{equation}
\label{eq:hexagon}
 e_1 - e_{15} - e_5 -e_{35} - e_3 - e_{34} - e_4 - e_{24} - e_2 - e_{26} - e_6 - e_{16} - e_1. 
\end{equation}
The group $\Ggroup$ acts on $\flagspace$,
and this action is transitive on the non-degenerate hexagons.
The labels in Figure~\ref{fig:shape_hexagon} show how the Weyl group $\Dtwelve$ acts on the $7$ linear equations in (\ref{eq:7linear}).
Note that $\Phi$ in (\ref{eq:trilinear}) is the sum over two triangles and three diagonals 
that cross the hexagon.

A key point of tropical geometry is that tropicalization and intersection do not commute:
 \begin{equation}
 \label{eq:grass2}
  \trop(\linespace) \,\, = \,\, \trop\bigl( \grassmannian(2,7) \,\cap\, \mathcal{L} \bigr)\,\,\,
  \subsetneq \,\, \, \trop( \grassmannian(2,7) ) \,\cap \, \trop( \mathcal{L} ). 
\end{equation}  
Both objects in (\ref{eq:grass2}) are polyhedral fans in $\RR^{21}/\one$.
Only that on the left is a tropical variety. Namely, it
is the {\em tropical  $\Ggroup$-Grassmannian}. The object on the right is 
 a  tropical prevariety,  called the {\em $\Ggroup$-Dressian}.
 This follows the nomenclature in~\cite{MaclaganSturmfels:2015}.
We also study the {\em tropical $\Ggroup$ flag variety}
$\trop(\flagspace) $, which is a tropical variety in 
$\RR^{7}/\one \times \RR^{21}/\one$.
The {\em $\Ggroup$ flag Dressian} is the tropical prevariety
defined by the $78$ generators of the ideal of $\flagspace$ listed in Section~\ref{sec2}.

We now discuss the organization and main results in this paper.
In Section~\ref{sec2} we introduce {\em $\Ggroup$-trees}. These are
phylogenetic trees on $7$ leaves that lie on the $\Ggroup$-Dressian.
The $\Ggroup$ flag Dressian is the corresponding universal family. Its elements are
$\Ggroup$-trees with a marked point. 
Complete descriptions of our two tropical prevarieties 
are given in Theorems~\ref{thm:dressian} and~\ref{thm:flag_dressian}.

In Section~\ref{sec3} we study the two tropical varieties
associated with $\Ggroup$. Theorem~\ref{thm:grassmannian} characterizes the Grassmannian
$\trop(\linespace)$, and Theorem~\ref{thm:flag} characterizes the flag variety $\trop(\flagspace)$.
A tropical basis for $\linespace$ is presented in Corollary \ref{cor:tropicalbasis} and explained in
Figure \ref{fig:tropbasis}. These objects were found by algebraic computations 
over the field $\QQ$ of rational numbers.
In Section~\ref{sec4} we turn from characteristic $0$ to characteristic $p > 0$.
The tropical Grassmannian and flag  variety remain as in Section~\ref{sec3} for $p \geq 3$. 
However, they have different tropicalizations in characteristic $p=2$.
This is proved in Theorems \ref{thm:char2} and \ref{thm:flagchar2}.
Figure \ref{fig:tetrahedra} offers a visualization.

Section \ref{sec5} introduces the $\Ggroup$-forest.
This is our gallery of $\Ggroup$-trees, one for each  cell in $\trop(\linespace)$. 
The rays are shown in Figure~\ref{fig:atlas_rays}.
The maximal cells appear in Figures~\ref{fig:atlas_toric_cells_prime} and~\ref{fig:atlas_toric_cells_nontoric}.

In Section~\ref{sec6} we examine toric degenerations of
the Grassmannian $\linespace$ and its flag variety
$\flagspace$. Such degenerations and their Newton--Okounkov polytopes are of
considerable interest in the representation theory of $\Ggroup$.
Specific instances appeared earlier in the articles~\cite{CI:2014, FeiginFourierLittelmann:2017, SW:2024}.
Theorem~\ref{thm:NObody} gives the complete classification.
The census for $\linespace$ is presented
in Table~\ref{tab:g2-grassmannian-newton-okounkov}.

Section \ref{sec7} introduces the $\Ggroup$-orchard.
This is our gallery of decorated $\Ggroup$-trees, one for each 
cell in $\trop(\flagspace)$. We focus on the fiber complexes.
These are the preimages in $\trop(\flagspace)$ of individual maximal cells in $\trop(\linespace)$.
We also study the Newton--Okounkov polytopes for~$\flagspace$.

The present paper serves as a blueprint for importing reductive groups into tropical geometry.
In Section~\ref{sec:homogeneous} we explore what it means to tropicalize a group and its homogeneous spaces. 
Theorem \ref{thm:homogeneous} shows how our approach to $\Ggroup$ fits into the general theory.
Section~\ref{sec:hexagons} is motivated by the theory of buildings \cite{Tits:1974, VanMaldeghem:1998}.
We study the $12$-dimensional variety $\hexagonspace$ of all hexagons
on $\flagspace$, i.e., the closure of the orbit of (\ref{eq:hexagon}) under the $\Ggroup$-action.
Theorem \ref{thm:completing_hexagons} builds tropicalized hexagons.
These are  the loop trails for enjoying a stroll in the $\Ggroup$-orchard.

Our inspiration for this article was the theory
of algebraic curves of genus $10$. Their moduli space $\modulispace_{10}$
has dimension $27$. It has a distinguished divisor $\divisor$, given
by curves that lie on a K3 surface.
According to Mukai~\cite{Mukai, Mukai3},
the canonical model of a curve in~$\divisor$ is an intersection
$\Lambda \cap \linespace$, where $\Lambda$ is a  linear
space of codimension~$4$ in  $\PP^{20}$, so that $\Lambda \cap \mathcal{L} \simeq \PP^9$.
In Section~\ref{sec10}, we study  intersections of $\trop(\linespace)$
with tropical linear spaces  $\trop(\Lambda)$.
Our {\em tropical canonical curves} 
are derived in Proposition \ref{prop:mukai}
and shown in
Figures \ref{fig:planar_curve} and \ref{fig:nonplanar_curve}.

We now  comment on our usage of AI tools.
From the very start and throughout the entire project, tools based on large language models (mostly Claude Opus 5.5) were used for literature search, drawing pictures, writing code, drafting results and their proofs.
Most importantly, Claude was invaluable in maintaining the overall consistency of this  long paper.
This is particularly important because several of our main results are based on nontrivial computations, mostly with the computer algebra system \oscar~\cite{OSCAR-book,Oscar}.
To check that our computational arguments are really complete, the reader is invited to use similar tools.
Our Zenodo repository at \href{https://doi.org/10.5281/zenodo.24104003}{\texttt{10.5281/zenodo.24104003}} is organized in a way to facilitate linking the results in this article with \oscar output and computations; see Appendices~\ref{sec:data} and~\ref{sec:computation}.

\section{Trees and Dressians}
\label{sec2}

The  tropical Grassmannian $\trop(\grassmannian(2,7))$ is the set of points $d=(d_{ij})_{1\leq i<j\leq 7}\in \RR^{21}$ such~that
 \begin{equation}\label{eq:4point}
  \max(d_{ij}+d_{kl},d_{ik}+ d_{jl},d_{il}+ d_{jk})
\end{equation}
is attained at least twice~\cite{speyer_sturmfels}; see also \cite[Theorem 10.47]{ETC} and \cite[Lemma 4.3.6]{MaclaganSturmfels:2015}.
This is the {\em four-point condition} for phylogenetic trees.
Hence $d$ is a tree metric, provided its coordinates are positive and satisfy the triangle inequality.
In algebraic geometry, we usually tropicalize with respect to the min-convention \cite{MaclaganSturmfels:2015}.
This is also \oscar's default \cite{Oscar}. However, for our trees we pass to the 
max-convention by negating the valuation, so $d_{ij} = -\val(p_{ij})$ represents the tree distance between leaves $i$ and $j$.
This is because (\ref{eq:4point}) is a max-tropical constraint.

A \emph{$\Ggroup$-tree} is a phylogenetic tree $d$ with set of leaves $[7] = \{1,2,\ldots,7\}$ that 
satisfies \begin{equation}\label{eq:dressian}
	\begin{gathered}
		\begin{aligned}
      d_{17}&=d_{56}, & d_{27}&=d_{46}, & d_{37}&=d_{45},\\
			d_{47}&=d_{23}, & d_{57}&=d_{13}, & d_{67}&=d_{12},
		\end{aligned}\\
		\text{and } \max\{d_{14},d_{25},d_{36}\} \text{ is attained at least twice.}
	\end{gathered}
\end{equation}
The seven linear equations in \eqref{eq:7linear} form a tropical basis,
and (\ref{eq:dressian}) records their tropicalizations.
These are the constraints that define $\trop(\mathcal{L})$,
and hence the intersection on the right of~(\ref{eq:grass2}).
In conclusion, the $\Ggroup$-Dressian is the space of
phylogenetic trees $d$ that happen to be $\Ggroup$-trees.

To make the notion of $\Ggroup$-trees more tangible,
we use the heart and the six colored shapes in
Figure~\ref{fig:shape_hexagon} to represent the seven leaves of each tree.
The seven constraints in~\eqref{eq:dressian} are then:
\begin{itemize}
  \item[(i)] the distance from any shape to the heart is the same as the distance between the two other shapes of the 
  opposite color;
  \item[(ii)] the maximum of the three distances between equal shapes is attained at least twice.
\end{itemize}

 \begin{figure}[ht]
  \centering
  \figurepdf{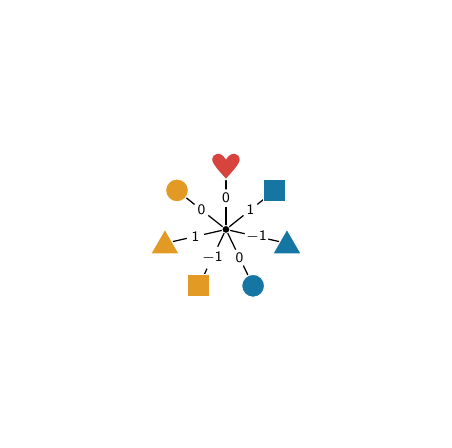}\hspace{19.58pt}\figurepdf{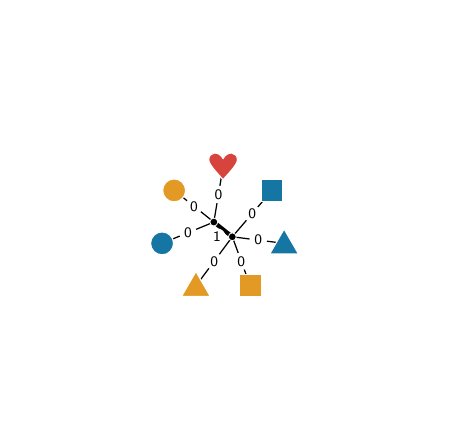}\hspace{19.58pt}\figurepdf{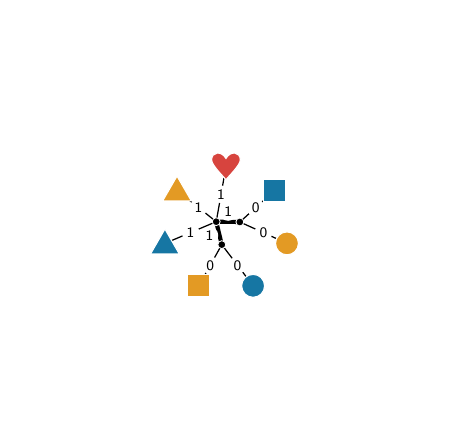}%
  \hspace{19.58pt}\figurepdf{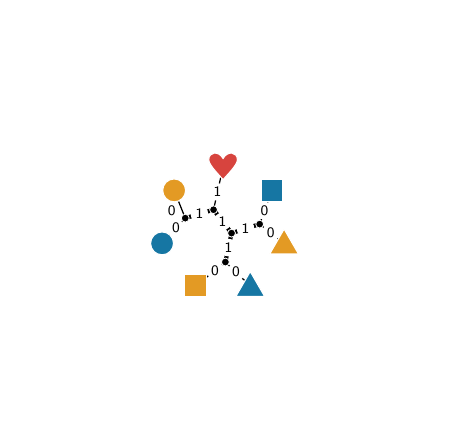}\hspace{19.58pt}\figurepdf{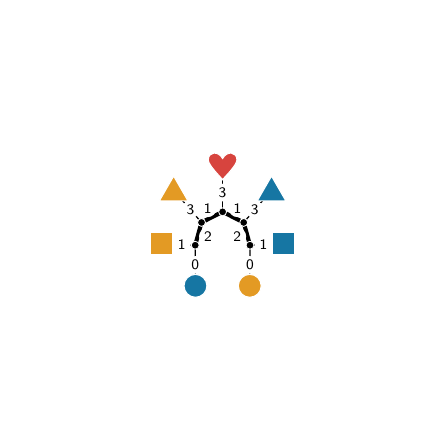}
  \caption{From left to right, a $\Ggroup$-tree with no inner edge, one inner edge, two inner edges, a snowflake $\Ggroup$-tree, and a caterpillar $\Ggroup$-tree. Check the conditions (i) and (ii) for each tree.
  }  \label{fig:ex_G2_trees}
\end{figure}

\begin{example} \label{ex:fivetrees}
Five $\Ggroup$-trees are shown in Figure~\ref{fig:ex_G2_trees}.
The rightmost tree gives $d \in \RR^{21}$~with
\[ \begin{matrix} d_{17}=d_{56} = 7, \,\, d_{27}=d_{46}=7, \,\,  d_{37}=d_{45} = 6 ,\,\,
d_{47}=d_{23} = 7, \,\, d_{57}=d_{13} = 7 ,\,\,  d_{67}=d_{12} = 6 ,\\
d_{36} = 6 \,<\, d_{14} = d_{25} = 8,\quad
d_{15} = 8,\,\,
d_{16} = 1, \,\,
d_{24} = 8,\,\,
d_{26} = 5, \,\,
d_{34} = 1,\,\,
d_{35} = 5. 
\end{matrix}
\]
The last six distances are the edges of our hexagon 
in Figure~\ref{fig:shape_hexagon}. They
do not appear in (\ref{eq:dressian}).
\end{example}

The lineality space of the tropical Grassmannian 
$\trop(\grassmannian(2,7))$ consists of all star trees~$d$. 
It has dimension $6$ in $\RR^{21}/\one$. Namely,
if $s_i$ is the length of the $i$th pendant edge in the star tree
then  $d_{ij}=s_i+s_j$ for all $i,j $.
We now describe the lineality space of the $\Ggroup$-Dressian.

\begin{figure}[ht]
  \centering
  \figurepdf{lineality_star_1}\hspace{29.6pt}\figurepdf{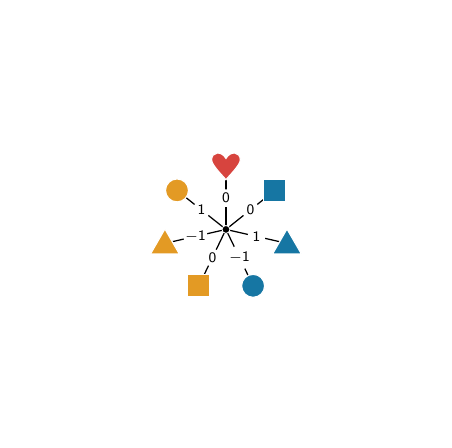}\hspace{29.6pt}\figurepdf{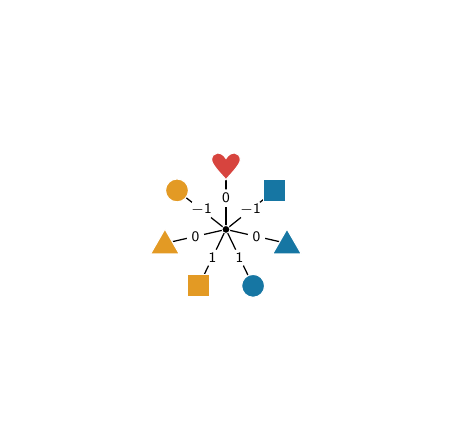}
  \caption{These three stars are  $\Ggroup$-trees. They
  span the lineality space of the $\Ggroup$-Dressian.}
  \label{fig:lineality_stars}
\end{figure}

\begin{lemma}\label{lem:lineality-space}
    The lineality space of the $\Ggroup$-Dressian has dimension $2$ in $\RR^{21}/\one$.
     It is the intersection of
    the lineality space of $\trop(\grassmannian(2,7))$ with the six hyperplanes  in~\eqref{eq:dressian}. 
    A basis~is given by  any two of the star trees in Figure~\ref{fig:lineality_stars}.
    Their edge lengths are the rows of the matrix
\begin{equation}
\label{eq:lineality}        
    \bordermatrix{  & \,\,s_1 & \,\,s_2 & \,\,s_3 & \,\,s_4 & \,\,s_5 & \,\,s_6 & \,\,s_7 \cr
 &   \phantom{-}1 & -1 & \phantom{-}0 & -1 & \phantom{-}1 & \phantom{-}0 & \phantom{-}0 \,\cr
 & \phantom{-} 0 &   \phantom{-}1 & -1 &   \phantom{-}0 & -1 &   \phantom{-}1 & \phantom{-} 0\,\cr
&		-1 & \phantom{-}0 & \phantom{-}1 & \phantom{-}1 & \phantom{-}0 & -1 & \phantom{-}0 \, }.		
\end{equation}
\end{lemma}
\begin{proof}
The $\Ggroup$-Dressian is the intersection in $\RR^{21}/\one$ 
of the two fans on the right in (\ref{eq:grass2}).
Its lineality space is therefore the intersection 
of the lineality spaces of the two intersectands.
\end{proof}

\begin{figure}[th]
  \centering
  \figurepdf{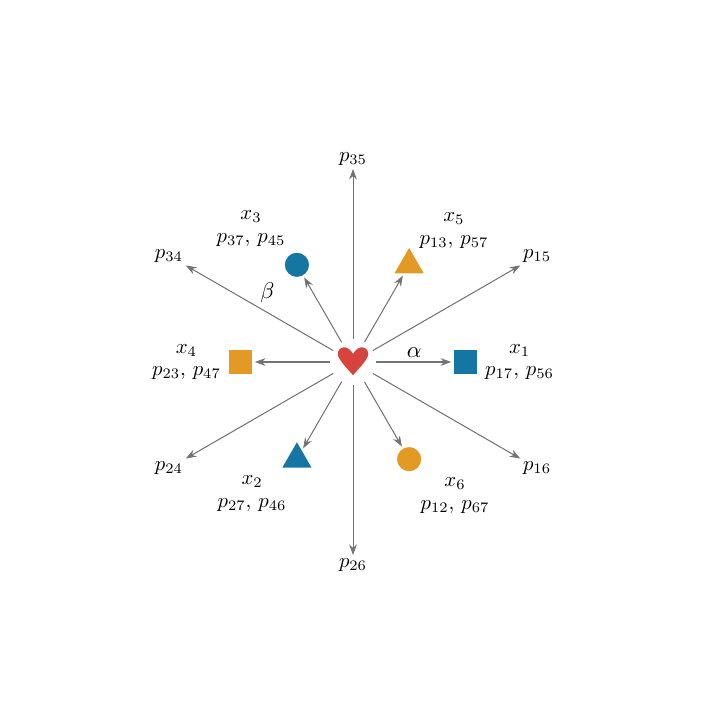}
  \caption{Root system of $\Ggroup$.
    Next to each root are the variables of that weight.
    Variables $x_7, p_{14}, p_{25}, p_{36}$ have weight $0$.
    The maximal torus $(\KK^*)^2$ acts on $\pointspace, \linespace,\flagspace$ with these weights.}
    \label{fig:root_lattice}
\end{figure}

Our lineality space mirrors the root system of $\Ggroup$, which is shown in Figure~\ref{fig:root_lattice}.
Leaves of $\Ggroup$-trees are labeled by the six short roots, like $\alpha=(1,0,-1)$, or by the origin (heart).
The six long roots include $\beta =(-1,-1,2)$.
The roots are the degrees of our variables $x_i$ and $p_{ij}$.

A $\Ggroup$-tree is \emph{normalized} if its edge lengths are relatively prime in $\NN$, its shortest pendant edge has length $0$, and its vector $s \in \RR^7$ of pendant edge lengths is in the kernel of~\eqref{eq:lineality}.
Two $\Ggroup$-trees are {\em equivalent} if some positive multiples are congruent modulo lineality space.

\begin{proposition}\label{prop:normalized}
  Every rational $\Ggroup$-tree is equivalent to a unique normalized $\Ggroup$-tree.
\end{proposition}

\begin{proof}
  By Lemma~\ref{lem:lineality-space}, lineality translates the vector $s$ of pendant edge lengths
  by the span~of $\one = (1,1,\ldots,1)$ and the rows of~\eqref{eq:lineality}, and changes nothing else.
  Exactly one translate of $s$ is orthogonal to the plane spanned by the three rows.
  Their entries sum to zero, so adding a multiple of $\one$ preserves orthogonality.
  Exactly one multiple makes the shortest pendant edge $0$.
  Both conditions are invariant under positive scaling, and the resulting tree $d$ is rational and nonzero.
  Hence exactly one positive multiple has integer edge lengths with gcd $1$.
\end{proof}

The Weyl group of~$\Ggroup$ is the dihedral group of order $12$, denoted $\Dtwelve$.
For us, $\Dtwelve$ is the permutation group 
 generated by the rotation $(1\,5\,3\,4\,2\,6)$ and the reflection $(2\,3)(5\,6)$.
Geometrically, these are symmetries of the hexagons in
Figures~\ref{fig:shape_hexagon} and~\ref{fig:root_lattice}.
It permutes the coordinates $p_{ij}$ of $\PP^{20}$, up to signs, and the coordinates $d_{ij}$ of $\RR^{21}/\one$.
See also the proof of Theorem \ref{thm:homogeneous} (b).
In this manner, $\Dtwelve$ acts on the varieties $\pointspace$, $\linespace$, $\flagspace$  and their tropical varieties and prevarieties.
Equipped with this $\Dtwelve$-action, we can now describe the $\Ggroup$-Dressian.

\begin{theorem}\label{thm:dressian}
  Modulo lineality, the $\Ggroup$-Dressian is a $3$-dimensional polyhedral complex with $f$-vector $(51, 241, 248, 3)$.
  Its $2$-dimensional cells are $245$ triangles and three quadrilaterals, and its three $3$-cells are tetrahedra.
  Modulo the $\Dtwelve$-action, the $f$-vector is $(11,33,31,1)$.
  The reduced integral homology of the $\Ggroup$-Dressian is  $\,\ZZ^{54}$, concentrated in dimension $2$.
\end{theorem}

\begin{proof}
  This is a computation in the computer algebra system \oscar~\cite{OSCAR-book,Oscar}.
  We compute the tropical prevariety of the ideal $I$ given by the $35$ Pl\"ucker relations and the $7$ linear relations~\eqref{eq:7linear}.
  The output is a fan $\Sigma$ over a polyhedral complex with $f$-vector $(51, 241, 248, 3)$.
  That fan lies in a 15-dimensional subspace cut out by the six linearly independent equations in (\ref{eq:dressian}).
  Among the $248$ polygons, we find $3$ quadrilaterals and $245$ triangles.
  After removing the $12$  facets of the tetrahedra, we have $236$ polygons that are maximal cells.
  In total, there are $239$ maximal cells, namely $233$ triangles, $3$ quadrilaterals and $3$ tetrahedra.
  Their cones span different linear spaces, so they give the coarsest fan structure on the support of $\Sigma$.

  The action of $\Dtwelve$ permutes the $51$ vertices and the $239$ facets of the polyhedral complex
  underlying the fan  $\Sigma$. The vertices come in $11$ orbits.
  Their $\Ggroup$-trees are shown in Figure~\ref{fig:atlas_rays}.
  The maximal triangles and quadrilaterals form $28$ orbits.
  Their $\Ggroup$-trees are the ones in Figures~\ref{fig:atlas_toric_cells_prime} and~\ref{fig:atlas_toric_cells_nontoric}, except those labeled by ``char~$2$'' or  ``char $\neq2$''.
  The tetrahedra form one orbit, represented by the $\Ggroup$-tree labeled \raytype{\epsilon} in Figure~\ref{fig:atlas_rays}.
  All trees in these figures are normalized as in Proposition~\ref{prop:normalized}.
  The homology groups are also computed in $\oscar$.
\end{proof}

The $\Ggroup$ flag Dressian is a tropical
prevariety in $\RR^7/\one \times \RR^{21}/\one$.
We now state its definition in precise terms.
The bihomogeneous prime ideal of the flag variety $\flagspace \subset \PP^6 \times \PP^{20}$
lives in the polynomial ring $\KK[x,p]$  with $28$ variables.
We use the following $1+7+\binom{7}{4} + \binom{7}{3} = 78$ generators of this ideal.
Namely, in bidegree $(2,0)$ we have our quadric
 $q$, in bidegree $(0,1)$ we have the $7$
linear forms in (\ref{eq:7linear}), in bidegree $(0,2)$ there are
$35$ Pl\"ucker quadrics, and in bidegree $(1,1)$ we have the
$35$ incidence equations $p_{ij} x_k - p_{ik} x_j + p_{jk} x_i = 0$. These
ensure that the point $x$ lies on the line $p$.
Modulo the products of the $x_i$ with the $7$ linear forms, only $34$ of them 
are linearly independent. This is the number given in the Introduction.
The {\em $\Ggroup$ flag Dressian} is the tropical prevariety 
defined by these $78$ polynomials. Its lineality space is $2$-dimensional,
given by concatenations of each row vector in (\ref{eq:lineality})
with the associated star tree $d$. The former with the vector $(1,1,1,-1,-1,-1,0)$ span the lineality of the $3$-dimensional tropical hypersurface
$\trop(\pointspace)$ and the latter span the lineality of the $\Ggroup$-Dressian.

The tropical incidence relations are  $\max(d_{ij} \!+\! d_{k8}, d_{ik} \!+\! d_{j8}, d_{jk} \!+\! d_{i8})$,
where $d_{k8} = - \val(x_k)$.
They define the tropical Grassmannian $\trop\grassmannian(2,8)$ over $\trop\grassmannian(2,7)$.
Therefore, the points $(d_{ij},d_{k8})$ of the $\Ggroup$ flag Dressian  
are $\Ggroup$-trees with one additional free leaf $8$.
The pendant edge at $8$ can be normalized to $0$.
The resulting point in $\RR^7/\one \times \RR^{21}/\one$ is a \emph{decorated $\Ggroup$-tree}.

\begin{theorem}\label{thm:flag_dressian}
  Modulo lineality, the $\Ggroup$ flag Dressian is a $4$-dimensional polyhedral complex with $f$-vector $(141, 1352, 3584, 2729, 33)$.
  The $3$-cells are $2603$ tetrahedra, $108$ square pyramids and $18$ triangular prisms. The $4$-cells are $33$ $4$-simplices. 
  The reduced integral homology is $\ZZ^{324}$, concentrated in dimension $3$.
  Modulo $\Dtwelve$-action, the $f$-vector is $(27, 173, 388, 274, 6)$.
\end{theorem}

\begin{proof}
  We used \oscar as well as \gfan~\cite{gfan} to compute the polyhedral fan.
  This fan has $2726$ maximal cones.
  Exactly $99$ pairs of them span the same linear space, and each pair has a convex union.
  Merging them removes six rays and gives $2627$ maximal cones with pairwise distinct spans.
  The resulting fan is the coarsest, with $f$-vector $(141, 1352, 3584, 2729, 33)$.

Points in our flag Dressian are decorated $\Ggroup$-trees.
  The lineality space is the diagonal: it is made of the star $\Ggroup$-trees with decorated inner vertex.
  The rays are $\Ggroup$-trees where one inner vertex is decorated. From Figure~\ref{fig:atlas_rays}, we count $2\cdot 20+3\cdot 30 + 4=134$ rays. The $7$ remaining rays lie over the lineality space of $\trop(\linespace)$: the tree is the star tree, and the extra point $8$ lies on a pendant edge.
  Only the index of the pendant edge matters: the supporting star tree can be moved freely in the lineality space.
  We similarly obtain maximal cells by decorating edges
  of the $\Ggroup$-trees on the  maximal cells of the $\Ggroup$-Dressian.
  For example, the $33$ $4$-simplices are obtained by decorating one of the $11$ 
  edges of the $3$ tetrahedral $\Ggroup$-trees.
  Decorating an edge of a $2$-cell $\Ggroup$-tree produces $2714$ $3$-cells.
  The remaining $15$ are obtained by decorating one of the $5$ nodes of the $3$ tetrahedral $\Ggroup$-tree.
  The $\Dtwelve$-action permutes rays and preserves maximal cells.
The homology is  computed with \oscar as in Theorem~\ref{thm:dressian}.
\end{proof}

\section{Grassmannian and Flag Variety}
\label{sec3}

The tropical $\Ggroup$-Grassmannian is contained in the $\Ggroup$-Dressian.
This inclusion, shown in~(\ref{eq:grass2}), must be strict  by
Theorem~\ref{thm:dressian}. The cones over the three tetrahedra have dimension larger than
 $\dim(\linespace)=5$, whereas $\trop(\linespace)$ is pure of dimension $5$,
by the Structure Theorem~\cite[Theorem 3.3.5]{MaclaganSturmfels:2015}.
 Our next theorem shows that the tetrahedra are the only cells in the difference.
  For all results in this section, we assume that the field $\KK$ has characteristic~$0$.

\begin{theorem}\label{thm:grassmannian}
  Modulo lineality, the tropical $\Ggroup$-Grassmannian $\trop(\linespace)$
   is a $2$-dimensional polyhedral complex with $f$-vector $(51, 241, 248)$.
  Its $2$-cells are $245$ triangles and three quadrilaterals.
  The reduced integral homology is $\ZZ^{57}$, concentrated  in dimension $2$.
  Modulo the $\Dtwelve$-action, the $f$-vector  of  $\,\trop(\linespace)$ is $(11, 33, 31)$.
  Each cell is encoded by a $\Ggroup$-tree. 
  \end{theorem}
  
\begin{proof}
  We perform the computation in \oscar~\cite{OSCAR-book}.
  The input is the same ideal $I$ as in Theorem~\ref{thm:dressian}.
  The output is the tropical variety, with its Gr\"obner fan structure. 
  It has  the same lineality space as the Dressian.
  Modulo lineality, the fan is simplicial and has $f$-vector $(51,244,251)$.
  The $251$ maximal cones span $248$ distinct linear spaces: three pairs have the same span.
  We merge each pair of triangles into a quadrilateral, by deleting its diagonal.
  
  This leaves $245$ triangles and three quadrilaterals. The $248$ maximal cones 
  have pairwise distinct spans, so we have the coarsest fan structure.
  This gives the $f$-vector in Theorem~\ref{thm:grassmannian}.
  The three quadrilaterals are precisely those of the $\Ggroup$-Dressian.
  Thus, the Grassmannian is obtained from the Dressian by removing the
  three solid tetrahedra (cf.\ Theorem~\ref{thm:dressian}).
The computation of orbits and of homology is analogous to 
those in   the proof of Theorem~\ref{thm:dressian}.
\end{proof}

The $11$ types of rays are shown in Figure \ref{fig:atlas_rays}
and the $31$ types of maximal cells are shown in Figures \ref{fig:atlas_toric_cells_prime} and \ref{fig:atlas_toric_cells_nontoric}.
These catalogues represent many refined aspects. Here is a first illustration.

\begin{remark}
All cones in the tropical Grassmannian $\trop(\grassmannian(2,n))$ have weight $1$.
This is not the case for $\trop(\linespace)$. One of the $31$ types of maximal cones
has multiplicity $2$. This is type \celltype{21}, on the upper left in
Figure \ref{fig:atlas_toric_cells_nontoric}. 
See Example \ref{ex:twentyone} for an algebraic explanation.
All other cones have weights $1$. With this,
the weighted fan $\trop(\linespace)$ is balanced 
\cite[Theorem~3.3.5]{MaclaganSturmfels:2015}.
\end{remark}

We now turn to the flag variety $\flagspace$.
A direct \oscar computation of the tropical variety did not terminate on the ideal of $\flagspace$.
Instead we exploit symmetry with the following algorithm.

\begin{algorithm}[htbp]
  \caption{Tropical variety via symmetric pruning}
  \label{alg:balanced_support}
  \begin{algorithmic}[1]
    \Require A homogeneous prime ideal $J=\langle S\rangle$ of dimension $d$, the tropical prevariety $\Sigma$ of the set $S$,
    and a group $\Gamma$ of signed permutations of the variables that preserves $J$ and~$\Sigma$.
    \Ensure The maximal cones of a fan with support $\trop(V(J))$
    \ForAll{$\Gamma$-orbits of cones $\sigma \in \Sigma$ with $\dim \sigma \geq d+1$}
      \State Compute a Gr\"obner basis of $J$ for a weight in $\relint(\sigma)$.
      \State Find $g \in J$ with $\initial_w(g)$ a monomial for all $w \in \relint(\sigma)$; if none, \Return fail.\label{alg:dominant}
    \EndFor
    \ForAll{$\Gamma$-orbits of cones $\tau \in \Sigma$ with $\dim \tau = d$}
      \State Let $w$ be the sum of the rays of $\tau$.
      \State Keep the orbit of $\tau$ if $\initial_w(J)$ contains no monomial.
    \EndFor
    \State \Return the kept cones
  \end{algorithmic}
\end{algorithm}

The correctness of Algorithm \ref{alg:balanced_support}
follows from basics in tropical geometry, including the
Structure Theorem~\cite[Theorem~3.3.5]{MaclaganSturmfels:2015}
and the Fundamental Theorem~\cite[Theorem~3.2.3]{MaclaganSturmfels:2015}.


\begin{theorem}\label{thm:flag}
  Modulo lineality, the tropical $\Ggroup$ flag variety $\trop (\flagspace)$ is a $3$-dimensional polyhedral complex with $f$-vector $(141, 1352, 3584, 2714)$.
  Its $3$-cells are $2588$ tetrahedra, $108$ square pyramids and $18$ triangular prisms.
  The reduced integral homology is $\ZZ^{342}$, concentrated  in dimension $3$.
  Modulo the $\Dtwelve$-action, the $f$-vector is $(27, 173, 388, 270)$. Each cell in $\trop (\flagspace)$
  is encoded by a $\Ggroup$-tree that is decorated as in Section \ref{sec7}.
\end{theorem}

\begin{proof}
   We apply Algorithm~\ref{alg:balanced_support} to the ideal $J$ of $\flagspace\subset\pointspace\times\linespace$,
  with $S$ its $78$ generators and $\Gamma = \Dtwelve$.
  For some vector $w_0$, 
  the initial ideal $\initial_{w_0}(J)$ is generated by $83$ binomials and is prime over~$\overline{\QQ}$.
  Hence $J$ is prime as the algorithm requires.
  The prevariety of $S$ is the flag Dressian of Theorem~\ref{thm:flag_dressian}, having $2627$ maximal cells, of dimensions $3$ and $4$ modulo lineality.
  The first loop removes the $33$ $4$-simplices, in $6$ orbits.
  This leaves the $2729$ $3$-cells of the flag Dressian, in $274$ orbits.
  The second loop keeps $270$ of these orbits, with $2714$ cells.
  The output is the coarsest and has $f$-vector $(141,1352,3584,2714)$, with $2588$ tetrahedra, $108$ square pyramids and $18$ triangular prisms. 
  It has the same lineality space as its Dressian.
  Projection to the second factor maps $\trop(\flagspace)$ onto $\trop(\linespace)$.
  Each of the $2714$ maximal cells maps onto one of the $248$ maximal cells of $\trop(\linespace)$.
  Each maximal cell of $\trop(\flagspace)$ is  a maximal $\Ggroup$-tree decorated by a point in the relative interior of an edge.
  This removes the $15$ tetrahedra of the $\Ggroup$ flag Dressian corresponding to node-decorated tetrahedral $\Ggroup$-trees.
  The $141$ rays are the same as those in the $\Ggroup$ flag Dressian.
  Namely, $134$ lie over rays of $\trop(\linespace)$, and all $51$ of these occur.
  Each ray is a minimal $\Ggroup$-tree with a marked point $8$.
  The remaining $7$ rays lie over the lineality space of $\trop(\linespace)$: the tree is the star tree, with the
   point $8$  on an edge.
\end{proof}

Theorems~\ref{thm:grassmannian} and~\ref{thm:dressian} show that the
 Dressian and the tropical Grassmannian are different.
 Our $\Ggroup$-theory is richer than the well-known identification
of phylogenetic trees with points in $\trop(\grassmannian(2,n))$. Namely,
we learned that the
$7+35$  linear and quadratic generators for the prime ideal $I$ of the
Grassmannian $\linespace$ do not form a  {\em tropical basis}
\cite[Definition 2.6.3]{MaclaganSturmfels:2015}.
We now remedy this point by augmenting our given generators to a tropical basis for $\linespace$.

\begin{corollary} \label{cor:tropicalbasis}
  The $35$ Pl\"ucker relations, the $7$ linear relations~\eqref{eq:7linear}, and the $\,3$ cubics
\begin{equation}
\label{eq:threecubics}
\begin{matrix}
    2p_{45}p_{46}p_{56} + p_{23}p_{36}p_{56} - p_{14}p_{12}p_{45} + p_{26}p_{13}p_{45} - p_{35}p_{12}p_{46} \\
    2p_{23}p_{12}p_{13} - p_{46}p_{14}p_{13} + p_{25}p_{56}p_{23} + p_{16}p_{45}p_{23} - p_{34}p_{56}p_{12}\\
    2p_{46}p_{56}p_{45} + p_{12}p_{25}p_{45} - p_{36}p_{13}p_{46} + p_{15}p_{23}p_{46} - p_{24}p_{13}p_{56}
\end{matrix}    
\end{equation}
  form a tropical basis for the ideal $I$ that defines the $\Ggroup$-Grassmannian $\linespace$ in $\,\PP^{20}$.
\end{corollary}

\begin{figure}[th]
  \centering
  \figurepdf[scale=0.9]{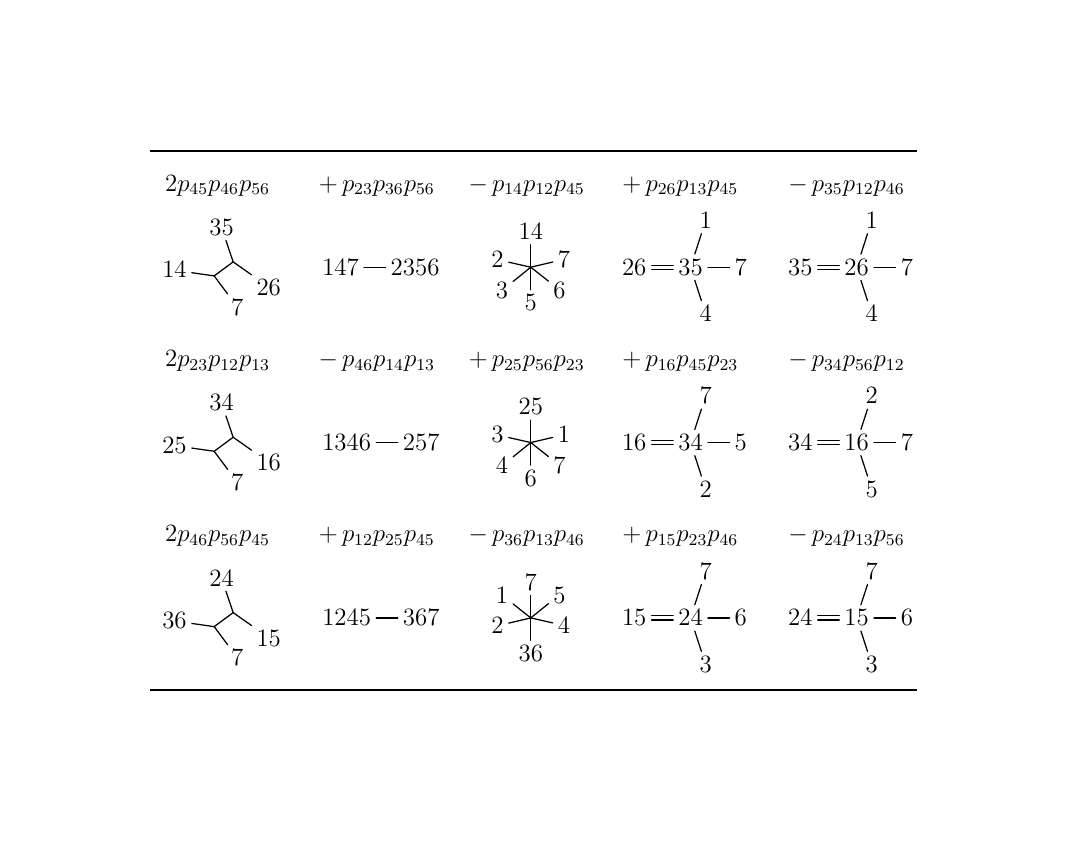}
  \caption{We associate a $\Ggroup$-tree to each monomial in the cubics (\ref{eq:threecubics}). 
  These are rays of $\trop_2(\linespace)$ and hence shown in Figure~\ref{fig:atlas_rays}.
  From left to right, their orbits are 
    \raytype{\epsilon}, \raytype{\alpha_2}, \raytype{\alpha_4}, \raytype{\alpha_3}, \raytype{\alpha_3}. }
      \label{fig:tropbasis}
\end{figure}

\begin{proof}
One first checks that each of the three cubics in (\ref{eq:threecubics}) lies in the ideal $I$ of $\linespace$.
To show that our $45 = 35+7+3$ polynomials form a tropical basis, we must prove the following statement
about an arbitrary element $d $ in $ \RR^{21}/\one$.
If $d$  does not lie in $\trop(\mathcal{G})$
then our list of $45$ candidates contains one polynomial $f$
whose initial form  $\initial_d (f)$ is a monomial.
This claim clearly holds for all $d$ that are not in the $\Ggroup$-Dressian:
here we can take $f$ from among the $42= 35+7$ original generators.
We may therefore assume that $d$ lies in the $\Ggroup$-Dressian.
  
From the description of the maximal cells in Theorem~\ref{thm:dressian},
we know that $d$ lies in the relative interior of one of the three special tetrahedra.
To examine the tetrahedra in  more detail,
we showcase $15$ special $\Ggroup$-trees, one for each of the 
$5+5+5$ monomials appearing in $(\ref{eq:threecubics})$.
These $\Ggroup$-trees are depicted in Figure~\ref{fig:tropbasis}.
Their edges have length $1$, except for the double-edges, which have length $2$.
For instance, the trees in the second column are splits. 
For the last tree in the first row, the tree metric $d$ satisfies
$d_{13} = 3$, $d_{23} = 2$, $d_{12} = 1$, $d_{35} = 0$.

Each row in Figure~\ref{fig:tropbasis} represents one of the three tetrahedra.
The last four $\Ggroup$-trees are the vertices of the tetrahedron,
and the leftmost $\Ggroup$-tree is in the interior of the tetrahedron.
One checks the following property for the cubic polynomial $f$ 
that is written above the trees. For each of the four trees $d$ on the right, the
initial form $\initial_d(f)$ is the sum of the four other monomials.
This verifies the fact, already known from Theorem~\ref{thm:grassmannian},
 that $d$ lies in $\trop(\linespace)$. 
 
 Consider the leftmost tree $d$. The initial form $\initial_d(f)$
is the monomial above it. This manifests $d \not\in \trop(\linespace)$. 
We now note that the initial form $\initial_d(f)$ remains unchanged if $d$ gets
replaced by any positive linear combination of the four $\Ggroup$-trees on the right.
For every $d$ in the interior of the tetrahedron, $\initial_d(f)$ is a monomial.
Hence the cubic $f$ is a witness for the entire tetrahedron.
This works for each of the three rows, and this completes the proof.
\end{proof}


\section{Characteristic Two}
\label{sec4}

The integer $2$ appears as the leading coefficient of the cubics (\ref{eq:threecubics}) in our tropical basis.
This suggests that the tropicalization of the  Grassmannian $\linespace$  depends on the field characteristic.
In this section we show that this is indeed the case. In what follows, $\KK$ denotes an arbitrary field with trivial 
valuation. Since $\linespace$ is defined over $\ZZ$,  its tropicalization depends only on the characteristic
$p = \characteristic(\KK)$, and we write $\trop_p(\linespace)$ for this tropical variety in $\RR^{21}/\one$.
The tropical variety for $p=0$ is the fan $\trop_0(\linespace)= \trop(\linespace)$,
which is described in Theorem~\ref{thm:grassmannian}.

\begin{theorem} \label{thm:char2}
  If $p = \characteristic(\KK) $ is an odd prime then $\trop_p(\linespace)= \trop(\linespace)$.
 Now assume~$p =  2$.
  Modulo lineality, the tropical $\Ggroup$-Grassmannian $\trop_2(\linespace)$
  is a $2$-dimensional polyhedral complex  with $f$-vector $(54, 253, 254)$.
  Its $2$-cells are $251$ triangles and three quadrilaterals.
  The reduced integral homology is $\ZZ^{54}$, concentrated  in dimension $2$.
  Modulo $\Dtwelve$-action, the $f$-vector  of  $\,\trop_2(\linespace)$ is $(12, 36, 32)$.
  Each cell is encoded by a $\Ggroup$-tree (see Figures \ref{fig:atlas_rays}, \ref{fig:atlas_toric_cells_prime}, \ref{fig:atlas_toric_cells_nontoric}).
\end{theorem}

The phenomenon of characteristic-dependence is familiar from \cite[Section 7]{speyer_sturmfels},
where it was shown for the tropicalization of $\grassmannian(3,7)$.  In the classical case, 
by \cite[Theorem 7.1]{speyer_sturmfels},
$\trop(\grassmannian(2,n))$ is characteristic-free: every tree lifts to a line in $\PP^{n-1}$.
We here show that the situation is different on the locus in $\grassmannian(2,7)$ that is
relevant for the exceptional group $\Ggroup$.

\begin{figure}[th]
  \centering
  \begin{tikzpicture}
  \node[inner sep=0pt,anchor=south west] (img) at (0,0)
    {\includegraphics[trim=96.7 64.7 71.5 53.0,clip,width=3.434in]{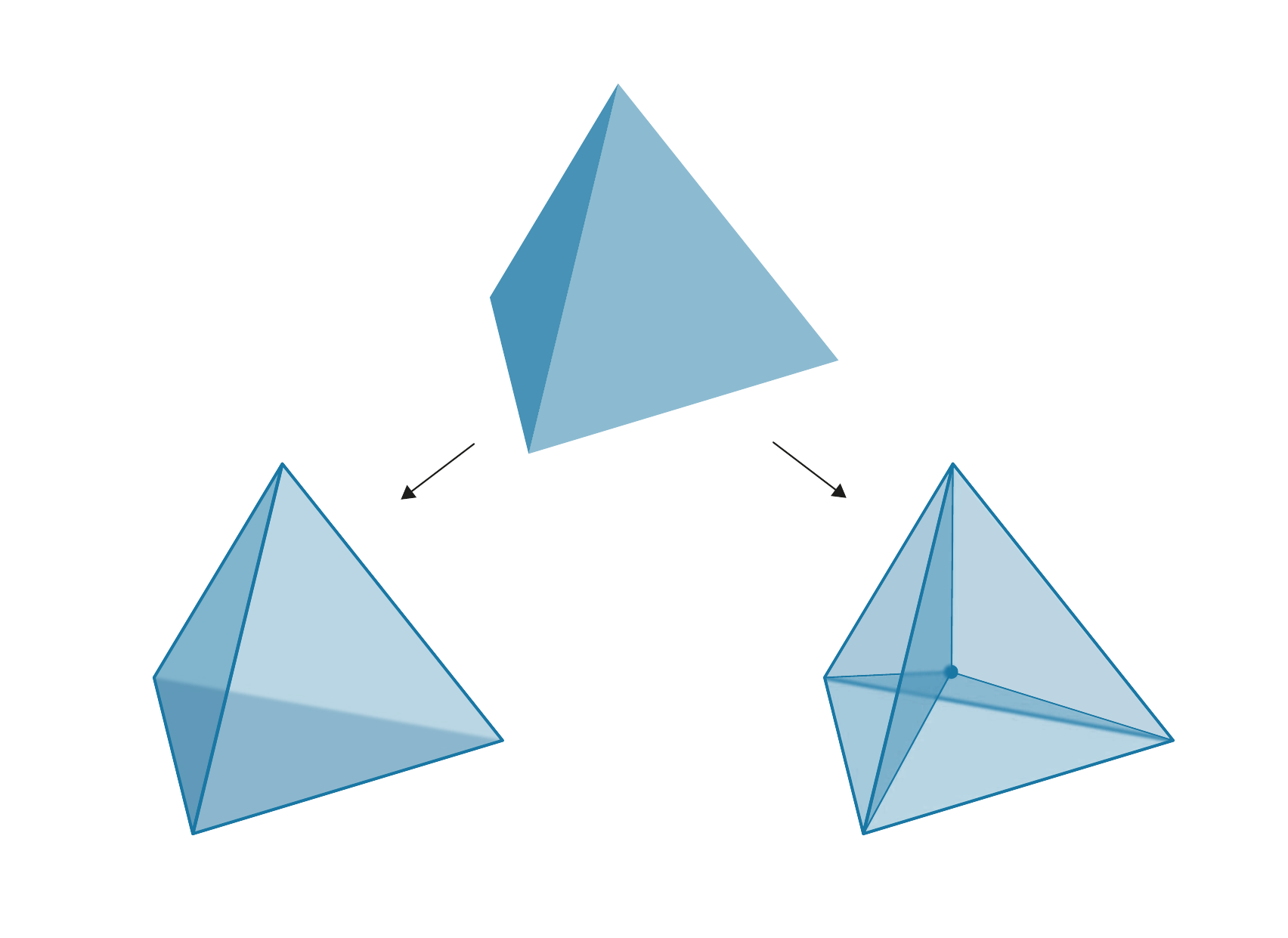}};
  \begin{scope}[x={(img.south east)},y={(img.north west)}]
    \node[above] at (0.455,1.000) {$\raytype{\alpha_2}$};
    \node[left] at (0.330,0.717) {$\raytype{\alpha_3}$};
    \node[right] at (0.671,0.632) {$\raytype{\alpha_3}$};
    \node[below left=-2pt] at (0.4,0.495) {$\raytype{\alpha_4}$};
    \node[above right=-1pt] at (0.775,0.22) {$\raytype{\epsilon}$};
    \node[below] at (0.241,0.051) {$\characteristic\KK\neq 2$};
    \node[below] at (0.909,0.051) {$\characteristic\KK=2$};
  \end{scope}
  \end{tikzpicture}
  \caption{The tetrahedron \raytype{\epsilon} shows how the
     tropical $\Ggroup$-Grassmannian depends on the field~$\KK$.}
  \label{fig:tetrahedra}
\end{figure}

Before proving Theorem~\ref{thm:char2},
we  explain its combinatorial essence. The characteristic~$2$ features emerge in the three tetrahedra of the $\Ggroup$-Dressian.
Each tetrahedron has a new ray, labeled  \raytype{\epsilon} in Figure~\ref{fig:atlas_rays}, that appears in its relative interior.
For $p=2$, the four facets of the tetrahedra are deleted and replaced by six triangles joining the center \raytype{\epsilon} to the six~edges.
This is shown in Figure~\ref{fig:tetrahedra}. The vertices of the tetrahedron are $\Ggroup$-trees of types
\raytype{\alpha_2}, \raytype{\alpha_4}, \raytype{\alpha_3}, \raytype{\alpha_3}. These are shown in
Figure~\ref{fig:tropbasis}, and -- more colorfully -- in Figure~\ref{fig:atlas_rays}.
The $4+6$ triangles in the second row of Figure \ref{fig:tetrahedra} are $\Ggroup$-trees of
types \celltype{32}, \celltype{32}, \celltype{33}, \celltype{34} for $p \not=2$ and of types
\celltype{18}, \celltype{18}, \celltype{19}, \celltype{20}, \celltype{35}, \celltype{35} for $p=2$.
Representatives for these orbits are shown in
Figures~\ref{fig:atlas_toric_cells_prime} and~\ref{fig:atlas_toric_cells_nontoric}.

\begin{proof}[Proof of Theorem~\ref{thm:char2}]
Let $I_\KK$ be the ideal in $\KK[p]$ that is generated by the $35$ Pl\"ucker quadrics and the seven linear forms in (\ref{eq:7linear}).
Since these generators have coefficients $\pm 1$, the tropical variety $\trop_p(\linespace)$ 
is contained in the $\Ggroup$-Dressian for every $p$.
For odd primes $p$, the proof of Corollary~\ref{cor:tropicalbasis} applies verbatim, as $2 \neq 0$ in $\KK$.
This shows the inclusion $\trop_p(\linespace) \subseteq \trop(\linespace)$.

For the reverse inclusion, we use the following criterion: if $C$ is a cone and $G$ a Gr\"obner basis of $I_\KK$
for an order refining some $w \in \relint(C)$, with $\initial_v(g) = \initial_w(g)$ for all $g \in G$ and $v \in \relint(C)$,
and the $\initial_w(g)$ have a common zero in the torus $(\overline{\KK}^{\,*})^{21}$, then $C \subseteq \trop_p(\linespace)$.

We apply this criterion to the $251$ maximal cones $C$ of the simplicial fan structure on $\trop(\linespace)$ in the proof of
Theorem~\ref{thm:grassmannian}, with $G$ a Gr\"obner basis over $\ZZ$, as in \cite[\S 4]{AdamsLoustaunau:1994}.
This was computed in \oscar.
We verify that $\initial_v(g) = \initial_w(g)$ on $\relint(C)$, and that the leading coefficients 
appearing in $G$ are $1$ or $2$.
They are $1$ unless $C$ is a facet of a tetrahedron.
For $p$ not dividing them, $G$ reduces modulo $p$ to a Gr\"obner basis of $I_\KK$ with the same initial forms.
These have a common zero in $\{\pm 1\}^{21}$ for $237$ cones.
For the $12$ facets of the tetrahedra and the two cones of type \celltype{21}, we verify that the ideal generated by
the $\initial_w(g)$, by $p_{ij} - 1$ for six coordinates $ij \in S$, and by $y \prod_{i'j' \notin S} p_{i'j'} - 1$ has a Gr\"obner basis over $\ZZ$
without constant whose leading coefficients are powers of $2$, all $1$ for type \celltype{21};
for $p$ not dividing them, this ideal is proper over $\KK$, and its zeros are torus points.
Hence $C \subseteq \trop_p(\linespace)$ for odd $p$, which proves the first assertion.
This also works for $p = 2$ unless $C$ is a facet of a tetrahedron.

Now let $p = 2$. Since $\trop_2(\linespace)$ is contained in the Dressian, it remains to examine the three tetrahedra.
Modulo $2$, the cubic $f$ in (\ref{eq:threecubics}) of a tetrahedron 
$T$  has four terms, of equal weight at the $\Ggroup$-tree \raytype{\epsilon} of $T$;
at each vertex of $T$, one term is strictly smaller than the other three, which are equal, and each term is the smaller one at exactly one vertex.
Hence $\initial_w(f)$ is a monomial for all $w \in T$ outside the six triangles spanned by \raytype{\epsilon} and an edge of $T$,
and $\trop_2(\linespace) \cap T$ lies in these triangles. Conversely, each of these $18$ triangles satisfies the criterion,
with a Gr\"obner basis of $I_{\FF_2}$ and the common zero $(1,\ldots,1)$.
Thus $\trop_2(\linespace)$ consists of the $239$ cones above and the $18$ triangles. Its $f$-vector, orbits and homology
are computed as in Theorem~\ref{thm:grassmannian}, after merging three pairs of triangles into quadrilaterals.
\end{proof}

An analogous result holds for the flag variety $\flagspace$.
We state this in a second theorem.

\begin{theorem} \label{thm:flagchar2}
  If $p = \characteristic(\KK) $ is odd  then $\trop_p(\flagspace)= \trop(\flagspace)$.
  If $p=2$ then, modulo~lineality, the tropical $\Ggroup$ flag 
  variety $\trop_2(\flagspace)$ is a $3$-dimensional polyhedral complex  with $f$-vector $(156, 1445, 3758,2792)$.   
  Its $3$-cells are $2630$ tetrahedra, $144$ square pyramids~and $18$ triangular prisms.
  The reduced integral homology is $\ZZ^{324}$, concentrated  in dimension $3$.
  Modulo  $\Dtwelve$-action,  $\,\trop_2(\flagspace)$ has $f$-vector $(31, 192, 416, 281)$.
  Each cell is a decorated $\Ggroup$-tree.
\end{theorem}

\begin{proof}
We proceed as in the proof of Theorem~\ref{thm:char2}. Let $p \geq 3$.
Since  the generators of $J$ have coefficients $\pm 1$,  $\trop_p(\flagspace)$ lies in the flag Dressian.
Outside $\trop(\flagspace)$, it consists of the $33$ open $4$-simplices and $15$ of their facets.
On each $4$-simplex, the initial form of a cubic in (\ref{eq:threecubics}) is its term with coefficient~$2$.
So $\trop_p(\flagspace) \subseteq \trop(\flagspace)$: the generators and the three cubics form a tropical basis for $p \neq 2$.
Conversely, we apply the criterion to one maximal cone in each of the $270$ $\Dtwelve$-orbits, as $\Dtwelve$ acts by signed permutations preserving $J$ over $\ZZ$,
or to the Gr\"obner cones in it if the initial forms vary ($44$ orbits).
Leading coefficients are $1$, $2$ or $4$; torus points are $\pm 1$-points or slices.
For $p=2$ we apply Algorithm~\ref{alg:balanced_support} over~$\FF_2$ to the $78$ generators and $12$ elements of a tropical basis of $\linespace$ over $\FF_2$.
Their prevariety's~$3038$ cones merge into $2792$ of dimension $8$ in $281$ orbits, all kept.
We locate the output in the flag Dressian of Theorem~\ref{thm:flag_dressian}.
Over each tetrahedron, a decorated version of
Figure \ref{fig:tetrahedra} emerges.
\end{proof}

\section{The $\mathbf{\Ggroup}$-Forest}
\label{sec5}

\begin{figure}[ht]
  \centering
  \figurepdf{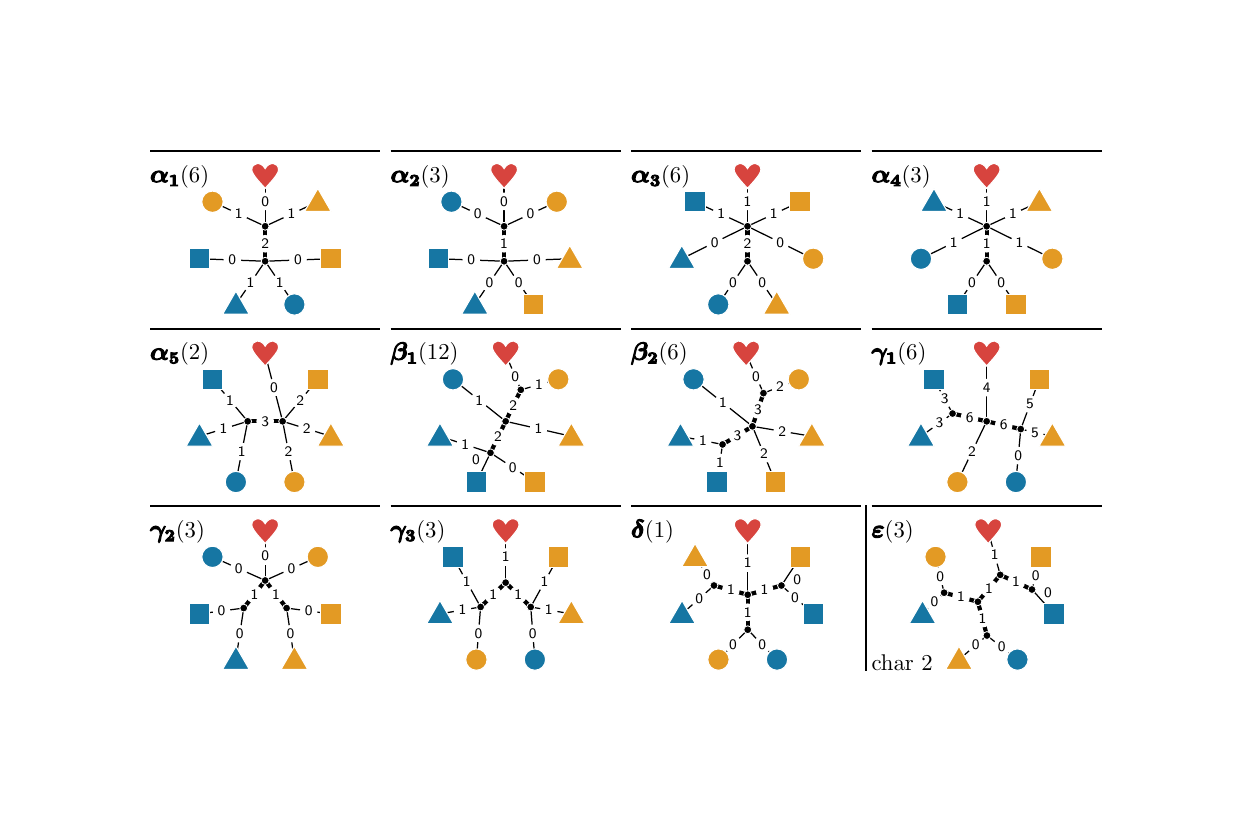}
  \caption{The $12$ $\Ggroup$-trees for the rays of the tropical Grassmannian and Dressian.}
  \label{fig:atlas_rays}
\end{figure}

\begin{figure}[tbp]
  \centering
  \figurepdf{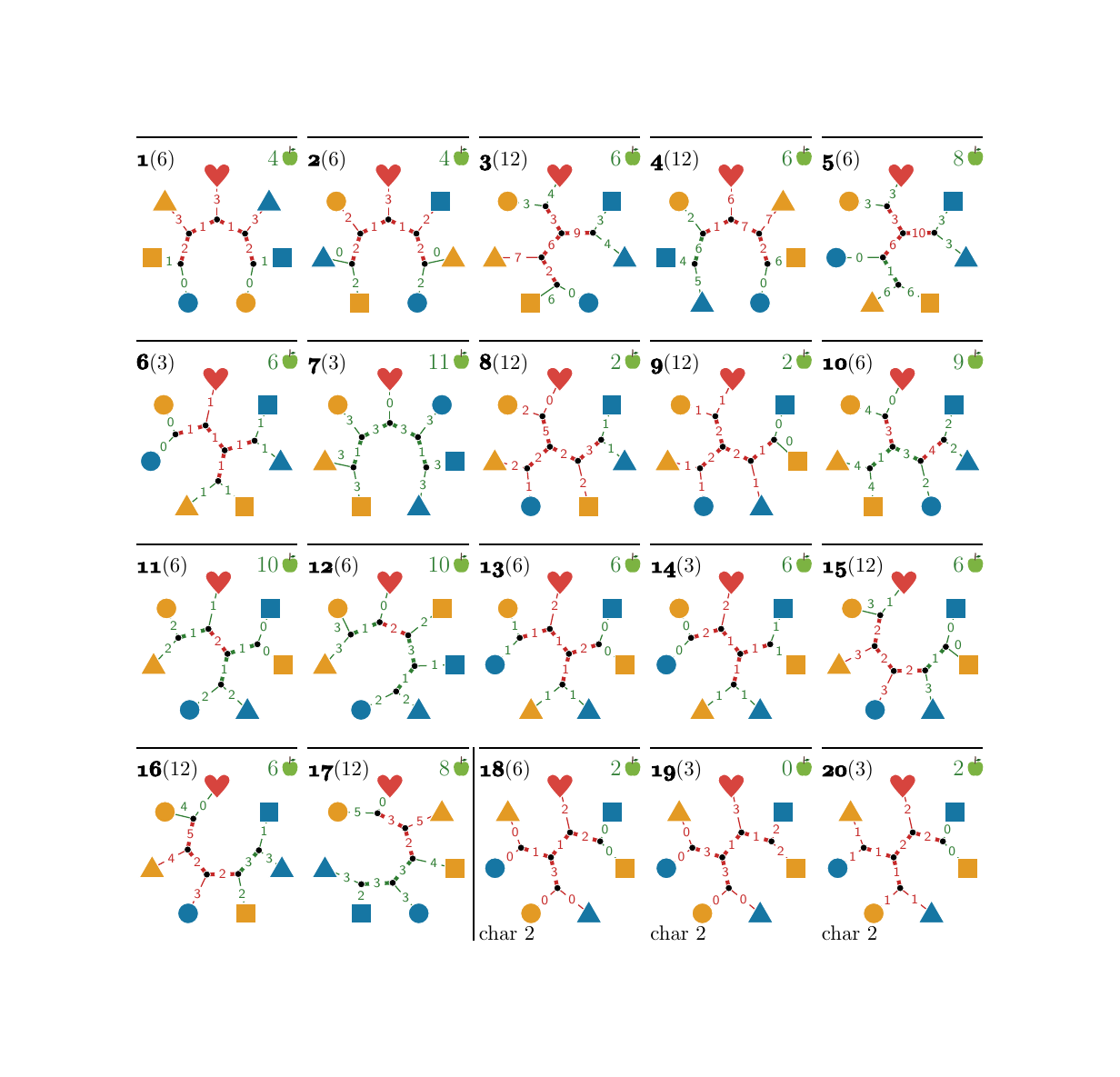}
  \caption{The $20$ types of maximal $\Ggroup$-trees that give
  toric degenerations in Sections \ref{sec6} and~\ref{sec7}.}
  \label{fig:atlas_toric_cells_prime}
\end{figure}

\begin{figure}[tbp]
  \centering
  \figurepdf{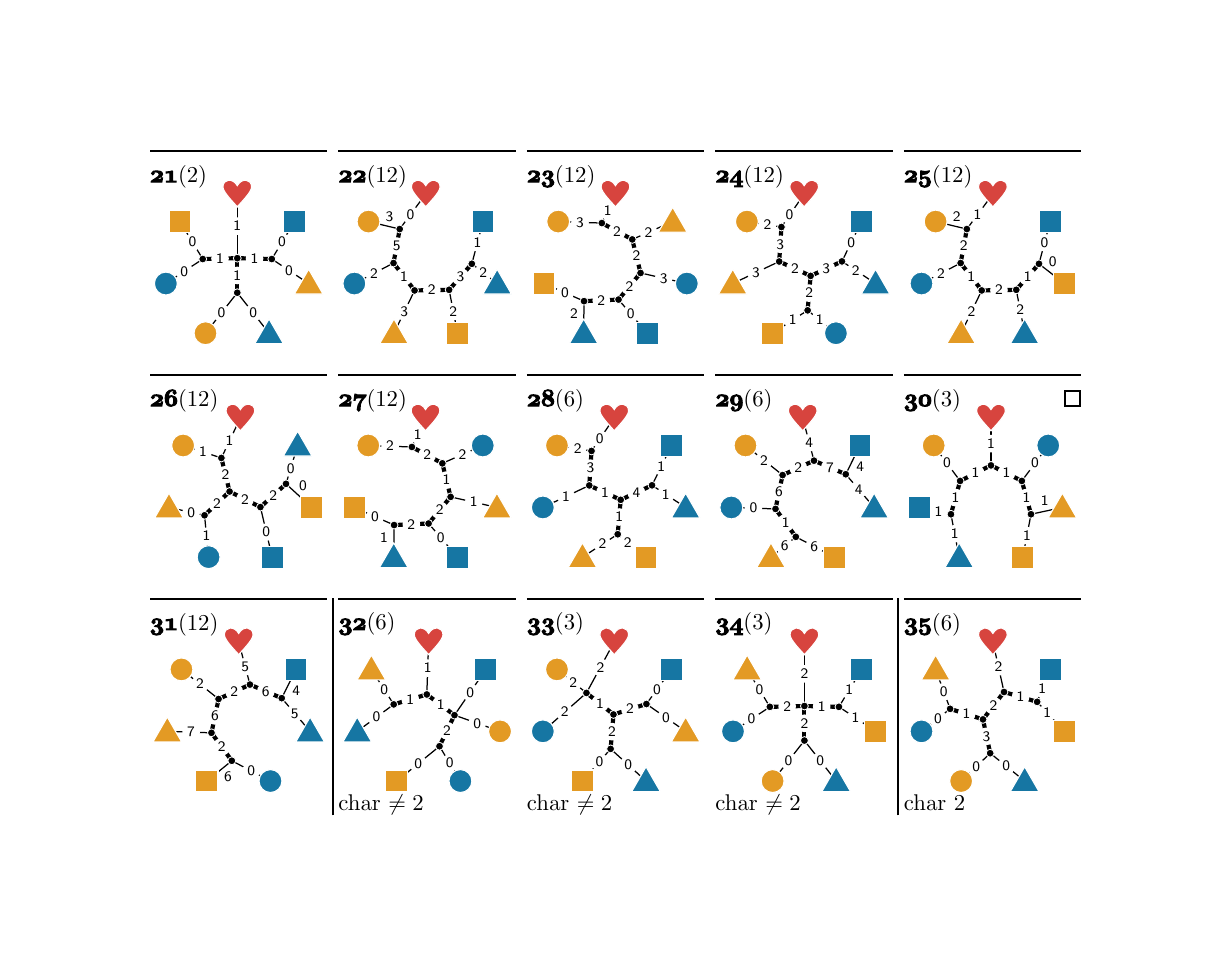}
  \caption{The $15$ types of maximal $\Ggroup$-trees that do not give toric degenerations.}
  \label{fig:atlas_toric_cells_nontoric}
\end{figure}

The $\Ggroup$-forest is composed of all $\Ggroup$-trees that correspond to rays and maximal cells of our $\Ggroup$-Grassmannians.
Here, we include both the general version $\trop(\linespace) = \trop_p(\linespace)$, $p \geq 3$,
and the characteristic $2$ version $\trop_2(\linespace)$.
These are described in Theorems  \ref{thm:grassmannian} and~\ref{thm:char2}.
The $\Dtwelve$-orbits of the $\Ggroup$-trees are shown in
 Figures~\ref{fig:atlas_rays},~\ref{fig:atlas_toric_cells_prime}, and~\ref{fig:atlas_toric_cells_nontoric}.
The diagrams are labeled by symbols in boldface, either Greek letters or numbers,
followed by the orbit size in parentheses.
The edge coloring and the green apples seen in Figure \ref{fig:atlas_toric_cells_prime} will be explained in Section~\ref{sec7}.

There are $12$ types of rays, $11$ for $\characteristic(\KK)\not=2$ 
plus one for $\characteristic(\KK) = 2$.
Their $\Ggroup$-trees are shown in Figure~\ref{fig:atlas_rays}.
The type \raytype{\alpha} represents splits, while \raytype{\beta} or \raytype{\gamma} represent caterpillars, with
the heart either at an end or in the middle.
Finally, \raytype{\delta} and \raytype{\epsilon} represent snowflakes, where \raytype{\delta} is fixed under $\Dtwelve$.
Note that \raytype{\epsilon} appears only in characteristic $2$; see Figures~\ref{fig:tropbasis} and~\ref{fig:tetrahedra}.
  
We group the maximal cells into types \celltype{1}--\celltype{20} that yield a toric degeneration, in Figure~\ref{fig:atlas_toric_cells_prime}, and types \celltype{21}--\celltype{35} that do not, in Figure~\ref{fig:atlas_toric_cells_nontoric}. 
Toric degenerations will be explained in~Section \ref{sec6},
where Table~\ref{tab:g2-grassmannian-newton-okounkov}
matches Figure~\ref{fig:atlas_toric_cells_prime}.
Types \celltype{18}--\celltype{20} and \celltype{35} appear only in characteristic $2$ while types \celltype{32}--\celltype{34} disappear in characteristic $2$. Type \celltype{30} is the orbit of the three quadrilaterals.

Each orbit size is easily found by looking at the $\Ggroup$-trees under the $\Dtwelve$-action in Figure~\ref{fig:shape_hexagon}. For instance, the split type \raytype{\alpha_2} has an orbit of size $3$: its three $\Ggroup$-trees correspond to grouping the heart with two same shapes and splitting them from the other four leaves.

Table~\ref{tab:cell_ray_orbits} lists which type of rays generates each type of maximal cells.
This table also helps to identify pairs of $\Ggroup$-trees that form edges of our $\Ggroup$-Grassmannian: those are generated by two rays of this table.
One necessary  (but not sufficient) condition for several rays to generate a cell is for their 
associated trees to have compatible splits.
For instance, \raytype{\delta} and \raytype{\alpha_4} are compatible, and they form an edge of the triangles \celltype{13} and \celltype{14},
as listed in Table~\ref{tab:cell_ray_orbits}.

\begin{table}[htbp]
  \centering
  \caption{The ray orbits of the maximal cells in Figures~\ref{fig:atlas_toric_cells_prime}
  and~\ref{fig:atlas_toric_cells_nontoric}, labeled as in Figure~\ref{fig:atlas_rays}.}
    \label{tab:cell_ray_orbits}
  \input{tables/cell_ray_orbits}
\end{table}

\begin{figure}[tbp]
  \centering
  \figurepdf{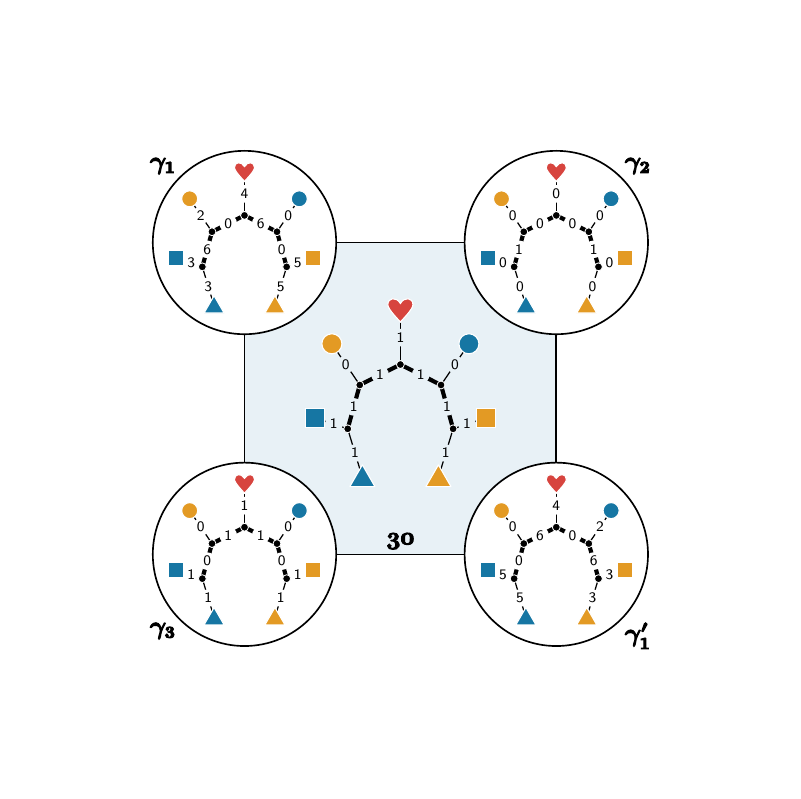}
  \caption{Detailed description of type \celltype{30}. This maximal cell of $\trop(\linespace)$ is a quadrilateral.}
  \label{fig:square_cell}
\end{figure}

Figure~\ref{fig:square_cell} describes all $4+4+1$ $\Ggroup$-trees on a quadrilateral (type \celltype{30}).
Two vertices of the quadrilateral belong to the same orbit \raytype{\gamma_1}.
Inner edges of length zero were artificially added  so the support trees agree.
The $\Ggroup$-tree $d$ on the quadrilateral 
 is normalized and expressed as $\frac 16(\raytype{\gamma_1}+\raytype{\gamma_1'} - 2 s) =\raytype{\gamma_2}+\raytype{\gamma_3}$, where $s$ is the star tree with all pendant edges of lengths $1$.
To be completely explicit, we have
$d_{12}{=}d_{16}{=} d_{26}{=} d_{34}{=} d_{35}{=} d_{36}{=} d_{37}{=} d_{45}{=} d_{67}  = 2$,
$\,d_{13}{=}d_{17}{=} d_{23}{=} d_{27}{=} d_{46}{=} d_{47}{=} d_{56}{=} d_{57} = 4$, and
 $d_{14}{=}d_{15}{=}d_{24}{=}d_{25}  = 6$ for the tree metric $d$.

\section{Toric Degenerations}\label{sec:toric_degenerations}
\label{sec6}

Toric degeneration is a  process which replaces a given projective variety by a toric variety.
This process retains key invariants of the original variety, such as the dimension
and the degree, and it offers new insights by means of polyhedral geometry  and combinatorics.
Associated with the toric variety is a lattice polytope, called the {\em Newton--Okounkov polytope}.
The lattice points in this polytope encode elements in a {\em Khovanskii basis} for the given coordinate ring,
and this in turn offers valuable algebraic information. This approach has
been especially fruitful in representation theory, where one studies
varieties with a group action, like Grassmannians and flag varieties.
In what follows, we explore this for our group $\Ggroup$.

Kaveh and Manon \cite{KM:2019} proved that the toric degenerations correspond to 
the maximal cones in the tropical variety whose associated initial ideal is prime.
From this data, one can then read off the toric variety and the Newton--Okounkov polytope.
We shall explain this process for the $\Ggroup$-Grassmannian, 
thereby classifying all toric degenerations of $\linespace$ inside $\PP^{20}$.
Our description is
for fields of every characteristic, and hence for both $\trop(\linespace)$ and $\trop_2(\linespace)$.

Let $C$ be a maximal cone of $\trop(\linespace)$. We call $C$ a {\em prime cone} if
the initial ideal $\initial_C(\linespace)$ is prime over the algebraic closure $\overline{\KK}$ of~$\KK$. Here we mean
the ideal $\initial_d(\linespace)$ for any $d \in \relint(C)$. 
If $C$ is a prime cone then $\initial_C(\linespace)$ is 
a {\em toric ideal}, and we also write $\initial_C(\linespace)$ for the corresponding
projective toric variety in $\PP^{20}$. The Newton--Okounkov polytope $P_C$,
which is associated to this toric variety, has dimension $5 = \dim(\linespace)$.
Its normalized volume is $18 = \degree(\linespace)$.
We compute $P_C$ as follows. Consider any
 $5 \times 21$ matrix whose rows form a $\ZZ$-basis for the
lattice points in the linear span of $C$.
Then $P_C$ is the convex hull in $\RR^5$ of the $21$ columns.

\begin{example}[Type \celltype{1}]
The first tree in Figure \ref{fig:atlas_toric_cells_prime} equals the rightmost tree in
Figure \ref{fig:ex_G2_trees}. Its metric $d$ is shown explicitly in
Example \ref{ex:fivetrees}. The initial ideal $\initial_d(\linespace)$ is toric. Its generators~are
\[ 
{\footnotesize \begin{matrix}
p_{37}-p_{45}, \,\, p_{12}+p_{67}, \,\, p_{27}+p_{46}, \,\, p_{23}+p_{47}, \,\, p_{17}-p_{56}, \,\, p_{13}-p_{57}, \,\, p_{14}+p_{25}, \,\, \\
p_{26} p_{57} + p_{36} p_{67}, \,\, p_{36} p_{45} - p_{35} p_{46}, \,\,
p_{45} p_{56} - p_{36} p_{57}, \,\, p_{36} p_{56} + p_{57} p_{67}, \,\, \\
p_{45} p_{46} - p_{36} p_{47}, \,\, p_{36} p_{46} + p_{47} p_{67}, \,\, 
p_{25} p_{35} + p_{45} p_{57}, \,\, p_{24} p_{35} + p_{45} p_{47}, \,\, \\
p_{25} p_{26} - p_{46} p_{67}, \,\, p_{15} p_{26} + p_{56} p_{67}, \,\, 
p_{47} p_{56} - p_{46} p_{57}, \,\, p_{46} p_{56} - p_{25} p_{67}, \,\, \\
p_{25} p_{45} + p_{47} p_{57}, \,\, p_{25} p_{36} + p_{46} p_{57}, \,\, 
p_{24} p_{36} + p_{46} p_{47}, \,\, p_{15} p_{36} - p_{56} p_{57}, \,\, \\
p_{25} p_{47} - p_{24} p_{57}, \,\, p_{15} p_{47} + p_{25} p_{57}, \,\,
p_{25} p_{46} - p_{24} p_{56}, \,\, p_{15} p_{46} + p_{25} p_{56}, \,\, \\
 p_{25}^2+p_{15} p_{24}, \,\,
 p_{36}^2+p_{45} p_{67}, \,\,
 p_{45}^2-p_{35} p_{47}, \,\,
 p_{46}^2-p_{24} p_{67}, \,\, \\
 p_{47}^2+p_{24} p_{45}, \,\,
 p_{56}^2+p_{15} p_{67}, \,\,
 p_{57}^2-p_{15} p_{45}, \,\,
 p_{67}^2-p_{26} p_{56} .
 \end{matrix}}
\]
The $7$ linear binomials say that $P_C$ is the convex hull of
$14 = 21-7$ distinct lattice points. The last $8$ quadrics show
that only $6 = 14-8$ of these are vertices. Hence $P_C$ is a $5$-simplex.
See the first row of Table  \ref{tab:g2-grassmannian-newton-okounkov}.
This is the {\em Chevalley polytope} studied by   Spacek and Wang~\cite{SW:2024}.
\end{example}

\begin{example}[Type \celltype{21}] \label{ex:twentyone}
The first tree in Figure \ref{fig:atlas_toric_cells_nontoric}
equals that on the triangle in Figure~\ref{fig:triangle_cell}.
We see that $d _{15} = d_{26} = d_{34} =0$
and $d_{ij} = 2$ for all other $i,j$. The initial
ideal $\initial_d( \linespace)$ is prime over $\QQ$, 
so \oscar calls this a prime ideal.
However, the ideal is not prime over $\CC$. To see this, we check that 
$(p_{12}^2 p_{13})^2 + (p_{12}^2 p_{13})( p_{14}^2 p_{16})
+ (p_{14}^2 p_{16})^2$ lies in $\initial_d( \linespace)$.
This polynomial factors over $\CC$, but neither factor
lies in $\initial_d( \linespace)$. In particular, 
$\initial_d( \linespace)$ is not a toric ideal.
We conclude the cones of type \celltype{21} are not prime cones
of the tropical Grassmannian $\trop(\linespace)$.
\end{example}

We compute the initial ideal for every maximal cone $C$ of $\trop(\linespace)$, 
we test whether it is toric, and -- if yes -- we compute $P_C$ as described above.
We then redo the same computation for the flag variety $\flagspace$.
Using this approach, we obtain the following classification result.

\begin{theorem} \label{thm:NObody}
  In any characteristic, the $\Ggroup$-Grassmannian has $135$ toric degenerations in $17$ orbits.
  In characteristic $2$, there are  $12$ further toric degenerations in $3$ orbits.
  All $20$ orbits are listed in Table~\ref{tab:g2-grassmannian-newton-okounkov}.
  In any characteristic, the $\Ggroup$ flag variety has $807$ toric degenerations in $81$ orbits.
  In characteristic $2$, there are $18$ further toric degenerations in $2$ orbits.
\end{theorem}

\begin{table}[htbp]
  \centering
    \caption{Newton--Okounkov polytopes for all $20$ toric
  degenerations of the $\Ggroup$-Grassmannian over $\QQ$ and over $\FF_2$. 
The orbit sizes in
Figure \ref{fig:atlas_toric_cells_prime} sum to $147$.
The last two columns give the number of lattice points $N$ in the polytope
and the field characteristic. See~\cite{BBS:2002} for multiplexes.}
  \label{tab:g2-grassmannian-newton-okounkov}
  \small
  \input{tables/toric_degeneration_table.tex}

\end{table}

\begin{remark}
A toric degeneration fails to be normal when $|P \cap \ZZ^5|$ is $15$ or $16$:
a lattice point of $P$ that is not among the $14$ points of the generators gives a monomial
in the integral closure of the toric ring that is not in the ring.
For instance, at a cell of type \celltype{6}, the toric ideal contains $p_{45} p_{56}^2-p_{57}^2 p_{67}$.
So $u = p_{57} p_{67}/p_{56}$ is integral, as $u^2 = p_{45} p_{67}$, but it is not in the toric ring:
its lattice point, the midpoint of those of $p_{45}$ and $p_{67}$, is not that of any $p_{ij}$.
\end{remark}

Several of the Newton--Okounkov polytopes 
in Table \ref{tab:g2-grassmannian-newton-okounkov}
 already appeared in the literature.
Type \celltype{1} is the Chevalley polytope of Spacek and Wang in~\cite[Section A.4]{SW:2024}.
Christophersen and Ilten \cite[Theorem 3.1]{CI:2014} degenerate $\linespace$
to the Stanley--Reisner schemes $\PP(K\ast\Delta_2)$ of simplicial
$2$-spheres $K$ with $11$ vertices, one for each of their term orders
\cite[Fig.~3]{CI:2014}. For two of them $K$ is the sphere $T_{11}$ in
\cite[Fig.~1]{CI:2014}. Each of these degenerations factors through our toric
degenerations among types \celltype{7}, \celltype{10}, \celltype{11}, \celltype{12} and \celltype{17}:
for suitable cells of these types. The canonical weight lies in the closed
Gr\"obner cone of their term order. Hence these Newton--Okounkov polytopes
admit a regular unimodular triangulation of the form $T_{11}\ast\Delta_2$.
In particular, they are normal. Cutting by general hyperplanes
\cite[Corollary 3.3]{CI:2014} degenerates the Fano threefold $V_{18}$ to
$\PP(T_{11}\ast\Delta_0)$ and the K3 surface of genus $10$ to $\PP(T_{11})$,
a connection which is important for our Mukai curves in Section~\ref{sec10}.
One of these toric degenerations, namely type \celltype{7}, also appeared in
work of Feigin, Fourier and Littelmann \cite[Definition 11.5]{FeiginFourierLittelmann:2017}.
 
\section{The $\mathbf{\Ggroup}$-Orchard}
\label{sec7}

\begin{figure}[htbp]
  \centering
  \figurepdf{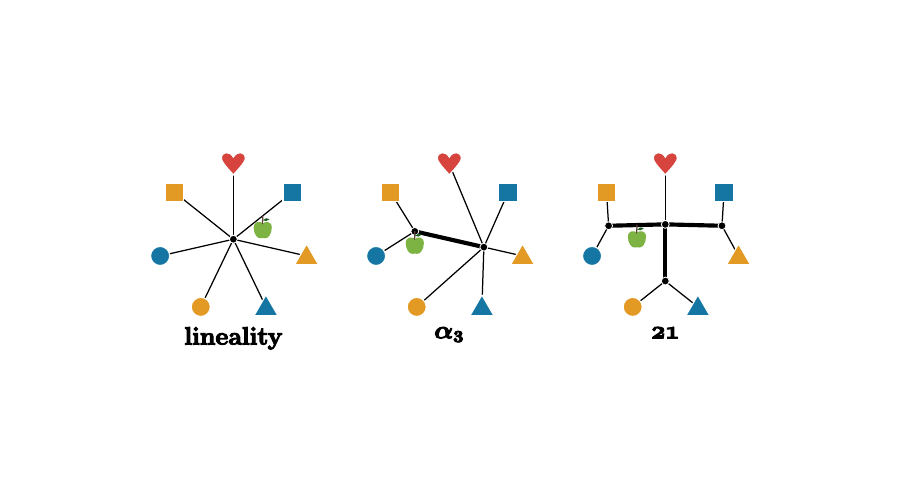}
  \caption{Two types of vertices and a $3$-cell of $\trop(\flagspace)$.
  }
  \label{fig:decorated_trees}
\end{figure}

The $\Ggroup$-orchard consists of $\Ggroup$-trees decorated by a fruit, 
namely the green apple in Figure~\ref{fig:decorated_trees}.
This represents the leaf $8$ discussed prior to
Theorem \ref{thm:flag_dressian}. The distance from the apple to leaf $k$ equals $d_{k8} = - \val(x_k)$.
Every decorated $\Ggroup$-tree is a point of the $\Ggroup$ flag Dressian.

The green apples explain Theorem \ref{thm:flag}.
They can grow on an edge or on an inner vertex of a $\Ggroup$-tree.
The $141$ vertices of $\trop (\flagspace)$ are 
$7$ star trees with an apple on an edge and $134$ decorated $\Ggroup$-trees with an apple on a node
in one of the $11$ trees in Figure~\ref{fig:atlas_rays}.
Two of them are on the left in Figure \ref{fig:decorated_trees}.
The $2714$ maximal cells of $\trop(\flagspace)$ are the edge-decorated $\Ggroup$-trees of type 
\celltype{1}--\celltype{35} (Figures~\ref{fig:atlas_toric_cells_prime} and~\ref{fig:atlas_toric_cells_nontoric}). 
One of them is shown on the right in Figure~\ref{fig:decorated_trees}.

This gives a way of identifying each decorated $\Ggroup$-tree, but it does not say much about~its surroundings
in $\trop(\flagspace)$.
We next provide  charts for guiding the reader through the orchard.
This rests on the projection $\trop(\flagspace) \rightarrow \trop(\linespace)$,
which comes from picking the~apple.

The fiber over the lineality space of $\trop(\linespace)$ is given by the $7$ decorated star trees.
These appear over each cell of $\trop(\linespace)$.
Fix a maximal (i.e., $2$-dimensional) cell $C$ of $\trop(\linespace)$.
Let $\flagspace_C$ denote the set of all cells in $\trop(\flagspace)$
that map to~$C$. This is a full-dimensional subcomplex of $\trop(\flagspace)$.
We call it the {\em fiber complex} over~$C$.
The tropical $\Ggroup$ flag variety $\trop(\flagspace)$ is the union of
the fiber complexes $\flagspace_C$, where $C$ runs over the maximal cells of $\trop(\linespace)$.

The cells of $\flagspace_C$ are all possible ways of decorating the $\Ggroup$-trees given by the faces of $C$. The fiber complex $\flagspace_C$ has either $10$ or $11$ maximal cells, of dimension $3$.
Indeed, an apple can grow from any of the edges of the $\Ggroup$-tree of $C$.
The $3$-cells are glued along facets.
These $2$-cells correspond to apples growing on a node. They all biject to $C$.
Every fiber complex is attached to the $7$ vertices corresponding to the edge-decorated star trees.
The remaining vertices are obtained by decorating the nodes of the $\Ggroup$-trees associated to the vertices of $C$.

The following description of the $\Ggroup$-orchard was obtained by
direct computations.

\begin{proposition}
There are $23$  combinatorial types of fiber complexes $\flagspace_C$
over the trees $C$ in Figures \ref{fig:atlas_toric_cells_prime} and
\ref{fig:atlas_toric_cells_nontoric}.  These are
 listed in Table~\ref{tab:local_flag_complex}.
Six distinct $f$-vectors are labeled by
$\mathfrak{A}, \ldots, \mathfrak{F}$.
\end{proposition}

\begin{table}[th]
    \caption{The $23$ combinatorial types of fiber complexes $\flagspace_C$.
  These cover $\trop(\flagspace)$.}
  \label{tab:local_flag_complex}
  \centering
\input{tables/local_flag_complexes.tex}
\end{table}

\begin{figure}[tbp]
  \centering
  \figurepdf{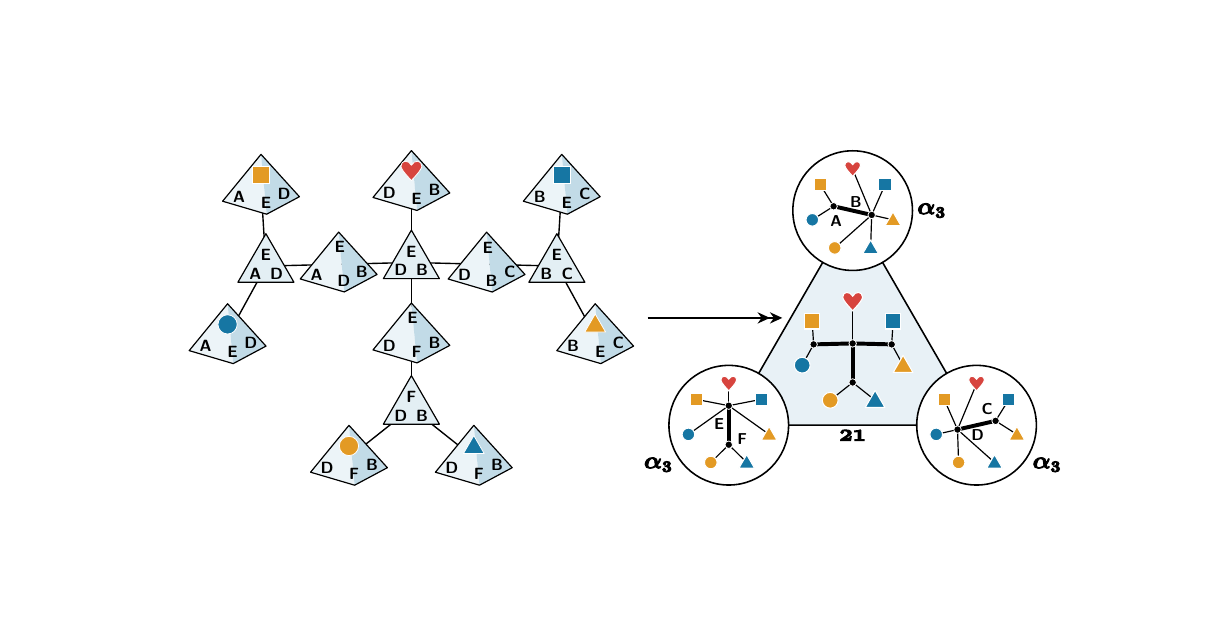}
  \caption{The polyhedral complex of type $\mathfrak A_1$ projecting to the triangle of type \celltype{21}.}
  \label{fig:triangle_cell}
\end{figure}

We now describe the $23$ fiber complexes in more detail. Pictures
will be helpful to appreciate the beauty of the $\Ggroup$-orchard.
For instance,  Figure \ref{fig:triangle_cell} shows the complex $\flagspace_{\mathfrak{A}_1}$. 
Here $\mathfrak{A}_1$ refers to the first line in Table \ref{tab:local_flag_complex}.
The triangle-tetrahedron incidence graph of this 
$3$-dimensional complex is the inner-edge subdivision of the $\Ggroup$-tree it projects to.
We see this in Figure~\ref{fig:triangle_cell}: each edge of the tree gives a tetrahedron, and each inner node a triangle.

The vertices of each fiber complex are labeled by our $7$ shapes in 
Figure \ref{fig:shape_hexagon} for the edge-decorated star trees, 
and by letters $A,B,\dots$ for internal nodes. For instance, in
Figure \ref{fig:triangle_cell}, we
have six nodes $A,B,C,D,E,F$
in the three trees appearing in the triangle of $\trop(\linespace)$
underlying $\mathfrak{A}_1$.
The maximal cells of $\flagspace_{\mathfrak{A}_1}$ are tetrahedra,
specified by the $4$ vertex labels. Likewise, triangles are given by $3$
vertex labels.
For instance, the three tetrahedra $ABDE$, $BCDE$ and $BDEF$
share the triangle $BDE$.
  The tetrahedron over an edge $e$ has as vertices the images of the endpoints of $e$
under the relevant contractions of the tree. For instance, for the tetrahedron $  ABDE$,
one contracts the edge $AB$ and retains the nodes $D$ and $E$.
For each pendant edge, one of the vertices of the tetrahedron is the
associated  colorful shape.

\begin{figure}[tbp]
  \centering
  \figurepdf{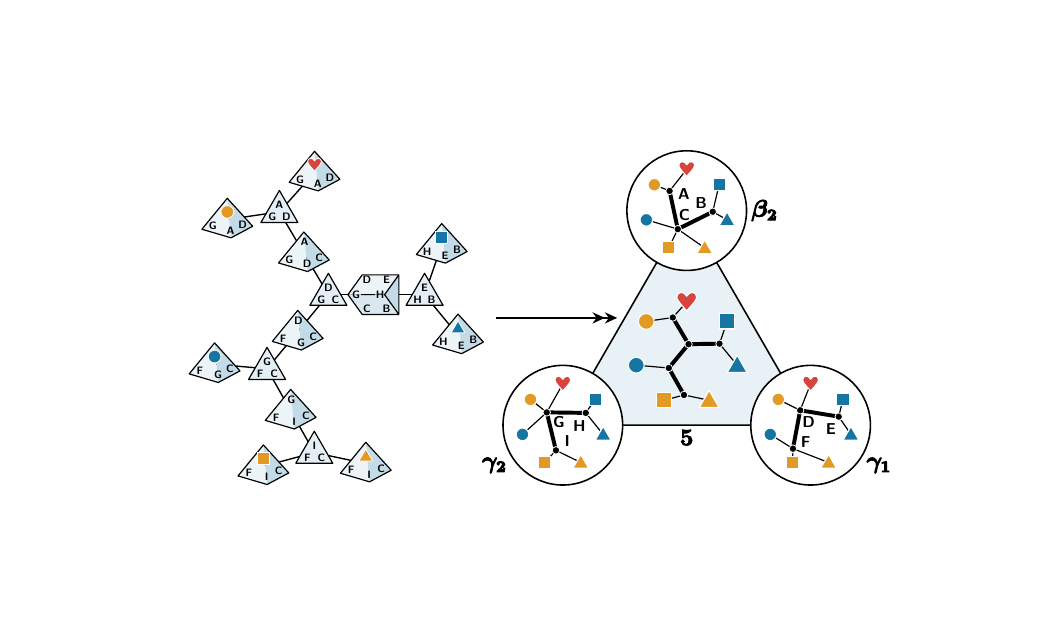}
  \caption{The polyhedral complex of type $\mathfrak D$ projecting to the triangle of type \celltype{5}.}
  \label{fig:prism_projection}
\end{figure}

In Figure~\ref{fig:prism_projection}, we present the fiber complex 
$\flagspace_{\mathfrak D}$ over the triangle of type \celltype{5} in Figure \ref{fig:atlas_toric_cells_prime}.
The fiber complex $\flagspace_{\mathfrak D}$ is not a simplicial complex. One of its $3$-cells
is a triangular prism, labeled by $BCDEGH$.
This appears because the blue square and the blue triangle  form a cherry 
in each of the three vertex $\Ggroup$-trees.
Hence the corresponding edges $BC$, $DE$ and $GH$ are consistent.
The two facets $CDG$ and $BEH$ of the triangular prism
are also in $\flagspace_{\mathfrak D}$.

Table \ref{tab:local_flag_complex} summarizes
the fiber complexes for all types \celltype{1}--\celltype{35}.
Only the first two rows yield simplicial complexes.
Pyramids and prisms appear in the other fiber complexes.


We close this section with a discussion of the toric degenerations of the
$\Ggroup$ flag variety~$\flagspace$.  Each of these induces a toric degeneration of
the $\Ggroup$-Grassmannian $\linespace$.
We know from  Theorem~\ref{thm:NObody} that $\trop(\flagspace)$ has
$807$ prime cones. Each of these $807$ projects onto one of the $135$ prime cones in
$\trop(\linespace)$. Hence, the $81$ types of toric degenerations
of $\flagspace$ map to the $17$ types of toric degenerations of $\linespace$.
There is an analogous correspondence in characteristic~$2$.

We examined this map for the types $\celltype{1}$--$\celltype{20}$.
For each $3$-cell $C$ of each fiber complex in Table~\ref{tab:local_flag_complex}, we computed the initial ideal
$\initial_C(\flagspace)$, and we recorded whether it is prime over $\overline{\QQ}$.
The resulting prime cells are marked by green edges in 
Figure~\ref{fig:atlas_toric_cells_prime}. The green number next to the green apple 
records the number of green edges. By multiplying each green number
with the parenthesized orbit size, and by adding up all $19$ products, 
we arrive at the total count:
\[ 807 \,+\, 18 \,\,= \,\, (6 \times 4\,+\,6 \times 4 \,+ \,12 \times 6 \,+ \cdots + 3 \times 11 \,+ \cdots  +12 \times 8)\, +\, (6\times 2\,+\, 3\times 2). \]
The edges that are colored in red in Figure~\ref{fig:atlas_toric_cells_prime}
do not give any toric degeneration of $\flagspace$.
We note that the characteristic $2$
  type \celltype{19} has no prime cone in its fiber complex
 inside $\trop_2(\flagspace)$.

\begin{remark}
Our discussion matches what is in the literature on $\Ggroup$.
The toric degeneration of $\flagspace$ due to Spacek and Wang~\cite{SW:2024} is obtained by decorating either of the two circle leaves of type \celltype{1}. The toric degeneration given by Feigin, Fourier and Littelmann~\cite[Definition 11.5]{FeiginFourierLittelmann:2017} is obtained by decorating one of the two edges supporting the cherries in type \celltype{7}.
\end{remark}

\section{Tropicalizing Homogeneous Spaces}
\label{sec:homogeneous}

Tropicalization depends on coordinate choices.
What does this mean for projective varieties that are homogeneous spaces?
When the representation of interest is minuscule~\cite{CDFS:2025,GKL:2025,Corey:2026},
a coordinate system is handed to us. This was the case when the
tropical Grassmannian first arose. The authors of \cite{speyer_sturmfels}
did not think about choices: they simply used Pl\"ucker coordinates; see also \cite{BO:2023} for other Grassmannians.
On the other hand, simple (or, more generally, reductive) groups can be tropicalized without coordinates, through their buildings~\cite{CU:2025,GanYu:2003}.
The present paper is one of the first to address this dichotomy as a practical issue.

The group $\Ggroup$ has three natural projective homogeneous spaces, namely the quotients by its two maximal parabolic subgroups
$P_1$ and $P_2$, and by its Borel subgroup $B=P_1\cap P_2$.
Tropicalizing $\Ggroup/P_1\cong\pointspace$ is easy; the weights of $\KK^7$ have multiplicity one.
However, for $\Ggroup/P_2\cong\linespace$, the representation has a zero weight plane; so a choice of coordinates seems necessary.
Coordinates for both then give coordinates for $\Ggroup/B\cong\flagspace\subset\pointspace\times\linespace$.

Forests and orchards grew in Sections~\ref{sec2}--\ref{sec7}.
Their trees came from embedding $\linespace$ into $\grassmannian(2,7)$.
The group only entered through the linear forms in~(\ref{eq:7linear}), the grading of Figure~\ref{fig:root_lattice},
the star trees that span the lineality space, and the Weyl group~$\Dtwelve$.
We will see that all of these come from the root system $R$ in Figure~\ref{fig:root_lattice},
through a maximal torus $T$ of $\Ggroup$.
In other words, we actually had no choice. A torus acting on a representation whose weight
spaces are lines \emph{forces} coordinates: its weight vectors, which are unique up to scaling.
After that forcing, the group still acts, but only through its torus and its Weyl group.
These make up the tropical reductive group $\Ggroup^{\trop}$ of~\cite[\S1.5]{GKUZ:2025}.
This is what tropicalization retains of the action of $\Ggroup$ on the three spaces.
The following theorem is our main result in this section.

\begin{theorem}\label{thm:homogeneous}
  Let $\LL \supset \KK$ be algebraically closed fields with a valuation onto $\RR$, trivial on~$\KK$.
  \begin{enumerate}[label=(\alph*)]
    \item The coordinates are forced: $x_1,\ldots,x_7$ by $T$ on $\KK^7$, and the Pl\"ucker
      coordinates $p_{ij}$ by the centralizer of $T$ in $\GL_7$ on $\wedge_2\KK^7$.
      The root system $R$ determines $\Phi$ up to rescaling $\Phi$ and the $x_i$.
      Hence $\trop(\pointspace)$, $\trop(\linespace)$ and $\trop(\flagspace)$ depend only on $R$ and on $\characteristic(\KK)$.
    \item The tropical group $\Ggroup^{\trop}$ in~\eqref{eq:tropical-reductive} acts on
      $\trop(\pointspace)$, $\trop(\linespace)$ and $\trop(\flagspace)$, induced by the
      elements of $\Ggroup(\LL)$ that act by monomial matrices.
      On $\trop(\linespace)$, the translations are the star trees in Figure~\ref{fig:lineality_stars}, and $\Dtwelve$ permutes the leaves of
      $\Ggroup$-trees as in Section~\ref{sec2}.
  \end{enumerate}
\end{theorem}

We start from the root system $R$ of $\Ggroup$. Let $\epsilon_1,\epsilon_2,\epsilon_3$ be the first three columns in (\ref{eq:lineality}), with $\epsilon_1+\epsilon_2+\epsilon_3=0$.
The six short roots in Figure~\ref{fig:root_lattice} are $\pm\epsilon_i$, and the six long roots are $\epsilon_i-\epsilon_j$ for $i\neq j$.
The roots $\alpha = \epsilon_1$ and $\beta = \epsilon_3 - \epsilon_1$ form a basis of $R$ and span the lattice $X^*=\ZZ\epsilon_1+\ZZ\epsilon_2$.
Since the root lattice of $R$ is its weight lattice, $X^*$ is the character lattice of the maximal torus $T$, see~\cite[\S II.1.6]{Jantzen:2003}.
Each character $\chi\in X^*$ is a homomorphism $T\to \KK^*$. Since $\epsilon_1,\epsilon_2$ is a basis of $X^*$, the map $t\mapsto(\epsilon_1(t),\epsilon_2(t),\epsilon_3(t))$ identifies $T$ with the triples $(t_1,t_2,t_3)\in(\KK^*)^3$ with $t_1t_2t_3=1$.
The Weyl group $\Dtwelve$ acts on~$X^*$ by its reflection representation, hence on the cocharacter lattice $X_*=\Hom(X^*,\ZZ)$ and on the Lie algebra of $T$, which is the Cartan subalgebra $\cartan=X_*\otimes\KK$. We now consider tropicalization.

Let $\LL$ be as in Theorem~\ref{thm:homogeneous}, with valuation $\val$.
For instance, $\LL$ can be the field $\overline{\KK}(\!(z^{\RR})\!)$ of generalized Puiseux series; cf.\ \cite[\S2.6]{ETC}.
Since $T$ is split, the map $\nu\colon T(\LL)\to X_*\otimes\RR=\Hom(X^*,\RR)$ that sends $t$ to the homomorphism $\chi\mapsto-\val(\chi(t))$ is surjective. This identifies $\trop(T)$ with $X_*\otimes\RR$.
From these data, Gross {\em et al.}~\cite{GKUZ:2025} define the {\em tropical reductive group}
 \begin{equation}\label{eq:tropical-reductive}
  \Ggroup^{\trop}\,=\,(X_*\otimes\RR)\rtimes \Dtwelve \, = \,   \trop(T) \rtimes \Dtwelve.
\end{equation}

The smallest nontrivial representation of $\Ggroup$ common to every field $\KK$ is $\KK^7$, with basis vector $e_i$ of weight $\epsilon_i$.
Here, $\epsilon_{i+3}=-\epsilon_i$ and $\epsilon_7=0$, and $T$ acts by $\diag(t_1,t_2,t_3,t_1^{-1},t_2^{-1},t_3^{-1},1)$.
Thus $x_i(t\cdot v)=\epsilon_i(t)\,x_i(v)$, which is the grading on $\KK[x_1,\ldots,x_7]$ in Figure~\ref{fig:root_lattice}.
The weights of $\KK^7$, the six short roots and $0$, have multiplicity one, and their labels are unique up to $\Dtwelve$.
Hence the coordinates of $\KK^7$ are forced.
Also, the centralizer of $T$ in $\GL_7$ is the diagonal torus $D$.
In this manner, we arrive at our point of departure in Section~\ref{sec1}: the space $\KK^7$ with its standard basis, which we implicitly identified with the imaginary part of the octonions. We did this by introducing a $T$-invariant form $\Phi\in\wedge_3\KK^7$.
Since $R$ has rank two, we want $T$ to be a maximal torus of the stabilizer of $\Phi$, and this leaves no choice.

\begin{proposition}\label{prop:trilinear}
  The $T$-invariant subspace of $\wedge_3\KK^7$ is
  spanned by $e_{123}$, $e_{456}$, $e_{147}$, $e_{257}$, $e_{367}$.
  Its elements with all five coefficients nonzero have $T$ as a maximal torus of their stabilizers and form one orbit $\KK^*\cdot D\cdot\Phi$.
  The stabilizer of $\Phi$ in $D$ is $T\times\mu_3$, where $\mu_3=\{\zeta I:\zeta^3=1\}$.
\end{proposition}

\begin{proof}
  The vectors $e_{ijk}$ form a basis of weight vectors for $\wedge_3\KK^7$, with weights $\epsilon_i+\epsilon_j+\epsilon_k$.
  The weight is zero, so that $e_{ijk}$ is $T$-invariant exactly for the five listed triples.
  A torus of the stabilizer that contains $T$ commutes with $T$, so it lies in  the $7$-dimensional torus  $D$.
  
  The five characters by which $D$ scales the coefficients $c_{ijk}$ are linearly independent, so the stabilizer in $D$ of
  a $T$-invariant form with $k$ nonzero $c_{ijk}$ has dimension $7-k$, which is $2$ only for $k=5$.
  The matrix $s=\diag(s_1,\ldots,s_7)\in D$ multiplies the coefficients
  $c_{ijk}$ by $s_1s_2s_3$, $s_4s_5s_6$ and $s_is_{i+3}s_7$, for $i=1,2,3$.
  For $r(c)=c_{147}c_{257}c_{367}/(c_{123}c_{456})$, we have $r(sc)=r(c)s_7^3$.
  Conversely, $D$ maps $c$ to every $c'$ such that $r(c')/r(c)$ is a cube.
  For any $\lambda\in\KK^*$ we have $r(\lambda c)=\lambda r(c)$, which proves the orbit statement.
  If $s$ fixes the form $\Phi$ in (\ref{eq:trilinear}), then $s_7^3=1$, and $\zeta=s_7$
  satisfies $\zeta^{-1}s\in T$. Conversely, $T$ and $\mu_3$ fix~$\Phi$.
\end{proof}

Proposition~\ref{prop:trilinear} holds over every field $\KK$. In characteristic $3$, the group $\mu_3$ is trivial, so the stabilizer of $\Phi$ in $D$ is $T$.

\begin{proof}[Proof of Theorem~\ref{thm:homogeneous}(a)]
We saw that $T$ forces $x_1,\ldots,x_7$ and has centralizer $D$ in $\GL_7$.
The $21$ weights of $D$ on $\wedge_2\KK^7$ are distinct, so $D$ forces the Pl\"ucker coordinates $p_{ij}=x_i\wedge x_j$, and $\flagspace\subset\pointspace\times\linespace$ inherits both.
Our varieties are built from $\Phi$ alone:
$\pointspace$ is cut out by the $6\times6$ Pfaffians of~\eqref{eq:Phi2}, and $\mathcal L=\PP(E)$, where $E$ is the kernel of the contraction
$\wedge_2\KK^7\to(\KK^7)^\vee$, $\omega\mapsto\Phi(\omega,\cdot)$, whose coordinates are the linear forms in~\eqref{eq:7linear}.
Hence $\pointspace$, $\linespace=\grassmannian(2,7)\cap\mathcal L$ and $\flagspace$ are $\Ggroup$-stable, and replacing $\Phi$ by $\lambda\,s\cdot\Phi$ moves them by $s$.
By Proposition~\ref{prop:trilinear}, any other choice of $\Phi$ and of the scalings of the $x_i$ is of this form, with $\lambda\in\KK^*$ and $s\in D(\KK)$.
So is any other identification of $R$ with Figure~\ref{fig:root_lattice}: it relabels the $x_i$ by some $w\in\Dtwelve$, which permutes the five $e_{ijk}$ of weight zero up to sign.
Either way, the tropicalizations are translated by the valuation of $s$, which vanishes as $\val$ is trivial on $\KK$.
Finally, only $\characteristic(\KK)$ matters: our varieties are defined over the prime field, and initial ideals commute with field extension.
\end{proof}

\begin{remark}
  Fix the basis 
   $\alpha = \epsilon_1$ and $\beta = \epsilon_3 - \epsilon_1$
  of $R$. The highest weights of $\KK^7$ and $\wedge_2\KK^7$ are
  $2\alpha+\beta=\epsilon_5$ and $3\alpha+2\beta=\epsilon_3+\epsilon_5$, with weight lines $\KK e_5$ and $\KK\,e_3\wedge e_5$. Their stabilizers are the parabolic subgroups $P_1$ and $P_2$ of $\Ggroup$. Hence $\pointspace=\Ggroup/P_1$ and $\linespace=\Ggroup/P_2$ are the orbits of these two lines, and $\flagspace=\Ggroup/B$ is the orbit of the flag $([e_5],[e_3\wedge e_5])$.
\end{remark}

Part (b) of Theorem~\ref{thm:homogeneous}
rests on a lemma that holds for any split reductive group $G$ over~$\KK$, with maximal torus $T$, character lattice $X^*$ and Weyl group $W=N(T)/T$. Here $N(T)$ is the normalizer of $T$ in $G$.
We follow~\cite[\S II.1.1--1.4]{Jantzen:2003}.
Since $G$ is split, every~$w\in W$ lifts to some $n_w\in N(T)(\KK)$, unique up to $T(\KK)$, and every element of $N(T)(\LL)$ is $t\,n_w$ for unique $t\in T(\LL)$ and $w\in W$.
Conjugation by $n_w$ is the action of $w$ on $T$, which~induces actions of $W$ on $X^*$ and on $\trop(T)=\Hom(X^*,\RR)$.
As for $\Ggroup$, we have the map $\nu\colon T(\LL)\to\trop(T)$.
It is surjective, vanishes on $T(\KK)$, and satisfies $\nu(n_wtn_w^{-1})=w(\nu(t))$.
We set $G^{\trop}=\trop(T)\rtimes W$ and $\trop(t\,n_w)=(\nu(t),w)$, which does not depend on the choice of~$n_w$.

\begin{lemma}\label{lem:tropical-group}
The map $\trop\colon N(T)(\LL)\to G^{\trop}$ is a surjective group homomorphism.
Let $G\to\GL_m(\KK)$ be a representation in which $N(T)$ acts by monomial matrices, so that $T$ scales the $i$th coordinate by a character $\chi_i\in X^*$.
Then, $G^{\trop}$ acts on $\RR^m/\one$ by tropical monomial matrices: the element $(u,w)\in G^{\trop}$ permutes the coordinates as $n_w$ does and then adds $(u(\chi_1),\ldots,u(\chi_m))$.
With coordinatewise $-\val$, this square commutes:
\[
  \begin{tikzcd}[column sep=large]
    N(T)(\LL)\times (\LL^*)^m/\LL^* \arrow[r] \arrow[d, "\trop\times(-\val)"'] & (\LL^*)^m/\LL^* \arrow[d, "-\val"] \\
    G^{\trop}\times(\RR^m/\one) \arrow[r] & \RR^m/\one.
  \end{tikzcd}
\]
Hence, if a subvariety $Y$ of $\PP(\KK^m)$ is $N(T)$-stable, then $\trop(Y)$ is $G^{\trop}$-stable.
\end{lemma}
\begin{proof}
The map $\trop$ is surjective since $\nu$ is.
For $w,w' \in W$ we write $n_wn_{w'}=s\,n_{ww'}$ with $s \in T(\KK)$.
Then $(t\,n_w)(t'n_{w'})$ equals $t\,(n_wt'n_w^{-1})\,s\,n_{ww'}$, and $\trop$ maps it to
\[
  \bigl(\nu(t)+\nu(n_wt'n_w^{-1})+\nu(s),\,ww'\bigr)\,=\,\bigl(\nu(t)+w(\nu(t')),\,ww'\bigr)\,=\,(\nu(t),w)\,(\nu(t'),w').
\]
By assumption, $T$ consists of diagonal matrices.
The matrix $n_w$ has entries in $\KK$ and maps the $j$th basis vector to a multiple of the $\sigma(j)$th, for a permutation $\sigma$ of $[m]$.
For $n=t\,n_w$ and $v\in(\LL^*)^m$, the coordinate of $n\cdot v$ at $\sigma(j)$ is $\chi_{\sigma(j)}(t)\,v_j$ times a scalar in $\KK^*$, so $-\val$ maps it to $u(\chi_{\sigma(j)})-\val(v_j)$, where $u=\nu(t)$.
This is the image of $-\val(v)$ under $(u,w)=\trop(n)$.
This shows that the square commutes.
Its vertical maps are surjective and $\trop$ is a homomorphism, so our rule is an action.

For the last claim, let $g\in G^{\trop}$ and $d\in\trop(Y)$.
Since $\val$ surjects onto $\RR$, we have $d=-\val(y)$ for some $y\in Y(\LL)$ with nonzero coordinates \cite[Theorem~3.2.3]{MaclaganSturmfels:2015}, and $g=\trop(n)$ for some $n\in N(T)(\LL)$.
Then $n\cdot y\in Y(\LL)$, and $g$ maps $d$ to $-\val(n\cdot y)\in\trop(Y)$.
\end{proof}

\begin{proof}[Proof of Theorem~\ref{thm:homogeneous}(b)]
Let $N=N(T)$ be the normalizer of $T$ in $\Ggroup$.
By~\cite[(2.11) and (3.4)]{Aschbacher:1987}, the Weyl group of $T$ is $N/T\cong\Dtwelve$.
The elements of $N$ permute the weight lines of $T$, so they are monomial matrices.
Conversely, a monomial matrix in $\Ggroup$ normalizes $\Ggroup\cap D$.
This group is $T$: by Proposition~\ref{prop:trilinear}, the diagonal matrices fixing $\Phi$ form $T\times\mu_3$, and $\det(\zeta I)=\zeta^7=\zeta$.
Hence $N(\LL)$ consists of the elements of $\Ggroup(\LL)$ that act by monomial matrices.
Let $\rho=(1\,5\,3\,4\,2\,6)$ and $\tau=(2\,3)(5\,6)$, with signed permutation matrices
\[
  n_\rho\colon\ e_i\mapsto e_{\rho(i)}\ \text{for } i\leq6,\quad e_7\mapsto-e_7,
  \qquad\quad
  n_\tau\colon\ e_i\mapsto-e_{\tau(i)}\ \text{for } i\leq6,\quad e_7\mapsto e_7.
\]
Both fix $\Phi$ and have determinant $1$, so they lie in $N$.
Sending a monomial matrix to its permutation of $[7]$ is a homomorphism $N\to \Sym_7$.
Its kernel contains $T$, so its image has at most $|N/T|=12$ elements, and it contains $\rho$ and $\tau$.
Hence it identifies $N/T$ with the group $\Dtwelve\subset \Sym_7$ of Section~\ref{sec2}.
Products $n_w$ of $n_\rho$ and $n_\tau$ lift all $w\in\Dtwelve$.
They map $x_i$ to $\pm x_{w(i)}$ and $p_{ij}$ to $\pm p_{w(i)w(j)}$.
The signs matter if $\characteristic(\KK)\neq2$: the permutation matrix of $\tau$ maps $p_{23}+p_{47}$ to $p_{47}-p_{23}$,
so it does not preserve the linear space $\mathcal L$ in~\eqref{eq:7linear}.

We now apply Lemma~\ref{lem:tropical-group}
to $\pointspace\subset\PP(\KK^7)$ and $\linespace\subset\PP(\wedge_2\KK^7)$, with $G^{\trop}=\Ggroup^{\trop}$ as in~\eqref{eq:tropical-reductive}.
Its proof applies verbatim to $\flagspace$ in their product, with one square for each factor.
The group $N$ acts on both spaces by monomial matrices, and $T$ scales $x_i$ and $p_{ij}$ by $\epsilon_i$ and $\epsilon_i+\epsilon_j$.
The three varieties are $\Ggroup$-stable by the proof of~(a).
So $\Ggroup^{\trop}$ acts on $\trop(\pointspace)$, $\trop(\linespace)$ and $\trop(\flagspace)$, induced by $N(\LL)$.

Explicitly, we identify $u\in\trop(T)$ with the vector $(u(\epsilon_i))_{i}=(a,-a,0)\in\RR^7$, where $a_1+a_2+a_3=0$.
Then $(u,w)$ acts on $\trop(\pointspace)$ by $d\mapsto w(d)+u$, and on $\trop(\linespace)$ by $d\mapsto w(d)+(u_i+u_j)_{ij}$, where $w(d)_{w(i)w(j)}=d_{ij}$.
On $\trop(\flagspace)$, both happen at once.
The translation adds the star tree with edge lengths $u$.
These $u$ span the lineality space, by Lemma~\ref{lem:lineality-space}: if $d$ and $-d$ satisfy the four-point condition~\eqref{eq:4point}, then $d$ is a star tree, so the lineality of $\trop(\linespace)$ lies in that of the Dressian.
The permutation relabels the leaves by $w$, as $\Dtwelve$ does in Section~\ref{sec2}.
The same holds for the two Dressians: $N$ preserves $\grassmannian(2,7)$ and $\mathcal L$, and
it permutes the $78$ tropical polynomials defining the flag Dressian.
Finally, the $n_w$ have entries $0,\pm1$, so all of this holds over every field, including $\characteristic(\KK)=2$ of Section~\ref{sec4}.
\end{proof}

We now discuss the zero weight plane $E_0$ of the representation $E$, where $\mathcal L=\PP(E)$; cf.\ proof of Theorem~\ref{thm:homogeneous}(a).
The linear forms in~\eqref{eq:7linear} have coefficients $\pm1$ and disjoint variables, so $E$ has dimension $14$ over every field.
By Figure~\ref{fig:root_lattice}, the weights $\epsilon_i+\epsilon_j$ of the $p_{ij}$ are the long roots once, for the hexagon edges, the short roots twice, such as $\epsilon_1$ for $\{p_{17},p_{56}\}$, and zero three times, for $\{p_{14},p_{25},p_{36}\}$.
Removing the weights of $(\KK^7)^\vee$, each root is a weight of $E$ once and zero twice.
So the root system chooses all coordinates of nonzero weight: $T$ alone forces the $18$ coordinates $p_{ij}$ whose weight is a root, where the two of a short root agree up to sign.
It is blind only on $E_0$, where $D$ supplies $p_{14},p_{25},p_{36}$, with sum zero.
Thus $p_{i,i+3}\mapsto\epsilon_i$ identifies $E_0^\vee$ with $X^*\otimes\KK=\cartan^\vee$, $\Dtwelve$-equivariantly.
Each line in $X^*_\QQ=X^*\otimes\QQ$ thus gives a coordinate on $E_0$: its primitive vector, up to sign.
Let $\lineset$ be a set of lines whose coordinates span $E_0^\vee$.
With one Pl\"ucker coordinate per root, it defines a tropical variety $\trop_{\lineset}(\linespace)$.
One is tempted to choose a basis, for instance $\lineset = \{\QQ\epsilon_1,\QQ\epsilon_2\}$.
This forgets $d_{36}$, and the projection $\trop(\linespace) \rightarrow \trop_{\lineset}(\linespace)$ of~\cite[Corollary~3.2.13]{MaclaganSturmfels:2015}
keeps $227$ of the $248$ maximal cells in Theorem~\ref{thm:grassmannian} and contracts the others.
Worse, the Weyl group action is lost since no basis of $E_0^\vee$ is $\Dtwelve$-stable.
This destroys the beauty of the forest and the orchard.

The $\Dtwelve$-orbits of lines in $X^*_\QQ$ have three or six elements.
The smallest $\Dtwelve$-stable choices are the \emph{short tripod} $\lineset_{\shortroots}=\{\QQ\epsilon_1,\QQ\epsilon_2,\QQ\epsilon_3\}$, the \emph{long tripod} $\lineset_{\longroots}=\{\QQ(\epsilon_2-\epsilon_3),\QQ(\epsilon_3-\epsilon_1),\QQ(\epsilon_1-\epsilon_2)\}$, and their union.
The short tripod gives the coordinates $p_{14},p_{25},p_{36}$ forced by $D$, so $\trop_{\lineset_{\shortroots}}(\linespace)$ is $\trop(\linespace)$ without the duplicate short root coordinates.
The long tripod replaces them by $p_{14}-p_{25}$, $p_{14}-p_{36}$, $p_{25}-p_{36}$, and $\trop_{\lineset_{\longroots}}(\linespace)$
has $f$-vector $(42,211,275)$ instead of $(51,241,248)$.
The union projects onto both: it adds $-\val(p_{i,i+3}-p_{j,j+3})$, which is $\max(d_{i,i+3},d_{j,j+3})$ unless $d_{i,i+3}=d_{j,j+3}$.
Neither the long tripod nor the union gives a space of trees, because
these differences are not Pl\"ucker coordinates.

The representation $E$ is a familiar one.
If $\characteristic(\KK)\neq3$, then $E$ is the adjoint representation $\Gliealgebra$ and $E_0$ is the Cartan subalgebra~$\cartan$.
In characteristic~$0$, this is~\cite[\S22.3]{FH:1991}.
In general, the highest root $\theta$ is the highest weight of~$E$, and a vector of weight~$\theta$ generates
a subrepresentation that maps onto the simple module of highest weight~$\theta$~\cite[\S II.2.13]{Jantzen:2003}.
This module has dimension $14$~\cite[Part~I, Corollary~2.7a]{Hogeweij:1982}, so it is isomorphic to both $E$ and $\Gliealgebra$.
In characteristic~$3$, the short roots generate a $7$-dimensional ideal of $\Gliealgebra$~\cite[Part~I, Table~1]{Hogeweij:1982}, so the simple module has dimension at most~$7$, and we leave open whether $E\cong\Gliealgebra$.
Our results do not depend on this identification.
Characteristic $3$ also affects the long tripod: it does not span $E_0^\vee$, since the long roots span a sublattice of index~$3$.
In characteristic $2$, every primitive vector of $X^*$ is congruent to some $\epsilon_i$ modulo $2X^*$, so every $\lineset$ yields coordinates among $p_{14},p_{25},p_{36}$.
Hence the discrepancy between $\trop_2(\linespace)$ and $\trop(\linespace)$ in Theorem~\ref{thm:char2}~is an
intrinsic feature of the group $\Ggroup$. 
 Characteristic-dependence did not come from a choice.

\section{Hexagons}
\label{sec:hexagons}

Section~\ref{sec:homogeneous} explained how the group $\Ggroup$ determines our coordinates.
In hindsight, this is~a story about buildings~\cite{Tits:1974, VanMaldeghem:1998}.
All relevant buildings have already appeared.
Linear subspaces of $\PP^6$ over $\KK$ form the spherical building of $\GL_7(\KK)$, of type $A_6$.
Its apartments correspond to the maximal split tori of $\GL_7$, and the diagonal torus $D$ gives the apartment spanned by the seven coordinate points.
Inside $\PP^6$, the points in $\pointspace$ and the lines in $\linespace$, with incidences given by $\flagspace$, form the spherical building $\hex(\KK)$ of $\Ggroup(\KK)$.
Its apartments correspond to the maximal split tori of $\Ggroup$.
The hexagon~\eqref{eq:hexagon} is the apartment of $T$: its points are the weight lines $\KK e_1,\ldots,\KK e_6$ of $T$ in $\pointspace$, and its lines are coordinate lines, so it lies in the apartment of $D$.
Hence the coordinates forced in Theorem~\ref{thm:homogeneous}(a) amount to the choice of an apartment.

Over the field $\LL$, the valuation adds the affine buildings of $\PGL_7(\LL)$ and of $\Ggroup(\LL)$~\cite{BruhatTits:1972}.
Their standard apartments are the tropical torus $\RR^7/\one$~\cite{JoswigSturmfelsYu:2007} and $\trop(T)=X_*\otimes\RR$~\cite{GanYu:2003}, which is the lineality space of $\trop(\linespace)$.
We leave the interplay between tropicalized spherical buildings and affine buildings, for $\Ggroup$ and for other split reductive groups, to future work.
Here we single out the spherical buildings $\hex(\KK)$ and $\hex(\LL)$. We tropicalize their apartments.

The incidence graph of $\hex(\KK)$ is the bipartite graph whose nodes are the points in $\pointspace$ and the lines in $\linespace$, and whose edges are given by the flags in $\flagspace$.
It has diameter $6$ and girth $12$: it is the \emph{split Cayley generalized hexagon} over the field $\KK$;
cf.~\cite[\S5]{Aschbacher:1987} and \cite[Definition 2.4.9]{VanMaldeghem:1998}.
The apartments of $\hex(\KK)$ are the $12$-cycles in this graph, which we call \emph{hexagons}; an example is~\eqref{eq:hexagon}.
Their Zariski closure is an irreducible subvariety $\hexagonspace$ of dimension $12$ in
$\flagspace^6 \subset  \pointspace^6\times\linespace^6$.

\begin{figure}[ht]
  \centering
  \figurepdf{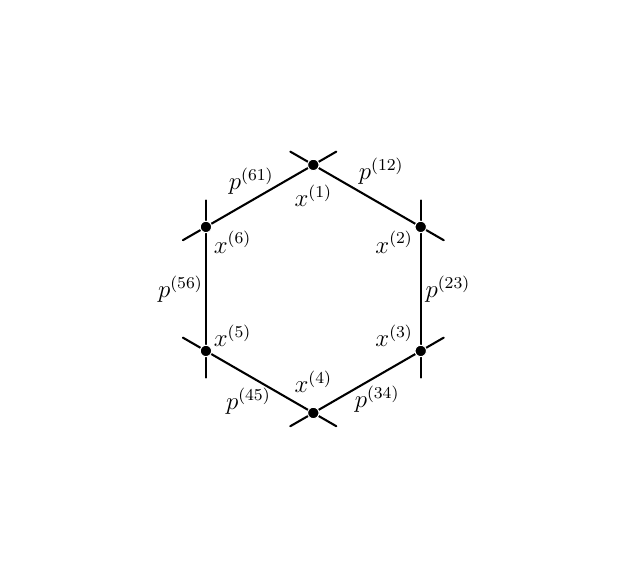}
  \caption{A hexagon for $\Ggroup$ consists of six points in $\pointspace$ and six lines in $\linespace$.}
  \label{fig:algebraic_hexagon}
\end{figure}

The prime ideal of $\hexagonspace$ lives in a polynomial ring
 with $6\times 7+ 6\times 21 = 168$ variables. We were not able
 to compute its generators. Among them are the ideal generators for
$\pointspace$ and $\linespace$, together with $12$ incidence relations given by $\flagspace$ 
as displayed in Figure~\ref{fig:algebraic_hexagon}. More precisely, all twelve pairs $(x^{(i)},p^{((i-1)i)})$ 
and $(x^{(i)},p^{(i(i+1))})$ sit in $\flagspace$, where labels are taken modulo $6$.

A hexagon is \emph{non-degenerate} when all its points and lines are distinct.
The group $\Ggroup$ acts transitively on $\pointspace$, $\linespace$, $\flagspace$, and on the non-degenerate part of~$\hexagonspace$~\cite[Theorem 5.2.9]{VanMaldeghem:1998}.
In particular, every non-degenerate hexagon is obtained as $g\cdot\eqref{eq:hexagon}$ for some $g\in\Ggroup$, with pointwise stabilizers $gTg^{-1}$~\cite[(5.8), (5.9) and \S5]{Aschbacher:1987}.
By Theorem~\ref{thm:homogeneous}(a), the torus $gTg^{-1}$ forces the coordinates moved by $g\in\Ggroup(\LL)$: in these coordinates a $\Ggroup$-stable variety $Y$ tropicalizes as $g^{-1}Y=Y$ does.
Hence, the tropical varieties of Sections~\ref{sec2}--\ref{sec7} do not depend on the apartment: every non-degenerate hexagon gives the same $\trop(\pointspace)$, $\trop(\linespace)$ and $\trop(\flagspace)$.

Consider the following projection, mapping a hexagon to its pair of \emph{opposite flags} \cite{VanMaldeghem:1998}:
\[
  \hexagonspace\,\rightarrow \,\flagspace^2,
  \,\, \left((x^{(1)},x^{(2)},\ldots, x^{(6)}),(p^{(12)},p^{(23)},\ldots,p^{(61)})\right) \,\mapsto\, \left((x^{(1)},p^{(12)}),(x^{(4)},p^{(45)})\right).
\]
In the case of the split Cayley hexagon, by Lemma \ref{lem:complete_hexagon}, this map is birational with inverse
\begin{equation}
\label{eq:sigma}
 \psi\,\colon \,\flagspace^2 \,\dashrightarrow \,\hexagonspace.
 \end{equation}
This is a parametric representation of the hexagon variety $\hexagonspace$, to be explained further below.
We here tropicalize the map $\psi$. Our main result
is the following correspondence theorem for hexagons of $\Ggroup$. 
For the sake of brevity, we assume $\characteristic(\KK)=0$ for the rest of this section.

\begin{theorem}\label{thm:completing_hexagons}
Fix fields $\KK \subset \LL$ of characteristic zero, where $\LL$ has a valuation as before.
  We give a piecewise-linear map $\trop(\psi)$ from the dominant region of $\trop(\flagspace)^2$ to $\trop(\hexagonspace)$ which 
  yields $-\val\circ\, \psi$ on points of $\flagspace(\LL)^2$ that are generic in their fiber under the valuation.
\end{theorem}

We now explain the statement.
A maximal cone in $\trop(\flagspace)$ is {\em dominant}
if its images in both $\trop(\pointspace)$ and $\trop(\linespace)$
are full-dimensional. 
Among the $2714$ maximal cones, precisely $1311$ are dominant.
Their union is the
{\em dominant region}  $\dominant \subset \trop(\flagspace)$,
and $\dominant^2$ is the dominant region of $\trop(\flagspace)^2$.
Each maximal cone of $\dominant^2$
is divided into smaller cones on which $\trop(\psi)$ is linear.
For $v$ in the relative interior of such a cone,
$-\val(\psi(z)) = \trop(\psi)(v)$ holds for
generic points $z \in \flagspace^2$ in the fiber $\trop^{-1}(v)$.
We focus on $\dominant^2$ because of Lemma \ref{lem:complete_hexagon}.

\begin{figure}[tbp]
  \centering
  \figurepdf{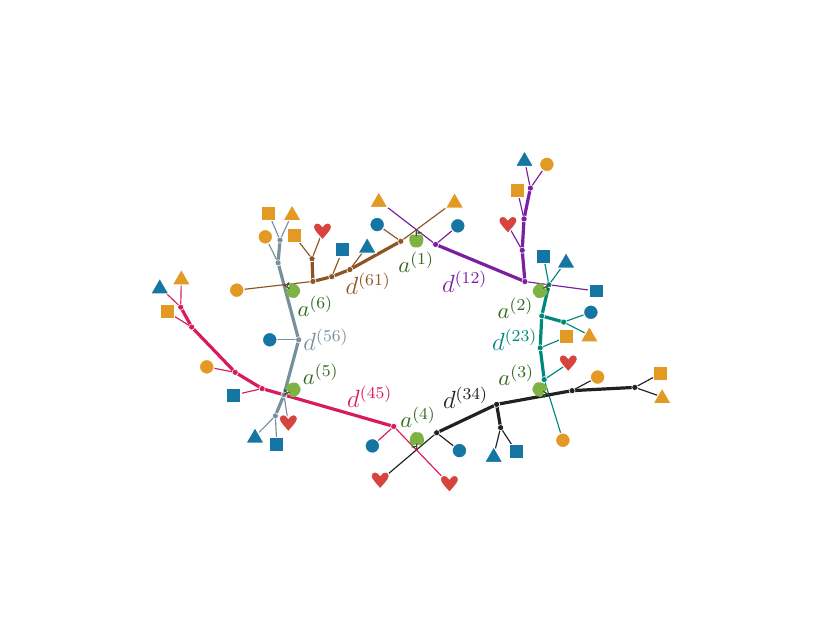}
  \caption{A tropicalized hexagon consists of six apples in $\trop(\pointspace)$ and six trees in $\trop(\linespace)$.}
  \label{fig:tropical_hexagon}
\end{figure}

The image of $\trop(\psi)$ is a polyhedral subset of $\trop(\hexagonspace)$.
Its points are \emph{tropicalized hexagons}, like that in Figure \ref{fig:tropical_hexagon}.
These have points in $\trop(\pointspace)$, lines in $\trop(\linespace)$, and incidence given by $\trop(\flagspace)$.
Thus, they are cycles of six $\Ggroup$-trees,
with six apples, such that two consecutive trees share an apple whose distance to the leaves is the same in both trees.

To prepare for the proof of Theorem \ref{thm:completing_hexagons}, we
now describe the parametrization $\psi$ in the classical case
(\ref{eq:sigma}) where domain and target are algebraic varieties.
We use the polar form
\begin{equation}
\label{eq:polarform}
\langle x,y \rangle \,=\, q(x+y)-q(x)-q(y) \,=\,
x_1y_4+x_2y_5+x_3y_6+x_4y_1+x_5y_2+x_6y_3- 2x_7y_7\, .
\end{equation}
With this, we compute the projection of a point onto a line at maximal distance \cite[\S1.5]{VanMaldeghem:1998}.

\begin{lemma}\label{lem:complete_hexagon}
Let $p^{(12)} \in \linespace$ and $x^{(4)} \in \pointspace$ be such that
the line $p^{(12)}$ is not contained in the hyperplane $ (x^{(4)})^{\perp}$.
There exist unique points $x^{(2)},x^{(3)} \in \pointspace$ and lines $p^{(23)},p^{(34)} \in \linespace$ such~that
\[p^{(12)} \ni x^{(2)} \in p^{(23)} \ni x^{(3)} \in p^{(34)} \ni x^{(4)}\, .\]
The entries of $x^{(2)},x^{(3)}$, and hence those of $p^{(23)},p^{(34)}$, depend polynomially on $p^{(12)},x^{(4)}$.
\end{lemma}

\begin{proof}
  By \cite[Remark 2.4.5]{VanMaldeghem:1998}, two points have distance $6$ if and only if they are not perpendicular 
  under (\ref{eq:polarform}).
  Since $p^{(12)} \not\subseteq (x^{(4)})^{\perp}$, there exists a point on $p^{(12)}$ that is opposite to $x^{(4)}$. This proves that the distance of $p^{(12)}$ and $x^{(4)}$ is exactly $5$.
  By \cite[Lemma 1.3.5]{VanMaldeghem:1998}, they have a unique shortest path. Let $x^{(2)},x^{(3)} \in \pointspace$ and $p^{(23)},p^{(34)} \in \linespace$ be the vertices of this path.

We compute formulas for $x^{(2)},x^{(3)}$, and hence for $p^{(23)},p^{(34)}$,
in terms of $p^{(12)}=(p_{ij}) \in \PP^{20}$ and $x^{(4)}=(x_1{:}\cdots{:}x_7) \in \PP^6$. 
The intersection $p^{(12)} \cap (x^{(4)})^{\perp}$ is the point 
$x^{(2)}$ with coordinates
\[ (x^{(2)})_i \,\,=\,\, x_4p_{1i}+x_5p_{2i}+x_6p_{3i}+x_1p_{4i}+x_2p_{5i}+x_3p_{6i}-2x_7p_{7i} 
  \quad \text{for} \,\,i=1,2,\ldots,7. \]

The \emph{cross product} with respect to our trilinear form $\Phi$ in \eqref{eq:trilinear} is defined as
  \[x \times y \;\coloneqq\; \Phi(x,y,\,\cdot\,) \;\in\; (\KK^7)^\vee\, \cong\, \KK^7 \, . \]
  The last identification is via $\langle \cdot , \cdot\rangle$.
  Since $\Ggroup$ preserves $\Phi$ and $\langle\cdot,\cdot\rangle$, the cross product is $\Ggroup$-stable.
  We set $x^{(3)} = x^{(2)} \times x^{(4)}$. The coordinates of $x^{(3)}$ are polynomials of bidegree~$(2,1)$ in $(x,p)$.
For example, $x^{(3)}_1 = 
p_{12}x_5x_7 + p_{13}x_6x_7 +  \cdots + 
 p_{67}x_1x_3 - 2p_{67}x_5x_7$ has $24$ terms. 
 
  The line spanned by two points $x,y \in \pointspace$ lies in $\linespace$ if and only if $x \times y=0$. 
  One checks that $x^{(2)}\times x^{(3)} $ $= x^{(3)}\times x^{(4)} = 0$,  that $x^{(3)}\in \pointspace$, and these points are distinct.
  Hence the lines $p^{(23)} = x^{(2)} \wedge x^{(3)}$  and $p^{(34)} = x^{(3)} \wedge x^{(4)}$ are in $\linespace$.
    Their Pl\"ucker coordinates are polynomials
    of bidegree $(3,2)$ and $(3,1)$ respectively.
 Our formulas for coordinates of $x^{(2)},x^{(3)},p^{(23)},p^{(34)}$
are not unique. They can be changed by adding 
polynomials that vanish on
 $\flagspace^2$. 
\end{proof}

Lemma \ref{lem:complete_hexagon} implies the birationality of the map $\psi$ in (\ref{eq:sigma}).
To see this, consider two flags $(x^{(1)},p^{(12)}),(x^{(4)},p^{(45)}) \in \flagspace$ whose points
$x^{(1)},x^{(4)}$ are opposite, and apply Lemma \ref{lem:complete_hexagon}  to the mixed pairs
 $(x^{(1)},p^{(45)}),(x^{(4)},p^{(12)})$. Our formulas now yield the unique hexagon completion. 

\begin{proof}[Proof of Theorem \ref{thm:completing_hexagons}]
We first split up $\psi$ into the two functions corresponding to completing the half-hexagons from mixed pairs, as in Lemma \ref{lem:complete_hexagon}. Let $I_{\pointspace\times \linespace}$ be the ideal of $\pointspace \times \linespace$,
in the polynomial ring with $28$ variables $x_k,p_{ij}$. For $v \in \trop(\pointspace)\times \trop(\linespace)$,
the initial ideal $ \initial_v(I_{\pointspace\times \linespace})$ contains no monomials.
Let $f$ be any of the coordinates of the functions constructed in
the proof of Lemma \ref{lem:complete_hexagon}.
For any $\tilde{f}\in f+I_{\pointspace\times \linespace}$ and any lift $z \in \pointspace \times \linespace$ with $-\val(z) = v$ we have $\trop(\tilde{f})(v)\geq -\val(f(z))$.
We then replace $f$ by another representative $\tilde{f} \in f+I_{\pointspace\times \linespace}$ such that
$\initial_v(\tilde{f})$ does not vanish on any
component of the initial variety $V(\initial_v(I_{\pointspace\times \linespace})) \cap (\overline{\KK}^*)^{28}$. We call such $\tilde{f}$ \emph{tropically faithful} at $v$. 
Consider a lift $z \in \pointspace \times \linespace$ and a representative
 $\alpha \in (\overline{\KK}^{*})^{28}$ of its residue class such that $\initial_v \tilde{f}(\alpha)\neq 0$.
 After rescaling the coordinates of $z$ to  valuation $0$, the equality $\trop(\tilde{f})(v)=-\val(f(z))$ holds.
The lift $z$ is \emph{generic} in the sense of the claim if it also avoids the forbidden locus $\{p \subset x^{\perp}\}$ of $\pointspace \times \linespace$ from Lemma \ref{lem:complete_hexagon}.

We performed an exhaustive computation that subdivides each of the $6\times 251 = 1506$ maximal Gr\"obner cones of $\trop(\pointspace)\times \trop(\linespace)$.
The subdivision is constructed such that each subcone $\tau$ admits a tropically faithful representative $f_{\tau}\in f+I_{\pointspace\times \linespace}$ on its relative interior. 

We now pass from half-hexagons to hexagons. 
Let a tropical flag $(a^{(1)},d^{(12)})$ be contained in the relative interior of a maximal Gr\"obner cone $\sigma$ of $\dominant$. Then any component of $V({\initial_{(a^{(1)},d^{(12)})}}I_{\flagspace})$ in the torus is an orbit of the $\dim(\sigma)$-dimensional subtorus $T_\sigma \subset (\overline{\KK}^*)^{28}$ with $\trop(T_\sigma) = \operatorname{span}(\sigma)$. 
Its images in $V(\initial_{a^{(1)}}I_{\pointspace})$ and $V(\initial_{d^{(12)}}I_{\linespace})$ are the orbits of the projections of $T_{\sigma}$, of full dimension 
by assumption, so these projections are dominant.

Thus, for a pair of general flags in the dominant region $\dominant$, a generic lift of the flags will be a generic lift of the mixed pairs $(a^{(1)},d^{(45)})$ and $(a^{(4)},d^{(12)})$. 
Hence, we can additionally require incidences $x^{(1)}\in p^{(12)}$ and $x^{(4)}\in p^{(45)}$ without changing valuations of the computed half-hexagons. 
In particular, on a dense open subset of $\dominant^2 \subset (\trop(\pointspace)\times \trop(\linespace))^2$, the valuation of a generically lifted hexagon is computed by $\trop(f_{\tau})$ applied to the mixed pair. 
These tropicalizations specify the piecewise linear map $\trop(\psi)$ promised in
Theorem \ref{thm:completing_hexagons}. 
\end{proof}

Theorem~\ref{thm:completing_hexagons} holds more generally.
The restriction to the dominant region $\dominant^2$ of $\trop(\flagspace)^2$  should not be needed.
We leave the details to a future project on tropicalizing morphisms.

\begin{example}\label{ex:tropical_hexagon}
The hexagon in Figure \ref{fig:tropical_hexagon}
was constructed by applying our piecewise-linear map
$\trop(\psi)$ to the following pair of tropical flags in $\dominant$.
As always for tree metrics, we use the max-convention, and write $a^{(i)}=-\val(x^{(i)}) \in \trop(\pointspace)$ and $d^{(ij)}=-\val(p^{(ij)})\in \trop(\linespace)$:
\[
{\footnotesize \begin{array}{ll}
  a^{(1)} = (11,12,4,15,0,14,13),
  &d^{(12)} = (9,9,12,7,11,10,10,5,8,0,7,13,0,12,11,11,7,10,10,9,9) \\
  a^{(4)} = (9,12,5,15,12,11,0),
  &d^{(45)} = (11,12,14,11,10,9,15,5,0,9,12,18,15,14,5,5,12,15,9,12,11).
\end{array}}
\]
The apples and trees in their tropical hexagon are computed by tropicalizing representatives  $\tilde{f}$ of coordinate functions $f$ of $\psi$ that are faithful in the sense of the proof  above.
We find
\[
{\footnotesize \begin{array}{ll}
   a^{(2)} = (0,3,3,6,1,5,4),&d^{(23)} = (3,3,6,1,5,4,6,9,4,8,7,7,0,6,5,5,7,6,4,3,3) \\
  a^{(3)}= (2,5,3,4,1,0,0),&d^{(34)} = (6,7,12,9,8,4,10,15,12,11,7,13,10,9,0,0,7,10,4,7,6) \\
  a^{(5)} = (0,3,3,6,3,2,0), &d^{(56)} = (2,3,6,3,2,0,6,9,6,5,3,7,4,3,3,3,3,6,0,3,2) \\
  a^{(6)} = (2,5,3,4,1,0,2),&d^{(61)} = (11,9,12,7,11,10,10,15,8,14,13,13,0,12,11,11,13,10,10,9,11).
\end{array}}\]
A visualization of the six trees $d^{(ij)}$ and their six apples $a^{(k)}$
is shown in Figure \ref{fig:tropical_hexagon}. We observe that all $\Ggroup$-trees 
in this hexagon lie in maximal cones of $\trop(\linespace)$.
Their orbit types~are:
\[
d^{(12)}:\celltype{1},\quad d^{(23)}:\celltype{24},\quad d^{(34)}:\celltype{10},\quad d^{(45)}:\celltype{15},\quad d^{(56)}:\celltype{12},\quad d^{(61)}:\celltype{23}.
\]
Thus, the hexagon appears in Figures \ref{fig:atlas_toric_cells_prime}--\ref{fig:atlas_toric_cells_nontoric},
 after relabeling the taxa via the Weyl group $\Dtwelve$.
\end{example}

\section{Curves of Genus Ten}
\label{sec10}

In this section we discuss an application of the
$\Ggroup$-Grassmannian $\linespace$ in algebraic geometry,
namely the construction of canonical curves of genus $10$
in Mukai's theory~\cite{Mukai3}.  To describe this, we write
$\Lambda_c$ for a general linear subspace
of codimension $c$ in the ambient projective space $\PP^{20}$.
By Bertini's Theorem,
the intersection $\linespace \cap \Lambda_c$ is
a smooth variety of dimension $5-c$.
The intersection, which is irreducible if $c \leq 4$,
 lives in the projective space
$\mathcal{L} \cap \Lambda_c \simeq \PP^{13-c}$.
This construction, which is due to Mukai, leads to the 
following cascade of interesting varieties:
\begin{itemize}
\item[(a)] If $c=2$, then $\linespace \cap \Lambda_c$  is a Fano threefold of genus $10$ in $\PP^{11}$. \vspace{-0.3cm}
\item[(b)] If $c=3$,  then $\linespace \cap \Lambda_c$ is a K3 surface of genus $10$ in $\PP^{10}$. \vspace{-0.3cm}
\item[(c)] If $c=4$,  then $\linespace \cap \Lambda_c$  is a canonical curve of genus $10$ in $\PP^9$. \vspace{-0.3cm}
\item[(d)] If $c=5$, then $\linespace \cap \Lambda_c$  is a self-dual configuration of $18$ points in $\PP^8$. 
\end{itemize}

Analogous constructions for $g=7,8,9$ yield birational models for the moduli space $\modulispace_g$.
For further context see~\cite{Kap1,Kap2}.
Our curves in (c) form a divisor $\divisor$ in $\modulispace_{10}$, i.e.~they
obey one constraint  in the $27$-dimensional moduli space of all curves.
The canonical divisor of the curve is the configuration in~(d).
The varieties (a)--(d) are arithmetically Gorenstein of degree $18$.
Their common Betti table records the
degrees of the higher syzygies.  It equals
\begin{equation}
\label{eq:bettitable}
\begin{bmatrix} 
    28 & 105  & 162 & 84  &    .    &    .   &    .   \\
    .     &    .   &    .   & 84  & 162 & 105 & 28  \\
\end{bmatrix}.
\end{equation}
This means that the ideal is minimally generated by $28$ quadrics, 
these quadrics satisfy $105$ independent linear syzygies, and
those in turn satisfy $162$ independent second syzygies, etc.

Here we examine these objects through the lens of tropical geometry.
The tropical moduli space $\trop(\modulispace_{10})$ is a
(stacky) polyhedral fan of dimension $27$; see~\cite{Chan:2021}.
It has $2317172$ maximal cones, one for each trivalent
graph of genus $10$, with $18$ vertices and $27$ edges.
Lower-dimensional cones of $\trop(\modulispace_{10})$ are indexed by non-trivalent graphs,
some with vertex labels that remember
contracted loops. All our graphs are undirected, connected, and they can have
loops and multiple edges. Only $41301$ of the trivalent graphs of genus $10$ are simple.

We aim to parametrize the tropical divisor $\trop(\divisor)$ in
$\trop(\modulispace_{10})$ by intersecting our  Grassmannian
$\trop(\linespace)$ with tropical linear spaces $\mathfrak{L}_4$ of codimension $4$.
These are points in $\trop(\grassmannian(17,21))$. This poses a challenge: it
is infeasible to compute this tropical Grassmannian, and we know of no way to properly sample from it.
Note that such an $\mathfrak{L}_4$ can be very large.
According to Speyer's formula~\cite{Spe},  it  can have up
to $735471$ maximal cells.
We used the following heuristics to sample from $\grassmannian(21\!-\!c,21)$. Here, the $\min$-convention was used.

\smallskip

\noindent {\bf Method 1}: {\em Shift the uniform Bergman fan.}
The Bergman fan $\mathfrak{B}_c$ consists of all $\binom{21}{c+1}$ cones
of codimension $c$ that are spanned by unit vectors in $\RR^{21}/{\bf 1}$.
We pick a random point $q$ in $ \RR^{21}/{\bf 1}$, and, 
as in \cite[Corollary 3.6.16]{MaclaganSturmfels:2015}, we take
our tropical linear space to be $\mathfrak{L}_c = q + \mathfrak{B}_c$.

\begin{figure}[H]
  \centering
  \figurepdf[scale=0.85]{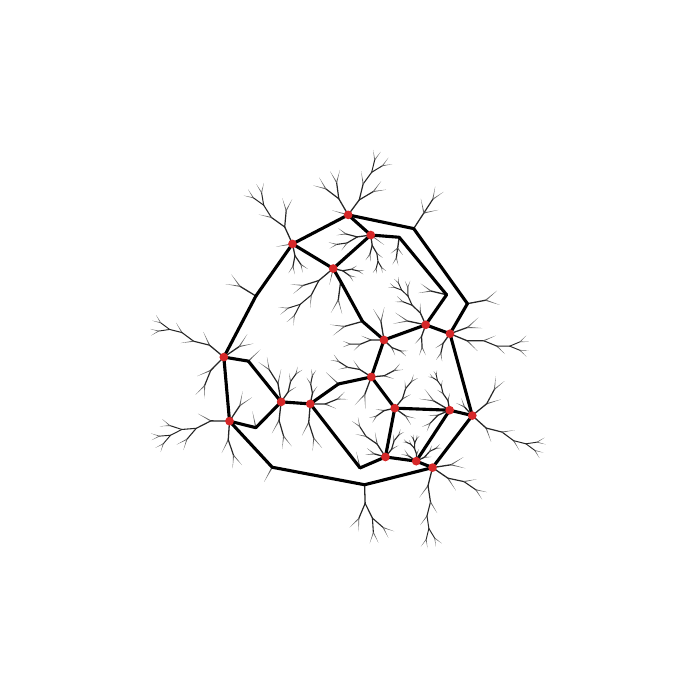}\hspace{1.5cm}%
  \figurepdf[scale=0.8]{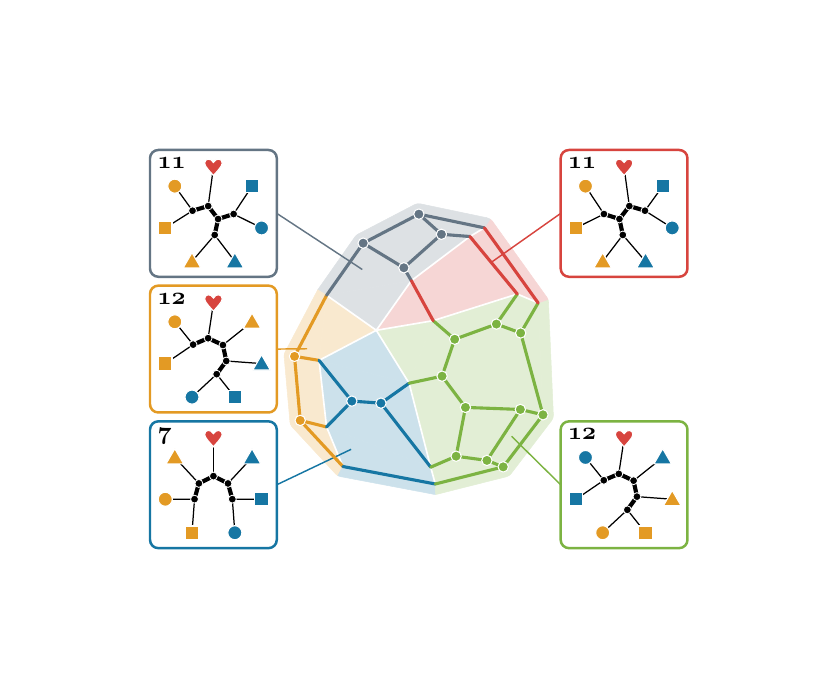}
  \caption{A planar graph obtained with Method 1 (left). Its trace on $\trop(\linespace)$ (right).}
  \label{fig:planar_curve}
\end{figure}

\smallskip

\noindent {\bf Method 2}: {\em Intersect successively with hyperplanes.}
We pick  hyperplanes  $\mathfrak{H}_i = q_i + \mathfrak{B}_{S_i}$,
where $S_i$ is a subset of the $21$ coordinates. These define
$\mathfrak{L}_c=\mathfrak{H}_1\cap\cdots\cap \mathfrak{H}_c$.
We compute
$$ \qquad  (((\trop(\linespace) \cap \mathfrak{H}_1) \cap \mathfrak{H}_2 ) \cap \mathfrak{H}_3) \cap\,\, \cdots
\quad \hbox{to get to} \quad \trop(\linespace) \,\cap \, \mathfrak{L}_c.  $$

\begin{figure}[H]
\vspace{-0.15in}
  \centering
  \figurepdf[scale=0.72]{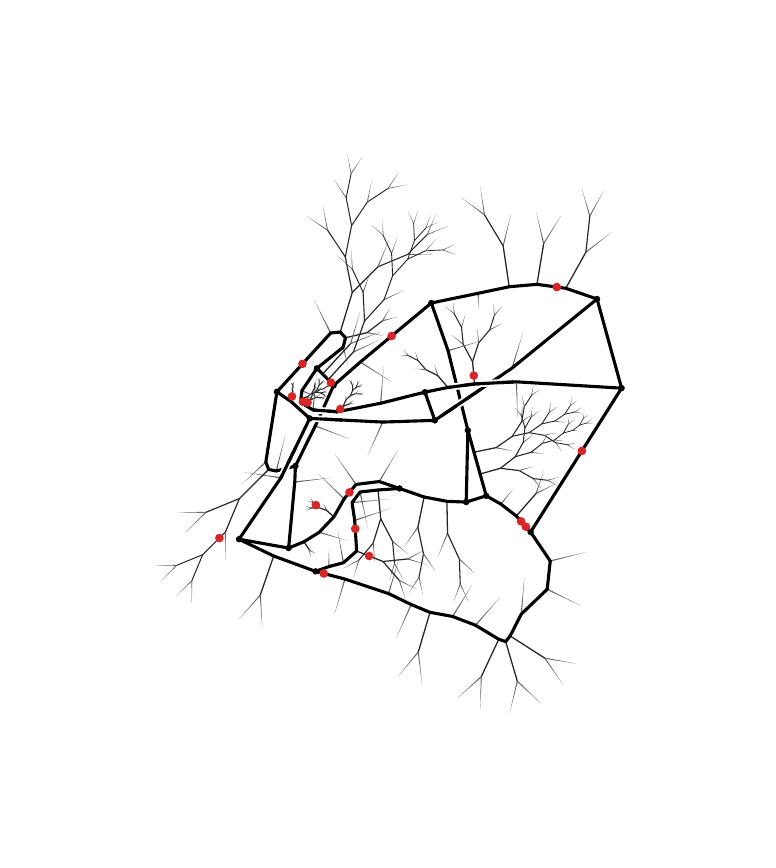}\hspace{2cm}%
  \figurepdf[scale=0.75]{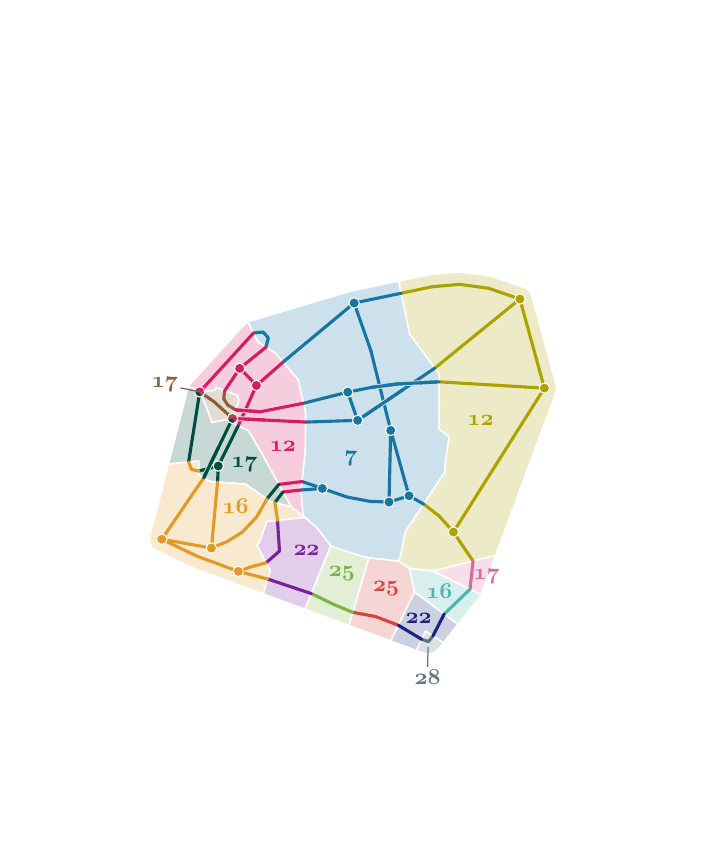}
  \caption{A nonplanar graph obtained with Method 2 (left). Its trace on $\trop(\linespace)$ (right).}
  \label{fig:nonplanar_curve}
\end{figure}

\smallskip

Using random sampling, we obtained the following result, 
illustrated in Figures~\ref{fig:planar_curve} and~\ref{fig:nonplanar_curve}.

\begin{proposition} \label{prop:mukai}
 Methods 1 and 2 can both yield
curves in maximal cones of $\trop(\modulispace_{10})$.
We obtain a sample of the $2317172$
trivalent metric graphs with $18$ vertices and $27$ edges.
\end{proposition}

\begin{proof}
We present input data for each of the two methods that leads
to a trivalent graph. For Method 1, we select the following shift vector to 
define $\mathfrak{L}_c = q + \mathfrak{B}_c$ for $c=1,2,3,4,5$:
\[ {\small \begin{array}{r@{\;\;}c@{\,}r*{5}{@{\;}c@{\;}r}}
 q \, = & - & 204958\,e_{12} & - & 1849487\,e_{13} & - & 1813412\,e_{14} & + & 1389486\,e_{15} & - & 611356\,e_{16} & - & 155086\,e_{17} \\
        & - & 128988\,e_{23} & - & 1698926\,e_{24} & + & 1863394\,e_{25} & - & 602010\,e_{26} & - & 2239681\,e_{27} & - & 1904628\,e_{34} \\
        & + & 124320\,e_{35} & - & 1025736\,e_{36} & - & 87111\,e_{37} & - & 1360439\,e_{45} & + & 705715\,e_{46} & - & 619752\,e_{47} \\
        & + & 2004970\,e_{56} & + & 1344973\,e_{57} & + & 232336\,e_{67}\rlap{.}
\end{array}}
\]
For Method 2, our four hyperplanes are defined by the following shift vectors $q_1,q_2,q_3,q_4$:
\[  {\small \begin{array}{r@{\;\;}c@{\,}r*{4}{@{\;}c@{\;}r}}
 q_1 \, = & - & 217936\,e_{13} & - & 84861\,e_{37} & + & 34008\,e_{46}\rlap{,} \\[2mm]
 q_2 \, = & & 1013667\,e_{14} & + & 132815\,e_{15} & + & 744071\,e_{25} & - & 794071\,e_{26} & + & 721948\,e_{37} \\
          & - & 60165\,e_{46} & + & 360635\,e_{47} & + & 705771\,e_{67}\rlap{,} \\[2mm]
 q_3 \, = & - & 806183\,e_{15} & - & 732664\,e_{26} & - & 107361\,e_{27} & - & 721205\,e_{34} & - & 846796\,e_{35} \\
          & - & 183135\,e_{36} & + & 542064\,e_{46} & + & 278271\,e_{47} & + & 78144\,e_{67}\rlap{,} \\[2mm]
 q_4 \, = & & 722953\,e_{12} & + & 1004831\,e_{14} & - & 757760\,e_{15} & - & 217944\,e_{16} & + & 733597\,e_{17} \\
          & + & 589026\,e_{23} & + & 213921\,e_{24} & + & 650364\,e_{26} & - & 890415\,e_{34} & - & 891447\,e_{35} \\
          & - & 636797\,e_{36} & + & 914965\,e_{37} & + & 315068\,e_{46} & - & 862815\,e_{56} & + & 692897\,e_{57}\rlap{.}
\end{array}}
\]
The divisor (d) is found by intersecting with one additional tropical hyperplane $q_5$.
For each sample, we compute the graph given by
 $\trop(\linespace) \cap \mathfrak{L}_4 $, and we check whether it is trivalent.
\end{proof}

A picture says more than a thousand words.
Our two genus $10$ curves are shown in Figures~\ref{fig:planar_curve} and~\ref{fig:nonplanar_curve}.
The left panel in each figure is the actual tropical curve. This curve has degree $18$ in dimension $15$, so it has up to $15\times 18 = 270$ unbounded rays.
These form trees that are attached to the {\em Berkovich skeleton}.
This is our metric graph, with $2$-valent vertices removed, shown here with dark edges.
Miraculously, that graph is trivalent, with $27$ edges.
 
The $18$ red points form the {\em canonical divisor} (d).
With Method 1, they are the $18$ vertices.
With Method 2, the $18$ self-dual points sit on edges of the curve.
We see this in Figure \ref{fig:nonplanar_curve}.
  
We now explain the panels on the right. The picture shows the
 Berkovich skeleton and explains how it
arises from intersecting 
$\trop(\linespace)$
with the linear space $\mathfrak{L}_4$.
First consider Figure \ref{fig:planar_curve},
where the graph is planar. Only five of the
$248$ maximal cones of $\trop(\linespace)$
meet the skeleton.
These are the cones of types
\celltype{11},
\celltype{12},
\celltype{7},
\celltype{11},
\celltype{12}, whose $\Ggroup$-trees are displayed.
The five shaded regions
indicate the intersections of these five cones with
the shifted Bergman fan $\mathfrak{L}_4$.
For instance, the blue region shows that 
the cone of type \celltype{7} intersects
$\mathfrak{L}_4$ in two pieces on the skeleton: one edge
and one tree with five edges. The two blue vertices
are the contribution of the $\Ggroup$-tree \celltype{7} to the
canonical divisor, which equals $\trop(\linespace) \cap \mathfrak{L}_5$.
Broken lines between two colored regions are $3$-valent vertices of the curve that 
become $2$-valent in the Berkovich skeleton.
They are intersections of  $\mathfrak{L}_4$ with codimension-$1$ cones in $\trop(\linespace)$.
 
The trivalent genus $10$ graph in Figure \ref{fig:nonplanar_curve} is not planar.
Here, the skeleton meets many more maximal cones of $\trop(\linespace)$.
The numbers
\celltype{17},
\celltype{12}, \celltype{16},
$\ldots$
 refer to their types of
  $\Ggroup$-trees  in Figures  \ref{fig:atlas_toric_cells_prime} and \ref{fig:atlas_toric_cells_nontoric}.
 For both approaches,
the intersection $\trop(\linespace) \cap \mathfrak{L}_3$ is the tropicalization of the K3-surface  in (b).
With Method 1, its Berkovich skeleton collapses onto a polyhedral $2$-sphere, its \emph{essential skeleton}, and hence onto the boundary of a $3$-polytope by Steinitz's Theorem.
In Figure \ref{fig:planar_curve}, our curve is the edge graph of that simple polytope.
We conclude that 
that it would be worthwhile to develop
{\em tropical Mukai theory}
for curves of genus $7,8,9,10$.

\bibliographystyle{mukai}
\bibliography{biblio}

\bigskip \bigskip
		
\footnotesize

\noindent {\bf Authors' addresses:}
		
\bigskip

\noindent

\begin{tabularx}{\linewidth}{@{}X r@{}}
  \textbf{Michael Joswig}                              & \url{joswig@math.tu-berlin.de} \\
  Technische Universit\"at Berlin, Chair of Discrete Mathematics / Geometry, Berlin \& MPI Mathematics in the Sciences, Leipzig, Germany & \\[4ex]
  \textbf{Kevin Kuehn}                                 & \url{kuehn@math.tu-berlin.de} \\
  Technische Universit\"at Berlin, Chair of Discrete Mathematics / Geometry, Germany & \\[4ex]
  \textbf{Bernd Sturmfels}                             & \url{bernd@mis.mpg.de} \\
  MPI Mathematics in the Sciences, Leipzig, Germany    & \\[2ex]
  \textbf{Mathieu Vall\'ee}                            & \url{vallee@mis.mpg.de} \\
  MPI Mathematics in the Sciences, Leipzig, Germany    &
\end{tabularx}

\bigskip

\noindent {\bf Acknowledgments:}
MJ and KK received funding by the Deutsche Forschungsgemeinschaft (DFG, German Research Foundation) under Germany's Excellence Strategy -- \enquote{The Berlin Mathematics Research Center MATH$^+$} (EXC-2046/1, EXC-2046/2, pro\-ject ID 390685689), \enquote{Symbolic Tools in Mathematics and their Application} (TRR 195, project ID 286237555), and the \enquote{Mathematical Research Data Initiative (MaRDI)} (project ID 460135501).
BS and MV were supported by
the European Research Council (ERC) through the synergy grant UNIVERSE+, 101118787.
$\!\!$ \begin{footnotesize}Views~and~opinions expressed
are however those of the authors only and do not necessarily reflect those of the European Union or the 
European
Research Council Executive Agency. Neither the European Union nor the granting authority
can be held responsible for them.
\end{footnotesize}

\bigskip

\appendix

\section{Data}\label{sec:data}
Our data are available in the Zenodo repository \href{https://doi.org/10.5281/zenodo.24104003}{\texttt{10.5281/zenodo.24104003}}; this repository comprises the $222$ files of Table~\ref{tab:data} and its folder \texttt{data}, the intermediate results of our computations, from which these files are copied and from which our code checks the statements of this article.
They were computed with \oscar~\cite{Oscar,OSCAR-book}, \polymake~\cite{GawrilowJoswig:2000} and \gfan~\cite{gfan}.
The files in the folders \texttt{00\_ideals} to \texttt{05\_curves} are written by \oscar~1.8.2 in its format \texttt{.mrdi}, and its function \texttt{load} reads them.

\begin{table}[tbp]
  \caption{Directory layout of the Zenodo repository.}
  \label{tab:data}
  \centering
  \footnotesize
  \renewcommand{\arraystretch}{0.9}
  \newcommand\qg{\char123 qq,gf2\char125}
  \begin{tabular*}{\linewidth}{@{\extracolsep{\fill}}lll@{}}
    \toprule
    file & content & reference \tabularnewline
    \midrule
    \texttt{00\_ideals/} & & \tabularnewline
    \quad\texttt{pointspace\_\qg} & the ideal of $\pointspace\subset\PP^6$: one quadric & \eqref{eq:quadric} \tabularnewline
    \quad\texttt{linespace\_\qg} & the ideal of $\linespace\subset\PP^{20}$: $7$ linear forms, $35$ quadrics & \eqref{eq:7linear}, \eqref{eq:grass1} \tabularnewline
    \quad\texttt{flagspace\_\qg} & the ideal of $\flagspace\subset\PP^6\times\PP^{20}$: $78$ generators & Section~\ref{sec2} \tabularnewline
    \texttt{01\_dressians/} & & \tabularnewline
    \quad\texttt{dressian} & the $\Ggroup$-Dressian: $51$ rays, $239$ maximal cones & Theorem~\ref{thm:dressian} \tabularnewline
    \quad\texttt{flag\_dressian} & the $\Ggroup$ flag Dressian: $141$ rays, $2627$ maximal cones & Theorem~\ref{thm:flag_dressian} \tabularnewline
    \texttt{02\_tropical\_varieties/} & & \tabularnewline
    \quad\texttt{grassmannian\_\qg} & $\trop(\linespace)$: $51$, $54$ rays; $248$, $254$ maximal cones & Theorems~\ref{thm:grassmannian}, \ref{thm:char2} \tabularnewline
    \quad\texttt{flag\_\qg} & $\trop(\flagspace)$: $141$, $156$ rays; $2714$, $2792$ maximal cones & Theorems~\ref{thm:flag}, \ref{thm:flagchar2} \tabularnewline
    \texttt{03\_no\_polytopes/} & & \tabularnewline
    \quad\texttt{grassmannian\_\qg/} & $17$, $20$ Newton--Okounkov polytopes \texttt{orbit\_NNN} of $\linespace$ & Theorem~\ref{thm:NObody} \tabularnewline
    \quad\texttt{flag\_\qg/} & $81$, $83$ Newton--Okounkov polytopes \texttt{orbit\_NNN} of $\flagspace$ & Theorem~\ref{thm:NObody} \tabularnewline
    \texttt{04\_hexagons/} & & \tabularnewline
    \quad\texttt{tropicalized\_hexagon} & the hexagon of Figure~\ref{fig:tropical_hexagon} & Example~\ref{ex:tropical_hexagon} \tabularnewline
    \texttt{05\_curves/} & & \tabularnewline
    \quad\texttt{genus\_10\_method\_1} & the curve of Figure~\ref{fig:planar_curve} & Proposition~\ref{prop:mukai} \tabularnewline
    \quad\texttt{genus\_10\_method\_2} & the curve of Figure~\ref{fig:nonplanar_curve} & Proposition~\ref{prop:mukai} \tabularnewline
    \quad\texttt{genus\_10\_method\_1\_pieces} & the polyhedra computed for Method~1 & Proposition~\ref{prop:mukai} \tabularnewline
    \texttt{grassmannian\_orbits.csv} & the orbit types of $\trop(\linespace)$ & Theorems~\ref{thm:grassmannian}, \ref{thm:char2} \tabularnewline
    \texttt{flag\_orbits.csv} & the orbit types of $\trop(\flagspace)$ & Theorems~\ref{thm:flag}, \ref{thm:flagchar2} \tabularnewline
    \texttt{dressian\_orbits.csv} & the orbit types of the $\Ggroup$-Dressian & Theorem~\ref{thm:dressian} \tabularnewline
    \texttt{flag\_dressian\_orbits.csv} & the orbit types of the $\Ggroup$ flag Dressian & Theorem~\ref{thm:flag_dressian} \tabularnewline
    \texttt{data/} & the intermediate results, as our code writes them & \tabularnewline
    \texttt{README.md} & a summary of this appendix & \tabularnewline
    \bottomrule
  \end{tabular*}
\end{table}

\paragraph{Conventions.}
The suffixes \texttt{qq} and \texttt{gf2} name the field $\QQ$ or $\FF_2$.
By Theorems~\ref{thm:char2}, \ref{thm:flagchar2} and~\ref{thm:NObody}, the fans and polytopes over $\QQ$ 
work for every characteristic other than~$2$.
The coordinates of $\PP^6$ are $x_1,\ldots,x_7$, those of $\PP^{20}$ are $p_{12},\ldots,p_{67}$ in lexicographic order, and those of $\PP^6\times\PP^{20}$ are $x_1,\ldots,x_7,p_{12},\ldots,p_{67}$.
The fans lie in $\RR^{21}$ and $\RR^{28}$.
All tropical objects in the files use the min-convention of \oscar.
A point $w$ in a file is  $-w$ in the max-convention of this paper, so the $\Ggroup$-tree of a ray $r$ is $-r$, up to lineality, and the coordinates $x_k$ give the $d_{k8}$ of the flag Dressian.
The fans are \texttt{PolyhedralFan}s with the coarsest fan structure of Sections~\ref{sec2}--\ref{sec4}.
Their rays are numbered as in \texttt{rays\_modulo\_lineality}, and their maximal cones as in \texttt{maximal\_cones}.

\paragraph{Orbit types.}
Our computation numbers the rays and maximal cones of each fan, and their orbits under $\Dtwelve$.
The four tables \texttt{*\_orbits.csv} translate these orbits into the types of this paper, one orbit per line.
The members of the orbits in \texttt{grassmannian\_orbits.csv} and \texttt{flag\_orbits.csv} are numbered as in the files of the same names in \texttt{02\_tropical\_varieties}, those in the other two tables as in \texttt{01\_dressians}.
An orbit over $\QQ$ and an orbit over $\FF_2$ share a line when they consist of the same rays or cones.
The columns are:
\begin{itemize}[nosep]
\item \texttt{kind}: \texttt{ray} or \texttt{cell}.
\item \texttt{orbit\_type}: the type \raytype{\alpha_1}, \ldots, \raytype{\epsilon} of a ray, as in Figure~\ref{fig:atlas_rays}, or the type \celltype{1}, \ldots, \celltype{35} of a maximal cone, as in Figures~\ref{fig:atlas_toric_cells_prime} and~\ref{fig:atlas_toric_cells_nontoric}, written \texttt{alpha\_1}, \ldots, \texttt{epsilon} and \texttt{1}, \ldots, \texttt{35}.
A ray or maximal cone of a flag fan has the type of the ray or maximal cone of $\trop(\linespace)$ below it, and a ray over the lineality space has type \texttt{star}.
The tetrahedra of the $\Ggroup$-Dressian, and the cones of the flag Dressian over them, have type \texttt{epsilon}; see Figure~\ref{fig:tetrahedra}.
\item \texttt{apple}, for the flag fans: the edges of the $\Ggroup$-tree of the cell type that carry the apple, as in Figure~\ref{fig:decorated_trees}, separated by \texttt{|}.
Here \texttt{pendant k} is the pendant edge of the leaf $k$, and \texttt{inner a-b} is the inner edge between the inner nodes $a,b\in\{8,\ldots,12\}$.
\item \texttt{size}: the size of the orbit.
\item \texttt{orbit}, \texttt{qq\_orbit}, \texttt{gf2\_orbit}: for a cell, the number of its orbit, which is the number \texttt{NNN} in \texttt{03\_no\_polytopes} for a toric cell type; for a ray, the smallest index of a ray in the orbit. The entry is empty if the orbit does not occur over that field.
\item \texttt{members}, \texttt{qq\_members}, \texttt{gf2\_members}: the indices of the rays or maximal cones of the orbit.
\end{itemize}
For example, the line \texttt{ray,alpha\_2,3,8,1,8;33;34,1;19;22} of \texttt{grassmannian\_orbits.csv} says that the rays of type \raytype{\alpha_2} form an orbit of size $3$.
This orbit consists of the
 rays $8$, $33$, $34$ of \texttt{grassmannian\_qq}.
 In characteristic $p=2$, it consists of the rays $1$, $19$, $22$ of \texttt{grassmannian\_gf2}.

\paragraph{Newton--Okounkov polytopes.}
File \texttt{03\_no\_polytopes/grassmannian\_qq/orbit\_NNN} is the polytope $P_C$ of Section~\ref{sec6} for a maximal cone $C$ in the orbit \texttt{NNN} of \texttt{02\_tropical\_varieties/grassmannian\_qq}, and likewise for the other three folders.
There is one file for each orbit of toric degenerations in Theorem~\ref{thm:NObody}.
The polytopes are lattice polytopes of dimension $5$ in $\RR^5$ for $\linespace$, and of dimension $6$ in $\RR^6$ for $\flagspace$, where $P_C$ is the Minkowski sum of the polytopes of the coordinates $x_1,\ldots,x_7$ and of the coordinates $p_{12},\ldots,p_{67}$.

\paragraph{Hexagon and curves.}
The file \texttt{tropicalized\_hexagon} is a named tuple with the fields \texttt{points}, the six vectors $\val(x^{(1)}),\ldots,\val(x^{(6)})$ in $\QQ^7$, and \texttt{lines}, the six vectors $\val(p^{(12)}),\val(p^{(23)}),\ldots,\val(p^{(61)})$ in $\QQ^{21}$, each shifted to the minimum $0$.
For such a vector $v$, the vector $a^{(i)}$ or $d^{(ij)}$ of Example~\ref{ex:tropical_hexagon} is $\max(v)\one-v$.
The files \texttt{genus\_10\_method\_1} and \texttt{genus\_10\_method\_2} hold the tropical curves of genus $10$ obtained with Methods~1 and~2 of Section~\ref{sec10}, as one-dimensional \texttt{PolyhedralComplex}es in $\RR^{21}$
 with one point over each point of $\RR^{21}/\one$.
They have $449$ and $447$ maximal cells, of which $239$ and $210$ are unbounded, and their bounded cells form graphs of genus $10$ on $201$ and $228$ vertices.
The first file is a named tuple with the field \texttt{complex}.
The second also has the field \texttt{multiplicities}, the multiplicities $1$, $2$, $4$ or $5$ of the maximal cells in the order of \texttt{maximal\_polyhedra}.
The file \texttt{genus\_10\_method\_1\_pieces} is the vector of the $454$ one-dimensional polyhedra computed for Method~1 from its shift $q$; they make up the $449$ maximal cells of \texttt{genus\_10\_method\_1}, five of which consist of two collinear pieces.

\section{Computation}\label{sec:computation}

\paragraph*{Public repository.} Our code is available at \href{https://github.com/dmg-lab/TropicalGeometryG2}{github.com/dmg-lab/TropicalGeometryG2}.

\paragraph{Software and hardware.}
We used \oscar~1.8.2~\cite{OSCAR-book,Oscar} in Julia~1.13, with \polymake~\cite{GawrilowJoswig:2000} inside \oscar for polyhedral geometry, \gfan~0.8~\cite{gfan} for tropical prevarieties, and, for the curves of Section~\ref{sec10}, the LP solver \texttt{HiGHS}~1.15~\cite{HuangfuHall:2018} as a floating-point filter whose answers are checked over~$\QQ$, with \texttt{SoPlex}~\cite{GleixnerSteffyWolter:2016} (via \oscar) as an exact alternative.
The fans of Sections~\ref{sec2}--\ref{sec4} were computed on a laptop with an Intel Core i7-12700H processor ($14$ cores) and $16$~GB of memory.
The curves of Section~\ref{sec10} were searched on a cluster whose nodes have two Intel Xeon Platinum 8268 processors ($48$ cores at $2.9$~GHz) and $384$~GB of memory.

\paragraph{Timings.}
Table~\ref{tab:timings} lists the wall-clock times to compute each fan from its ideal, before coarsening.
The Dressians are the prevarieties of the generators; they do not depend on the field.
\oscar~1.8.2 tropicalizes prime ideals with the trivial valuation only over $\QQ$; over $\FF_2$, we use the same method, which calls the same function of \texttt{Singular}, and it returns the fan of Theorem~\ref{thm:char2}.
Our census computes it differently: it certifies each orbit of cells of the prevariety of the generators by a Gr\"obner basis, and adds $6$ polynomials of these Gr\"obner bases to the generators whenever some cells fail; after two such rounds, the $54$ polynomials form a tropical basis.
For the flag Dressian, \oscar's \texttt{tropical\_prevariety} ran on the cluster on a node with two Intel Xeon Gold 6226R processors ($32$ cores at $2.9$~GHz) and $754$~GB of memory, after eliminating the six coordinates that the binomial linear forms in~\eqref{eq:7linear} identify with others.

\begin{table}[htbp]
  \caption{Time to compute the tropical fans of Sections~\ref{sec2}--\ref{sec4}; loading \oscar takes $14$ s more.}
  \label{tab:timings}
  \centering
  \footnotesize
  \renewcommand{\arraystretch}{0.9}
  \begin{tabular*}{\linewidth}{@{\extracolsep{\fill}}llr@{}}
    \toprule
    fan & method & time \tabularnewline
    \midrule
    $\Ggroup$-Dressian & \gfan, $1$ / $8$ threads & $1.6$ s / $0.7$ s \tabularnewline
    & \oscar, \texttt{tropical\_prevariety} & $13$ min \tabularnewline
    $\Ggroup$ flag Dressian & \gfan, $1$ / $8$ threads & $34$ s / $9.8$ s \tabularnewline
    & \oscar, \texttt{tropical\_prevariety} & stopped after $89$ min; $4.9$ h on the cluster \tabularnewline
    $\trop(\linespace)$ over $\QQ$ and $\FF_2$ & \oscar, \texttt{tropical\_variety} & $2.5$ min each \tabularnewline
    $\trop(\linespace)$ over $\FF_2$ & tropical basis, \gfan with $8$ threads & $2.6$ min \tabularnewline
    $\trop(\flagspace)$ over $\QQ$ / $\FF_2$ & Algorithm~\ref{alg:balanced_support}, \gfan with $8$ threads & $2.5$ min / $16.5$ min \tabularnewline
    \bottomrule
  \end{tabular*}
\end{table}

\paragraph{Curves of genus $10$.}
Table~\ref{tab:curve_search} summarizes the random searches behind Proposition~\ref{prop:mukai}.
Method~1 intersects $\trop(\linespace)$ with $\mathfrak{L}_4=q+\mathfrak{B}_4$ exactly, cell by cell.
The part of $q$ orthogonal to the $15$-dimensional linear span of $\trop(\linespace)$ is fixed, taken from an earlier sample, and the part in this span is random.
With uniformly random shifts $q$ instead, $19200$ samples gave $791$ curves of genus $10$, only one of them with a trivalent core.
Method~2 cuts $\trop(\linespace)$ with four tropical hyperplanes in turn, shifted by $q_1,\ldots,q_4$ with supports of random sizes.
It decides which cells meet in floating point, and recomputes the curves of genus $10$ exactly, in $7$ to $33$ minutes each.
In both methods, the random vectors have integer entries in $[-2^{20},2^{20}]$, and sample $k$ draws them with the random number generator \texttt{Xoshiro(k)} of Julia.
In Method~1, such a vector is projected onto the span, and the shift $q$ in the proof of Proposition~\ref{prop:mukai} is twice that of sample $2722$, so that it is integral.
A curve has $18$ distinct points in Table~\ref{tab:curve_search} when its intersection with $\mathfrak{L}_5$ consists of $18$ points of multiplicity~$1$.

\begin{table}[htbp]
  \caption{The random searches on the cluster for the tropical curves in Proposition~\ref{prop:mukai}.}
  \label{tab:curve_search}
  \centering
  \footnotesize
  \renewcommand{\arraystretch}{0.9}
  \begin{tabular*}{\linewidth}{@{\extracolsep{\fill}}lrr@{}}
    \toprule
    & Method~1 & Method~2 \tabularnewline
    \midrule
    time per sample & $40$ s on one core & $20$ min on two cores \tabularnewline
    total time & $260$ core-hours & $2000$ core-hours \tabularnewline
    curves of genus $10$ & $892$ of $22380$ samples & $6$ of $3023$ samples \tabularnewline
    \quad with trivalent core & $7$ & $6$ \tabularnewline
    \quad and $18$ distinct points & $7$ & $1$ \tabularnewline
    curve of the figure & sample $2722$, Figure~\ref{fig:planar_curve} & sample $100128$, Figure~\ref{fig:nonplanar_curve} \tabularnewline
    \bottomrule
  \end{tabular*}
\end{table}

\end{document}

%% file: tables/cell_ray_orbits.tex
\begingroup
\definecolor{g2CellLabel}{HTML}{000000}
\definecolor{g2RayLabel}{HTML}{000000}
\renewcommand{\arraystretch}{0.9}
\begin{tabular*}{\linewidth}{@{\extracolsep{\fill}}lllll@{}}
  \toprule
  \textcolor{g2CellLabel}{\bfseries\boldmath\oldstylenums{1}}: \textcolor{g2RayLabel}{\boldmath$\alpha_3$}, \textcolor{g2RayLabel}{\boldmath$\alpha_3$}, \textcolor{g2RayLabel}{\boldmath$\gamma_3$} & \textcolor{g2CellLabel}{\bfseries\boldmath\oldstylenums{8}}: \textcolor{g2RayLabel}{\boldmath$\alpha_3$}, \textcolor{g2RayLabel}{\boldmath$\beta_1$}, \textcolor{g2RayLabel}{\boldmath$\beta_2$} & \textcolor{g2CellLabel}{\bfseries\boldmath\oldstylenums{15}}: \textcolor{g2RayLabel}{\boldmath$\alpha_1$}, \textcolor{g2RayLabel}{\boldmath$\alpha_4$}, \textcolor{g2RayLabel}{\boldmath$\beta_1$} & \textcolor{g2CellLabel}{\bfseries\boldmath\oldstylenums{22}}: \textcolor{g2RayLabel}{\boldmath$\alpha_2$}, \textcolor{g2RayLabel}{\boldmath$\beta_1$}, \textcolor{g2RayLabel}{\boldmath$\beta_2$} & \textcolor{g2CellLabel}{\bfseries\boldmath\oldstylenums{29}}: \textcolor{g2RayLabel}{\boldmath$\alpha_1$}, \textcolor{g2RayLabel}{\boldmath$\gamma_1$}, \textcolor{g2RayLabel}{\boldmath$\gamma_2$} \tabularnewline
  \textcolor{g2CellLabel}{\bfseries\boldmath\oldstylenums{2}}: \textcolor{g2RayLabel}{\boldmath$\alpha_3$}, \textcolor{g2RayLabel}{\boldmath$\alpha_3$}, \textcolor{g2RayLabel}{\boldmath$\gamma_3$} & \textcolor{g2CellLabel}{\bfseries\boldmath\oldstylenums{9}}: \textcolor{g2RayLabel}{\boldmath$\alpha_3$}, \textcolor{g2RayLabel}{\boldmath$\alpha_4$}, \textcolor{g2RayLabel}{\boldmath$\beta_1$} & \textcolor{g2CellLabel}{\bfseries\boldmath\oldstylenums{16}}: \textcolor{g2RayLabel}{\boldmath$\alpha_1$}, \textcolor{g2RayLabel}{\boldmath$\beta_1$}, \textcolor{g2RayLabel}{\boldmath$\beta_2$} & \textcolor{g2CellLabel}{\bfseries\boldmath\oldstylenums{23}}: \textcolor{g2RayLabel}{\boldmath$\alpha_1$}, \textcolor{g2RayLabel}{\boldmath$\alpha_3$}, \textcolor{g2RayLabel}{\boldmath$\beta_1$} & \textcolor{g2CellLabel}{\bfseries\boldmath\oldstylenums{30}}: \textcolor{g2RayLabel}{\boldmath$\gamma_1$}, \textcolor{g2RayLabel}{\boldmath$\gamma_2$}, \textcolor{g2RayLabel}{\boldmath$\gamma_1$}, \textcolor{g2RayLabel}{\boldmath$\gamma_3$} \tabularnewline
  \textcolor{g2CellLabel}{\bfseries\boldmath\oldstylenums{3}}: \textcolor{g2RayLabel}{\boldmath$\alpha_3$}, \textcolor{g2RayLabel}{\boldmath$\beta_2$}, \textcolor{g2RayLabel}{\boldmath$\gamma_1$} & \textcolor{g2CellLabel}{\bfseries\boldmath\oldstylenums{10}}: \textcolor{g2RayLabel}{\boldmath$\alpha_5$}, \textcolor{g2RayLabel}{\boldmath$\beta_2$}, \textcolor{g2RayLabel}{\boldmath$\gamma_2$} & \textcolor{g2CellLabel}{\bfseries\boldmath\oldstylenums{17}}: \textcolor{g2RayLabel}{\boldmath$\alpha_1$}, \textcolor{g2RayLabel}{\boldmath$\alpha_5$}, \textcolor{g2RayLabel}{\boldmath$\beta_2$} & \textcolor{g2CellLabel}{\bfseries\boldmath\oldstylenums{24}}: \textcolor{g2RayLabel}{\boldmath$\alpha_1$}, \textcolor{g2RayLabel}{\boldmath$\alpha_3$}, \textcolor{g2RayLabel}{\boldmath$\beta_2$} & \textcolor{g2CellLabel}{\bfseries\boldmath\oldstylenums{31}}: \textcolor{g2RayLabel}{\boldmath$\alpha_1$}, \textcolor{g2RayLabel}{\boldmath$\alpha_3$}, \textcolor{g2RayLabel}{\boldmath$\gamma_1$} \tabularnewline
  \textcolor{g2CellLabel}{\bfseries\boldmath\oldstylenums{4}}: \textcolor{g2RayLabel}{\boldmath$\alpha_3$}, \textcolor{g2RayLabel}{\boldmath$\gamma_1$}, \textcolor{g2RayLabel}{\boldmath$\gamma_3$} & \textcolor{g2CellLabel}{\bfseries\boldmath\oldstylenums{11}}: \textcolor{g2RayLabel}{\boldmath$\alpha_1$}, \textcolor{g2RayLabel}{\boldmath$\alpha_4$}, \textcolor{g2RayLabel}{\boldmath$\gamma_2$} & \textcolor{g2CellLabel}{\bfseries\boldmath\oldstylenums{18}}: \textcolor{g2RayLabel}{\boldmath$\alpha_3$}, \textcolor{g2RayLabel}{\boldmath$\alpha_4$}, \textcolor{g2RayLabel}{\boldmath$\epsilon$} & \textcolor{g2CellLabel}{\bfseries\boldmath\oldstylenums{25}}: \textcolor{g2RayLabel}{\boldmath$\alpha_2$}, \textcolor{g2RayLabel}{\boldmath$\alpha_4$}, \textcolor{g2RayLabel}{\boldmath$\beta_1$} & \textcolor{g2CellLabel}{\bfseries\boldmath\oldstylenums{32}}: \textcolor{g2RayLabel}{\boldmath$\alpha_2$}, \textcolor{g2RayLabel}{\boldmath$\alpha_3$}, \textcolor{g2RayLabel}{\boldmath$\alpha_4$} \tabularnewline
  \textcolor{g2CellLabel}{\bfseries\boldmath\oldstylenums{5}}: \textcolor{g2RayLabel}{\boldmath$\beta_2$}, \textcolor{g2RayLabel}{\boldmath$\gamma_1$}, \textcolor{g2RayLabel}{\boldmath$\gamma_2$} & \textcolor{g2CellLabel}{\bfseries\boldmath\oldstylenums{12}}: \textcolor{g2RayLabel}{\boldmath$\alpha_1$}, \textcolor{g2RayLabel}{\boldmath$\alpha_5$}, \textcolor{g2RayLabel}{\boldmath$\gamma_2$} & \textcolor{g2CellLabel}{\bfseries\boldmath\oldstylenums{19}}: \textcolor{g2RayLabel}{\boldmath$\alpha_3$}, \textcolor{g2RayLabel}{\boldmath$\alpha_3$}, \textcolor{g2RayLabel}{\boldmath$\epsilon$} & \textcolor{g2CellLabel}{\bfseries\boldmath\oldstylenums{26}}: \textcolor{g2RayLabel}{\boldmath$\alpha_3$}, \textcolor{g2RayLabel}{\boldmath$\alpha_3$}, \textcolor{g2RayLabel}{\boldmath$\beta_1$} & \textcolor{g2CellLabel}{\bfseries\boldmath\oldstylenums{33}}: \textcolor{g2RayLabel}{\boldmath$\alpha_2$}, \textcolor{g2RayLabel}{\boldmath$\alpha_3$}, \textcolor{g2RayLabel}{\boldmath$\alpha_3$} \tabularnewline
  \textcolor{g2CellLabel}{\bfseries\boldmath\oldstylenums{6}}: \textcolor{g2RayLabel}{\boldmath$\alpha_2$}, \textcolor{g2RayLabel}{\boldmath$\alpha_4$}, \textcolor{g2RayLabel}{\boldmath$\gamma_2$} & \textcolor{g2CellLabel}{\bfseries\boldmath\oldstylenums{13}}: \textcolor{g2RayLabel}{\boldmath$\alpha_2$}, \textcolor{g2RayLabel}{\boldmath$\alpha_4$}, \textcolor{g2RayLabel}{\boldmath$\delta$} & \textcolor{g2CellLabel}{\bfseries\boldmath\oldstylenums{20}}: \textcolor{g2RayLabel}{\boldmath$\alpha_2$}, \textcolor{g2RayLabel}{\boldmath$\alpha_4$}, \textcolor{g2RayLabel}{\boldmath$\epsilon$} & \textcolor{g2CellLabel}{\bfseries\boldmath\oldstylenums{27}}: \textcolor{g2RayLabel}{\boldmath$\alpha_2$}, \textcolor{g2RayLabel}{\boldmath$\alpha_3$}, \textcolor{g2RayLabel}{\boldmath$\beta_1$} & \textcolor{g2CellLabel}{\bfseries\boldmath\oldstylenums{34}}: \textcolor{g2RayLabel}{\boldmath$\alpha_3$}, \textcolor{g2RayLabel}{\boldmath$\alpha_3$}, \textcolor{g2RayLabel}{\boldmath$\alpha_4$} \tabularnewline
  \textcolor{g2CellLabel}{\bfseries\boldmath\oldstylenums{7}}: \textcolor{g2RayLabel}{\boldmath$\alpha_5$}, \textcolor{g2RayLabel}{\boldmath$\alpha_5$}, \textcolor{g2RayLabel}{\boldmath$\gamma_2$} & \textcolor{g2CellLabel}{\bfseries\boldmath\oldstylenums{14}}: \textcolor{g2RayLabel}{\boldmath$\alpha_2$}, \textcolor{g2RayLabel}{\boldmath$\alpha_4$}, \textcolor{g2RayLabel}{\boldmath$\delta$} & \textcolor{g2CellLabel}{\bfseries\boldmath\oldstylenums{21}}: \textcolor{g2RayLabel}{\boldmath$\alpha_3$}, \textcolor{g2RayLabel}{\boldmath$\alpha_3$}, \textcolor{g2RayLabel}{\boldmath$\alpha_3$} & \textcolor{g2CellLabel}{\bfseries\boldmath\oldstylenums{28}}: \textcolor{g2RayLabel}{\boldmath$\alpha_2$}, \textcolor{g2RayLabel}{\boldmath$\beta_2$}, \textcolor{g2RayLabel}{\boldmath$\gamma_2$} & \textcolor{g2CellLabel}{\bfseries\boldmath\oldstylenums{35}}: \textcolor{g2RayLabel}{\boldmath$\alpha_2$}, \textcolor{g2RayLabel}{\boldmath$\alpha_3$}, \textcolor{g2RayLabel}{\boldmath$\epsilon$} \tabularnewline
  \bottomrule
\end{tabular*}%
\endgroup

%% file: tables/toric_degeneration_table.tex
\begingroup\setlength{\tabcolsep}{2pt}%
\renewcommand{\arraystretch}{1.1}%
\begin{tabular*}{\linewidth}{@{\extracolsep{\fill}}l c >{\raggedright\arraybackslash}p{9.7cm} r c@{}}
  \toprule
  Type & $f$-vector & Newton--Okounkov polytope & \multicolumn{1}{c}{$N$} & Char \tabularnewline
  \midrule
  \textcolor[HTML]{000000}{\textbf{\oldstylenums{1}}(6)} & $(6,15,20,15,6)$ & Lecture hall simplex $P_{(1,3,2,3,1)}$, Chevalley polytope & 14 & all \tabularnewline 
  \textcolor[HTML]{000000}{\textbf{\oldstylenums{2}}(6)} & $(6,15,20,15,6)$ & Simplex & 16 & all \tabularnewline 
  \textcolor[HTML]{000000}{\textbf{\oldstylenums{3}}(12)} & $(7,19,26,19,7)$ & Triple pyramid over a quadrilateral & 14 & all \tabularnewline 
  \textcolor[HTML]{000000}{\textbf{\oldstylenums{4}}(12)} & $(7,19,26,19,7)$ & Triple pyramid over a quadrilateral & 14 & all \tabularnewline 
  \textcolor[HTML]{000000}{\textbf{\oldstylenums{5}}(6)} & $(9,26,35,24,8)$ & Pyramid over the wedge at an edge of the pyramid over a quadrilateral & 14 & all \tabularnewline 
  \textcolor[HTML]{000000}{\textbf{\oldstylenums{6}}(3)} & $(9,26,35,24,8)$ & Pyramid over the wedge at an edge of the pyramid over a quadrilateral & 15 & all \tabularnewline 
  \textcolor[HTML]{000000}{\textbf{\oldstylenums{7}}(3)} & $(13,38,48,30,9)$ & Wedge at a triangular face of the wedge at an edge of the pyramid over a quadrilateral & 14 & all \tabularnewline 
  \textcolor[HTML]{000000}{\textbf{\oldstylenums{8}}(12)} & $(8,23,32,23,8)$ & Double pyramid over the multiplex $M^{3,5}$ & 15 & all \tabularnewline 
  \textcolor[HTML]{000000}{\textbf{\oldstylenums{9}}(12)} & $(8,23,32,23,8)$ & Double pyramid over the multiplex $M^{3,5}$ & 16 & all \tabularnewline 
  \textcolor[HTML]{000000}{\textbf{\oldstylenums{10}}(6)} & $(11,33,44,29,9)$ & Pyramid over the wedge at an edge of the multiplex $M^{3,5}$ & 14 & all \tabularnewline 
  \textcolor[HTML]{000000}{\textbf{\oldstylenums{11}}(6)} & $(11,33,44,29,9)$ & Pyramid over the wedge at an edge of the multiplex $M^{3,5}$ & 14 & all \tabularnewline 
  \textcolor[HTML]{000000}{\textbf{\oldstylenums{12}}(6)} & $(11,33,44,29,9)$ & Wedge at a triangular face of the pyramid over $M^{3,5}$ & 14 & all \tabularnewline 
  \textcolor[HTML]{000000}{\textbf{\oldstylenums{13}}(6)} & $(10,33,50,37,12)$ & Double pyramid over a square antiprism & 15 & all \tabularnewline 
  \textcolor[HTML]{000000}{\textbf{\oldstylenums{14}}(3)} & $(10,33,50,37,12)$ & Double pyramid over a square antiprism & 15 & all \tabularnewline 
  \textcolor[HTML]{000000}{\textbf{\oldstylenums{15}}(12)} & $(9,27,38,27,9)$ & Pyramid over the multiplex $M^{4,7}$ & 14 & all \tabularnewline 
  \textcolor[HTML]{000000}{\textbf{\oldstylenums{16}}(12)} & $(9,27,38,27,9)$ & Pyramid over a 4-polytope with ${f=(8,19,19,8)}$ & 14 & all \tabularnewline 
  \textcolor[HTML]{000000}{\textbf{\oldstylenums{17}}(12)} & $(10,30,41,28,9)$ & Pyramid over a 4-polytope with ${f=(9,21,20,8)}$ & 14 & all \tabularnewline 
  \textcolor[HTML]{000000}{\textbf{\oldstylenums{18}}(6)} & $(7,19,26,19,7)$ & Triple pyramid over a quadrilateral & 16 & $2$ \tabularnewline 
  \textcolor[HTML]{000000}{\textbf{\oldstylenums{19}}(3)} & $(7,19,26,19,7)$ & Triple pyramid over a quadrilateral & 14 & $2$ \tabularnewline 
  \textcolor[HTML]{000000}{\textbf{\oldstylenums{20}}(3)} & $(9,26,35,24,8)$ & Pyramid over the wedge at an edge of the pyramid over a quadrilateral & 15 & $2$ \tabularnewline 
  \bottomrule
\end{tabular*}%
\endgroup

%% file: tables/local_flag_complexes.tex
\begingroup
\definecolor{g2CellLabel}{HTML}{000000}
\renewcommand{\arraystretch}{0.9}
\begin{tabularx}{\linewidth}{@{}l l >{\raggedright\arraybackslash}X@{}}
  \toprule
  $f$-vector & Maximal cells & Subtypes and orbits \tabularnewline
  \midrule
  $(13,33,31,10)$ & 10 tetrahedra & \mbox{$\mathfrak A_{1}$: \textcolor{g2CellLabel}{\bfseries\boldmath\oldstylenums{21}}, \textcolor{g2CellLabel}{\bfseries\boldmath\oldstylenums{34}}}\quad \mbox{$\mathfrak A_{2}$: \textcolor{g2CellLabel}{\bfseries\boldmath\oldstylenums{32}}}\quad \mbox{$\mathfrak A_{3}$: \textcolor{g2CellLabel}{\bfseries\boldmath\oldstylenums{33}}} \tabularnewline
  $(14,36,34,11)$ & 11 tetrahedra & \mbox{$\mathfrak B_{1}$: \textcolor{g2CellLabel}{\bfseries\boldmath\oldstylenums{1}}, \textcolor{g2CellLabel}{\bfseries\boldmath\oldstylenums{2}}}\quad \mbox{$\mathfrak B_{2}$: \textcolor{g2CellLabel}{\bfseries\boldmath\oldstylenums{6}}}\quad \mbox{$\mathfrak B_{3}$: \textcolor{g2CellLabel}{\bfseries\boldmath\oldstylenums{7}}, \textcolor{g2CellLabel}{\bfseries\boldmath\oldstylenums{12}}, \textcolor{g2CellLabel}{\bfseries\boldmath\oldstylenums{17}}}\quad \mbox{$\mathfrak B_{4}$: \textcolor{g2CellLabel}{\bfseries\boldmath\oldstylenums{9}}, \textcolor{g2CellLabel}{\bfseries\boldmath\oldstylenums{26}}}\quad \mbox{$\mathfrak B_{5}$: \textcolor{g2CellLabel}{\bfseries\boldmath\oldstylenums{11}}, \textcolor{g2CellLabel}{\bfseries\boldmath\oldstylenums{24}}}\quad \mbox{$\mathfrak B_{6}$: \textcolor{g2CellLabel}{\bfseries\boldmath\oldstylenums{15}}, \textcolor{g2CellLabel}{\bfseries\boldmath\oldstylenums{23}}, \textcolor{g2CellLabel}{\bfseries\boldmath\oldstylenums{25}}, \textcolor{g2CellLabel}{\bfseries\boldmath\oldstylenums{27}}, \textcolor{g2CellLabel}{\bfseries\boldmath\oldstylenums{31}}} \tabularnewline
  $(15,38,35,11)$ & 10 tetrahedra, 1 pyramid & \mbox{$\mathfrak C_{1}$: \textcolor{g2CellLabel}{\bfseries\boldmath\oldstylenums{3}}}\quad \mbox{$\mathfrak C_{2}$: \textcolor{g2CellLabel}{\bfseries\boldmath\oldstylenums{4}}}\quad \mbox{$\mathfrak C_{3}$: \textcolor{g2CellLabel}{\bfseries\boldmath\oldstylenums{8}}}\quad \mbox{$\mathfrak C_{4}$: \textcolor{g2CellLabel}{\bfseries\boldmath\oldstylenums{10}}}\quad \mbox{$\mathfrak C_{5}$: \textcolor{g2CellLabel}{\bfseries\boldmath\oldstylenums{13}}}\quad \mbox{$\mathfrak C_{6}$: \textcolor{g2CellLabel}{\bfseries\boldmath\oldstylenums{14}}}\quad \mbox{$\mathfrak C_{7}$: \textcolor{g2CellLabel}{\bfseries\boldmath\oldstylenums{16}}, \textcolor{g2CellLabel}{\bfseries\boldmath\oldstylenums{22}}, \textcolor{g2CellLabel}{\bfseries\boldmath\oldstylenums{29}}}\quad \mbox{$\mathfrak C_{8}$: \textcolor{g2CellLabel}{\bfseries\boldmath\oldstylenums{28}}} \tabularnewline
  $(16,39,35,11)$ & 10 tetrahedra, 1 prism & \mbox{$\mathfrak D$: \textcolor{g2CellLabel}{\bfseries\boldmath\oldstylenums{5}}} \tabularnewline
  $(16,40,36,11)$ & 9 tetrahedra, 2 pyramids & \mbox{$\mathfrak E_{1}$: \textcolor{g2CellLabel}{\bfseries\boldmath\oldstylenums{18}}}\quad \mbox{$\mathfrak E_{2}$: \textcolor{g2CellLabel}{\bfseries\boldmath\oldstylenums{19}}}\quad \mbox{$\mathfrak E_{3}$: \textcolor{g2CellLabel}{\bfseries\boldmath\oldstylenums{20}}}\quad \mbox{$\mathfrak E_{4}$: \textcolor{g2CellLabel}{\bfseries\boldmath\oldstylenums{35}}} \tabularnewline
  $(19,52,45,11)$ & 7 pyramids, 4 prisms & \mbox{$\mathfrak F$: \textcolor{g2CellLabel}{\bfseries\boldmath\oldstylenums{30}}} \tabularnewline
  \bottomrule
\end{tabularx}%
\endgroup